\documentclass[a4paper,11pt,oneside]{book}

\usepackage[french,english]{babel}
\usepackage[utf8]{inputenc}
\usepackage[T1]{fontenc}

\usepackage{hyperref}
\usepackage{amsmath, amssymb, mathrsfs,amsthm}
\usepackage{graphicx,graphics,color}
\usepackage{lineno}
\usepackage{ulem}
\usepackage{amssymb}
\usepackage[linesnumbered,ruled]{algorithm2e}
\usepackage{diagbox}
\usepackage{setspace}

\usepackage[margin=27mm]{geometry}

\newtheorem{theorem}{Theorem}[chapter]
\newtheorem{proposition}{Proposition}[chapter]
\newtheorem{corollary}{Corollary}[chapter]
\newtheorem{definition}{Definition}[chapter]
\newtheorem{example}{Example}[chapter]
\newtheorem{remark}{Remark}[chapter]

\graphicspath{{./figures/}}

\begin{document}

\begin{titlepage}
  \begin{sffamily}
  \begin{center}

    \textsc{\LARGE Université Paris-Est Marne-la-Vallée}\\[2cm]

    \textsc{\LARGE Habilitation to supervise research}\\
    \textsc{\large (Habilitation à diriger des recherches)}\\[1.5cm]
    
    {\large Discipline: Applied Mathematics} \\[1cm]

    \textbf{ \huge Dynamic programming systems for modeling and control of the traffic in transportation networks \\[1.5cm]}

    {\large presented by}\\[1.5cm]
    
    {\LARGE Nadir Farhi} \\
    {\large Research associate (Chargé de recherche) at Ifsttar - Cosys / Grettia} \\[2cm]
    
  \end{center}    
    
    \noindent
    {\large Defended on December 5th 2018 in front of the Jury composed of:} \\
    
    \noindent
    \textbf{Rapporteurs:} \\
    Cécile Appert-Rolland - Directrice de recherche - Université Paris-Sud Saclay \\
    Said Mammar - Professeur - Université d'Evry Val-d'Essonnes \\
    Markos Papageorgiou - Professor - Technical University of Crete \\ [1mm]
    
    \noindent
    \textbf{Reviewers (Examinateurs):} \\
    Carlos Canudas-de-Wit - Directeur de recherche - CNRS - INRIA - Rhones-Alpes \\
    Florian De Vuyst - Professeur - Université de Technologie de Compiègne \\
    Jean-Patrick Lebacque - Ingénieur Général des Ponts, des Eaux et des Forêts - Ifsttar \\
    Pierre Rouchon - Professeur - Ecole des Mines-ParisTech \\    
    Pravin Varaiya - Professor - University of California at Berkeley \\  [1cm]

  \begin{center}
    {\large Ifsttar - Paris Marne-la-Vallée - 2018} 
  \end{center}

  \end{sffamily}
\end{titlepage}

\frontmatter

\setcounter{tocdepth}{1} 

\newpage
This page is intentionally left blank.

\tableofcontents

\chapter{A Summary in French}
\label{chap-intro-fr}
 
\section*{Introduction}
 
La modélisation et l'optimisation de la gestion des systèmes complexes restent parmi les thèmes de recherche les plus populaires.
Les systèmes complexes sont des systèmes à plusieurs composants, avec éventuellement des sous-systèmes, interagissant entre eux et avec
l'environnement extérieur. La modélisation des systèmes complexes est intrinsèquement difficile à cause de la complexité des relations et des interactions
entre différents composants et/ou sous-systèmes, pouvant inclure des dépendances, des boucles de rétroaction, des comportements émergents, des sous-systèmes 
auto-organisés, etc. La difficulté de modélisation des systèmes complexes complique la bonne compréhension de leurs fonctionnement, limite la possibilité
d'anticipation, et rend ainsi difficile l'optimisation de leur gestion.
Nous sommes concernés dans ce travail par les systèmes de transport complexes. 
Nous nous intéressons plus précisément à la modélisation et à la régulation (ou optimisation de la gestion) du trafic dans ces systèmes.

La modélisation mathématique de la dynamique d'un système complexe se fait par la détermination ou l'identification des systèmes dynamiques 
décrivant son évolution.
L'étude et l'analyse du fonctionnement du système complexe revient alors à l'étude du système dynamique le modélisant.
De même, l'optimisation de la gestion du système complexe revient à l'optimisation d'un ou de plusieurs critère(s) 
associé(s). Cette optimisation peut être statique ou dynamique (contrôle optimal temps-réel).
La modélisation mathématique de la dynamique d'un système peut être classée en se basant sur différents critères:
déterministe ou stochastique, selon que les incertitudes sont présentes et prises en compte ou non dans la modélisation;
discrète ou continue (espace et temps) selon que les variables d'état, de caommande, et/ou le temps sont discrets ou continus;
linéaire ou non linéaire dans une telle ou telle algèbre; etc.
La compréhension de la dynamique du système peut être complète, comme elle peut être partielle ou inatteignable, selon la nature 
et la complexité du système dynamique (exemples: systèmes linéaires, systèmes ergodiques, systèmes chaotiques, etc.)
Dans le cas d'un système chaotique, par exemple, où l'état est imprévisible bien que le système est déterministe, on se limite à
la recherche de l'existence de régimes stationnaires, étant donné que la compréhension des régimes transitoires est hors d'atteinte.
Nous présentons ici, aussi bien des modèles déterministes que stochastiques, bien que la plupart soient déterministes.
Cependant, nous verrons que certain modèles déterministes admettent des interprétations stochastiques intéressantes.
Pour certains de nos modèles, nous arrivons à montrer la convergence de la dynamique vers un régime stationnaire
et dériver les phases du système analytiquement. Dans d'autre cas, nous montrons que le système dynamique est 
instable. Dans des cas plus compliqués, nous nous contenterons de simuler la dynamique du système pour l'analyser.

Un point important dans la modélisation de la dynamique des systèmes est l'échelle de modélisation.
Cette dernière est déterminée par le choix des variables d'état pour un système dynamique.
Dans le cas où on s'intéresse à la dynamique fines des unités mobiles du système (positions de véhicules, leurs vitesses, etc. 
pour un système de transport), on parle de 
modélisation dynamique microscopique. Dans le cas où on s'intéresse à la dynamique de variables représentant des agrégations
d'autres variables plus fines (densités de véhicules, leurs débits, etc. pour un système de transport), on parle de modélisation macroscopique.
Des avantages et inconvénients existent pour chacune des deux échelles de modélisation. 
L'échelle est en général choisie selon les besoins de modélisation. En d'autres termes, selon le phénomène ou le comportement qui nous intéresse
à modéliser ou à reproduire, nous déterminons l'échelle de modélisation adéquate.
En général, nous utilisons l'échelle macroscopique pour modéliser la dynamique de systèmes de grande taille (de grands réseaux, etc.),
et on utilise l'échelle microscopique pour modéliser la dynamique de systèmes de taille réduite (une partie ou un axe du réseau, etc.)
Une ou des échelles intermédiaires pourraient également être considérées si besoin. On parle dans ces cas de modélisation
mésoscopique. L'objectif de cette échelle intermédiaire est de tirer bénéfice des avantages des deux échelles microscopique et macroscopique,
sans subir leurs inconvénients.

Un autre paramètre en lien direct avec l'échelle de modélisation de la dynamique des systèmes est le niveau de décision ou de gestion
envisagé pour le système.
Comme indiqué plus haut, nous pourrons modéliser un système pour comprendre et pouvoir reproduire et anticiper sa dynamique,
comme nous pourrons modéliser un système pour pouvoir agir sur-, orienter ou optimiser, voire contrôler en temps réel, sa dynamique.
Dans ce deuxième cas, le niveau auquel nous voudrons agir sur le système pourrait déterminer l'échelle de modélisation requise.
Il se peut également qu'on ait à choisir le niveau de décision selon l'échelle de modélisation possible ou pertinente, dans le cas où
les échelles de modélisation ne sont pas toutes possibles ou pas toutes pertinentes.
Trois principaux niveaux de décisions sont en général distingués pour la gestion d'un système.
Un niveau dit \textit{stratégique} qui concerne les décisions de long terme, et qui nécessite en général une échelle macroscopique de modélisation
du système. Un niveau dit \textit{tactique} qui concerne des décisions de moyen terme, et dont les stratégies pourraient aussi bien
être développées sur la base d'une modélisation macroscopique ou microscopique (voire mésoscopique) de la dynamique du système.
Un niveau dit \textit{opérationnel} qui concerne en général les décisions à prendre en temps-réel, et dont les stratégies pourraient aussi bien
être développées sur la base d'une modélisation macroscopique, microscopique ou mésoscopique.

Le contexte dans lequel sont réalisés les travaux présentés ici est celui où
la modélisation et la gestion du trafic dans les réseaux de transport et de la mobilité en général sont en pleine mutation,
suite notamment - à l'arrivée des technologies d'information et de communications (TIC) et du numérique, - à la disponibilité de données massives 
et au développement de nouvelles approches et méthodes pour leur analyse, et - à l'automatisation grandissante des véhicules.
Tous ou la majorité des composants et des sous-systèmes des systèmes de transport, ainsi que leurs interactions sont concernés
par ces développements. La mobilité des biens et des personnes est en pleine mutation suite à ces développements.
La modélisation ainsi que les méthodes de gestion de la mobilité et du trafic devraient donc être adaptées en conséquence.

Comme c'est bien connu, l'état du trafic (véhicules et passagers) sur un réseau de transport est résultat de l'interaction entre
la demande de déplacement des passagers/véhicules et l'offre de déplacement du réseau.
Les TICs, le numérique, les big data, et l'automatisation des véhicules ont des effets très importants aussi bien sur la demande que sur l'offre de mobilité.
Au niveau demande de déplacement, toutes les étapes des choix effectués par les usagers sont affectées par les nouvelles technologies,
partant de la génération des déplacements, en passant par la distribution et le choix modal, et jusqu'au choix d'itinéraires.
Au niveau de l'offre de transport, les infrastructures sont (ou seraient) de plus en plus équipées pour permettre
l'échange d'information entre les véhicules/passagers et l'infrastructure dans les deux sens. De plus, de nouvelles stratégies et algorithmes
de régulation prenant en compte toutes ces nouvelles technologies se développent actuellement.
Notons aussi que la demande de mobilité subit les effets des nouvelles technologies sur l'offre, par interaction; et vice-versa.

Les travaux que nous présentons ici concernent l'étude de l'offre de transport plutôt que de sa demande.
Certains travaux concernent la modélisation du trafic routier (autoroutier et urbain). D'autres travaux concernent le
transport collectif, et plus précisément la modélisation de la dynamique des trains sur une ligne de métro, avec prise en
compte de la demande de déplacement des passagers. Tous les modèles et les stratégies de contrôle présentées ici
sont développés dans l'esprit du contexte décrit ci-dessus.

Ce mémoire résume mes travaux de recherche depuis ma thèse de doctorat jusqu'à début 2018.
Sont exclus de ce mémoire les travaux confidentiels réalisés dans le cadre d'un projet de recherche industrielle avec la société Metrolab
(projet de 4 ans), ainsi que les travaux de recherche non encore publiés.
Quelques travaux de recherche sont des développements naturels de mes travaux antérieurs effectués durant ma thèse de doctorat
(généralisation et extension d'approches de modélisation du trafic routier).
D'autres travaux sont de nouvelles approches et modèles ainsi que de nouvelles stratégies de régulation du trafic routier et du 
trafic en transport collectif.

Ainsi, des modèles microscopiques du trafic routier développés dans ma thèse de doctorat~\cite{Far08,Far12} sont généralisés pour tenir compte
de l'anticipation dans la conduite~\cite{FHL13}. Cette généralisation est importante car l'anticipation dans la conduite devient de plus en plus
pertinente et efficace grâce aux communications entre véhicules et entre véhicules et infrastructure de transport.
D'autre part, en s'inspirant des modèles macroscopiques du trafic routier basés sur l'algèbre min-plus~\cite{Far08,FGQ11}, nous avons
développé une approche duale pour la modélisation de la dynamique des trains sur une ligne de 
métro~\cite{FNHL16, FNHL17a,FNHL17b,SFCLG17,SFLG18_acc,Far18},
permettant de décrire les phases du trafic et de comprendre sa physique.
D'autres modèles de la théorie du \textit{Calcul réseau} (Network Calculus) sont développés durant mon post-doc à
l'INRIA Rhones Alpes et à L'ENS Paris~\cite{FG10,BFG11,BFG12}.
J'ai ensuite adapté la théorie du Network Calculus pour le développement d'une nouvelle approche
(théorie des systèmes de trafic routier)
permettant une modélisation des systèmes du trafic routier sous forme de systèmes linéaires dans l'algèbre min-plus~\cite{FHL14a, FHL14}, et de
dériver ensuite des performances sur ces systèmes, en appliquant la théorie du Network Calculus.
D'autres approches de modélisation et de régulation du trafic sont développées indépendamment des travaux de thèse de doctorat.
Je cite ici le développement des travaux sur le guidage optimal et robuste effectué dans le cadre de la thèse de Mme Farida Manseur,
ainsi que les travaux sur la simulation du trafic et de stratégies de régulation du trafic urbain avec prise en compte des 
communication entre véhicules et infrastructure (travaux en collaboration avec Cyril Nguyen Van Phu du laboratoire Grettia). 

Ce mémoire est intitulé \textit{Systèmes de programmation dynamique pour la modélisation et la régulation du trafic dans 
les réseaux de transport}.
Deux parties sont distinguées dans ce mémoire: 1) les méthodes et approches basée sur l'algèbre min-plus ou max-plus, où les
dynamiques sont des systèmes de programmation dynamique déterministe; 2) les méthodes et approches dont les systèmes dynamiques sont non linéaires
mais s'interprètent comme des systèmes de programmation dynamique stochastique.
Chacune des deux parties comporte des chapitres principaux ainsi qu'un chapitre résumant d'autres contributions à moi sur le même thème de la partie concernée.

La partie 1 inclut un premier chapitre contenant une introduction et quelques rappels nécessaires; deux principaux chapitres, l'un sur la modélisation
en algèbre max-plus de la dynamique de trains sur une ligne de métro, l'autre sur l'approche calcul réseau (Network Calculus) pour la modélisation et
le calcul de bornes de performance sur les réseaux routiers; et un dernier chapitre résumant mes autres contributions sur le thème de cette partie. 
La partie 2 inclut un premier chapitre contenant une introduction et quelques rappels nécessaires; deux principaux chapitres, l'un sur la modélisation
microscopique du trafic prenant en compte l'anticipation dans la conduite, l'autre sur la modélisation de la dynamique de trains sur une ligne de métro 
avec prise en compte de la demande de déplacement des passagers; et un dernier chapitre résumant mes autres contributions sur le thème de cette partie. 
Ci-dessous, nous donnons une brève description de chacun des travaux présentés dans ce mémoire.

\newpage
\section*{Synthèse des travaux}

Nous donnons dans cette section de brèves descriptions des approches et résultats des différents chapitres de ce mémoire.

\subsection*{Partie I - Modèles basés sur la programmation dynamique déterministe}

\begin{itemize}
 \item Chapitre 1 - Ce chapitre donne les rappels nécessaires pour les deux chapitres suivants (chapitres 2 et 3).
   Ces rappels incluent les principales notions et résultats nécessaires à la compréhension des modèles et des stratégies
   de contrôle présentés aux Chapitres~2 et~3. Le chapitre est organisé en trois parties.
   \begin{itemize}
     \item La première partie présente les principales définitions et notions sur les algèbres max-plus et min-plus,
       ainsi que les principaux résultats nécessaires aux développement des modèles des deux chapitres suivants.
     \item La deuxième partie présente les principaux théorèmes sur l'existence et l'unicité de régimes stationnaires
       pour les systèmes linéaires max-plus.
     \item La troisième partie donne les rappels nécessaires sur la théorie du \textit{calcul réseau} (Network Calculus)
       incluant les principaux résultats et notions.
   \end{itemize}
 \item Chapitre 2 - Ce chapitre présente un modèle pour la dynamique des trains sur une ligne de métro.
   La dynamique des trains est décrite par un modèle à événements discrets prenant en compte des contraintes
   sur les temps de parcours inter-stations, sur les temps de stationnements, ainsi que sur les temps de séparation
   entre trains successifs. Nous montrons que le modèle s'écrit linéairement dans l'algèbre max-plus, et que
   le système dynamique admet un régime stationnaire avec un taux d'accroissement moyen unique (indépendant de
   la condition initiale). Nous dérivons le taux d'accroissement moyen analytiquement, en fonction du nombre
   de trains circulant sur la ligne, obtenons ainsi les diagrammes de phases du modèle.
   Finalement, nous interprétons ces diagrammes en terme de trafic ferroviaire, et déduisons 
   la capacité de la ligne, le nombre optimal de trains à faire circuler sur la ligne, ainsi que la dépendance de
   la fréquence moyenne asymptotique des trains des différent paramètres de la ligne (temps de parcours,
   temps de stationnement, temps de séparation des trains, etc.)
   Les principales références pour ce chapitres sont~\cite{FNHL16, FNHL17a, FNHL18}.
 \item Chapitre 3 - Ce chapitre présente une nouvelle approche de modélisation et de calcul de bornes de performance
   sur les réseaux routiers. Cette approche est basée sur l'un des modèles macroscopiques du trafic les plus connus (le modèle
   de Lighthill Whitham and Richards (LWR)~\cite{LW55,Ric56}), et sur le schéma numérique de transmission cellulaire
   (cell transmission model~\cite{Dag94, Leb96}). Nous montrons que la description du trafic suivant ce schéma numérique
   et sous l'hypothèse d'un diagramme fondamental trapézoïdal du trafic, la dynamique s'écrit linéairement 
   en algèbre des fonctions min-plus. Nous décrivons en détail le modèle du trafic sur une section de route, et nous montrons
   qu'il admet une représentation sous forme d'un système linéaire min-plus à deux entrées et deux sorties.
   Nous dérivons ensuite analytiquement la réponse impulsionnelle de ce système.
   Nous proposons une approche \textit{théorie des système du trafic} qui nous permet de construire de
   grands systèmes du trafic routier à base de systèmes élémentaires.
   Pour cela, nous définissons quelques opérateurs qui nous permettent, par exemple, d'obtenir le système
   linéaire min-plus modélisant une route à plusieurs sections, à partir du modèle d'une seule section,
   en appliquant une concaténation de toutes les sections.
   Un autre opérateur permet de mettre une boucle fermée sur un système existant, et permet, par exemple,
   de modéliser une route circulaire, à partir d'un modèle d'une route ouverte (non circulaire).
   Nous montrons également que la théorie du \textit{calcul réseau} (Network Calculus) s'applique
   à ce type de systèmes, et permet de dériver, pour un système linéaire min-plus, des bornes de performance
   interprétées ici comme des bornes sur les temps de parcours, sur les densités du trafic, etc.
   Finalement, nous donnons quelques idées pour l'extension de cette approche aux réseaux routiers en deux
   dimensions (avec présence d'intersections).
   L'approche pourrait dans ce cas être appliquée aux réseaux urbains avec prise en compte des stratégies
   de gestion d'intersections, ou aux réseaux autoroutiers, avec prise en compte des stratégies de gestion
   des accès. Les principales références pour ce chapitre sont~\cite{FHL14a,FHL14}.
 \item Chapitre 4 - Ce chapitre résume quatre de mes autres contributions sur la modélisation déterministe
   du trafic.
   \begin{itemize}
     \item Le premier travail concerne une nouvelle approche du calcul réseau (Network Calculus) qui nous permet
	de prendre en compte la variation de la taille des paquets des flux de données dans le calcul de
	bornes de performance. Pour cela, nous avons défini une nouvelle notion de \textit{courbe de paquet}.
	Les références principales pour ce travail sont~\cite{FG10, BFG11,BFG12}.
     \item Le deuxième travail résume une approche semi-décentralisée pour la régulation du trafic urbain, que 
	nous avons proposée~\cite{FNHL15}. Cette approche est basée sur un modèle de régulation centralisée du trafic
	(TUC: Traffic Urban Control)~\cite{Dia99,DPA02}, et introduit une fenêtre de temps dans les cycles des feux de
	signalisation durant laquelle la régulation est décentralisée. La taille de cette fenêtre est 
	optimisée par le niveau centralisé de régulation. La principale référence pour ce travail est~\cite{FNHL15}. 
     \item Le troisième travail est réalisé dans le cadre d'une collaboration avec Mme Souaad Lahlah de l'Université
        de Bejaia (Algérie), qui est invitée à plusieurs reprises au laboratoire Grettia durant la préparation 
        de sa thèse de doctorat~\cite{LSBF18}. Il s'agit d'un nouveau protocole de sélection optimal d'itinéraire pour les VANETs.
        La particularité de ce protocole est qu'il intègre un modèle du trafic routier.
        La principale référence pour ce travail est~\cite{LSBF18}.
     \item Le quatrième travail est réalisé en collaboration avec Cyril Nguyen Van Phu du laboratoire Grettia.
        Il s'agit du développement d'un algorithme pour la régulation du trafic sur une intersection urbaine,
        où le feux de signalisation est supposé pouvoir communiquer avec les véhicules. L'algorithme a été évalué 
        avec une simulation numérique combinant simulation du trafic urbain sous le simulateur SUMO~\cite{sumo}
        avec simulation de la communication sous le simulateur Omnet++~\cite{Omnet} à l'aide de l'outil
        open source Veins~\cite{Sommer11}. Les deux simulateurs sont mis en boucle fermée, et une interface
        de contrôle a été utilisée pour l'implémentation de l'algorithme.
        La principale référence pour ce travail est~\cite{NFHL17}.
   \end{itemize}
\end{itemize}

\subsection*{Partie II - Modèles basés sur la programmation dynamique stochastique} 

\begin{itemize}
 \item Chapitre 5 - Ce chapitre donne les rappels nécessaires aux chapitres 6 et 7. 
   Ces rappels sont présentés en trois parties.
   \begin{itemize}
     \item Systèmes dynamiques non expansifs. Ces systèmes dynamiques sont définis par des applications non expansives. 
       Nous rappelons un résultat principal qui donne les conditions sous lesquelles ces systèmes admettent un
       régime stationnaire.
     \item Systèmes de programmation dynamique associés aux problèmes de contrôle optimal de chaînes de Markov.
       Nous rappelons qu'il s'agit d'un cas particulier des systèmes non expansifs, et réinterprétons les conditions
       d'existence de régimes stationnaires en terme de chaînes de Markov.
     \item Systèmes de programmation dynamique associés au jeux stochastiques ergodiques sur une chaîne de Markov.
       Nous rappelons, comme dans le cas précédant, qu'il s'agit d'un cas particulier des systèmes non expansifs, 
       et réinterprétons les conditions
       d'existence de régimes stationnaires en terme de jeux stochastiques sur chaînes de Markov.
   \end{itemize}
 \item Chapitre 6 - Ce chapitre résume l'extension du modèle du trafic microscopique \textit{linéaire par morceaux} que nous avons proposée 
   dans~\cite{Far08} pour
   prendre en compte l'anticipation dans la loi de poursuite. Ce modèle assume une loi de poursuite où chaque véhicule réagirait à un stimulus dépendant 
   de l'état du trafic sur plusieurs véhicules le précédant, au lieu de prendre en compte uniquement un seul véhicule leader.
   Nous montrons que la dynamique des véhicules s'écrit comme un système de programmation dynamique associé à un jeu stochastique ergodique sur une 
   chaîne de Markov. Nous dérivons le comportement émergeant de cette dynamique par une dérivation analytique du diagramme fondamental correspondant
   au régime stationnaire. Finalement, nous comparons la dynamique des véhicules et le comportement émergeant par rapport au cas sans anticipation.
   La référence principale pour ce travail est~\cite{FHL13}.
 \item Chapitre 7 - Ce chapitre résume l'extension du modèle max-plus pour la dynamique des trains, présenté au chapitre 2. 
   Cette extension a pour objectif de prendre en compte l'effet de la demande de déplacement des passagers dans la dynamique des trains et en particulier
   dans les temps de stationnement des trains en stations. Nous montrons d'abord qu'une dynamique de trains non contrôlée, où les temps de stationnement
   augmentent mécaniquement avec l'augmentation des arrivées de passagers sur les quais, est naturellement instable.
   Nous proposons ensuite une loi de contrôle des temps de stationnement pour stabiliser la dynamique des trains.
   Nous montrons que la dynamique obtenue s'écrit comme un système de programmation dynamique d'un problème de contrôle optimal
   d'une chaîne de Markov. Nous caractérisons ensuite les conditions sous lesquelles la dynamique admet un régime stationnaire, avec un
   taux d'accroissement unique qui s'interprète en terme de trafic comme le temps inter-véhiculaire (inter trains) asymptotique moyen 
   (qui est l'inverse de la fréquence des trains). De plus, nous dérivons par simulation numérique la fréquence asymptotique moyenne des trains comme fonction
   du nombre de trains circulant sur la ligne, et comparons les phases du trafic obtenues à celles obtenues analytiquement par le modèle
   max-plus du chapitre 2. Les références principales pour ce travail sont~\cite{FNHL16, FNHL17b, FNHL18, Far18}.
 \item Chapitre 8 - Ce chapitre résume deux de mes autres contributions sur la modélisation stochastique
   du trafic sur les réseaux de transport.
   \begin{itemize}
     \item Le premier travail est sur le guidage optimal et robuste des usagers des réseaux routiers.
       C'est le travail de thèse de Mme Farida Manseur~\cite{Man17,MFHL17a,MFHL17b,FHL14_icnaam} que j'ai encadrée au laboratoire Grettia.
       Il s'agit d'une extension d'une approche stochastique existante~\cite{Sam12,Sam14,Sam14b} pour le routage optimal dans les
       réseaux; voir aussi~\cite{Fan05,FN06a,FN06b,Nie09,Nie12}.
       Nous avons étendu cette approche pour prendre en compte la robustesse des stratégies de routage contre d'éventuelles 
       défaillances des liens du réseau. Les références principales pour ce travail sont~\cite{Man17,MFHL17a,MFHL17b}.
     \item Le deuxième travail résume une contribution proposée en collaboration avec des collègues de l'ENS Cachan, sur la modélisation
       de l'affectation des véhicules sur les différentes voies d'une route, en fonction de l'état du trafic.
       Nous présentons le modèle, quelques schémas numériques associés, et quelques résultats de simulation.
       La principale référence pour ce travail est~\cite{Far13}.
   \end{itemize}
\end{itemize}

\section*{Conclusions et perspectives}

En conclusion, les travaux présentés dans ce mémoire sont en partie en continuation avec mes travaux antérieurs de thèse de doctorat et de post-doctorat,
et en partie orientés vers de nouvelles directions et approches incluant toutes les nouveautés du domaine.
Nous pensons que l'étude des systèmes de transport et de la mobilité en général prendra encore plus d'importance dans l'avenir.
La modélisation mathématiques, la régulation et la simulation numérique du fonctionnement de ces systèmes 
sont et resterons nécessaires à la compréhension, à l'optimisation, et à l'anticipation des différents phénomènes et comportements émergents de ces systèmes.
L'approche \textit{systèmes dynamiques} adoptée ici est très efficace et utile à la compréhension de la physique et de la dynamique des systèmes étudiés.
Cependant, comme entamés dans les travaux de ce mémoire, les modèles et stratégies de contrôle devraient s'adapter aux nouvelles
donnes des TICs, des big data, du numérique et de l'automatisation voire de la robotique.
En parallèle, de nouveaux modèles et stratégies de contrôle devraient repenser ce domaine de modélisation et de gestion optimale des systèmes
de transport complexes dans ses nouvelles dimensions, indépendamment des anciens modèles et stratégies de contrôle.

En terme d'échelle de modélisation, l'échelle microscopique prendrait plus de part pour la modélisation de la mobilité
grâce à l'équipement grandissant des unités mobiles et des infrastructures de transport.
Cependant, l'échelle macroscopique restera indispensable pour la modélisation des phénomènes complexes.
L'un des exemples concrets est la modélisation macroscopique de l'affectation des véhicules sur les voies que nous avons
présentée ici au chapitre~\ref{chap-otherstoch}, voir aussi~\cite{Far13}.
En effet, les changements de voies et dépassements des véhicules sur une route sont très difficiles à modéliser en microscopique,
à cause, entre autres, de la non-linéarité des dynamiques. Le comportement que nous comprenons et que nous sommes capables de reproduire
avec le moins d'erreurs est plutôt macroscopique (les flux s'affectant sur les voies selon les vitesses moyennes sur chaque voie, cherchant un équilibre).

En terme de niveau de décision pour la gestion du trafic et de la mobilité, nous aurons de plus en plus de possibilités de considérer le niveau opérationnel
grâce à l'augmentation de la capacité d'observation de l'état du système, rendu possible par la disponibilité et la variété des capteurs.
Ainsi, des flux de données intéressantes arriveraient en temps réel, et pourraient donc être utilisés pour une gestion opérationnelle.
D'autre part, l'agrégation de toutes ces informations combinée à des modèles prédictifs continuerons toujours à être intéressants pour
une gestion à des niveau plus hauts (tactique et stratégique).
Par conséquent, nous aurons besoin d'imaginer des gestion intelligentes combinant plusieurs niveaux de décisions, et utilisant des
informations disponibles à tous les niveaux. Un exemple concret est celui de la régulation \textit{semi-décentralisée} pour le trafic 
routier urbain, que nous avons proposée au chapitre~\ref{chap-other-det}, voir aussi~\cite{FNHL15}.

\newpage
This page is intentionally left blank

\mainmatter

\chapter{Introduction}
\label{chap-intro-fr}

Modeling and optimization of the management of complex systems remain among the most popular research topics.
Complex systems are multi-component systems, possibly with subsystems, interacting with each other and with the external
environment. The modeling of complex systems is intrinsically difficult because of the complexity of the relationships 
and interactions between different components and/or subsystems, which may include dependencies, feedback loops, 
emerging behaviors, self-organized subsystems, etc. The difficulty of modeling complex systems complicates the proper
understanding of their operation, limits the possibility of anticipation, and thus makes it difficult to optimize their
management. We are concerned in this work by complex transport systems. We are particularly interested in the modeling and 
regulation of traffic in these systems.

Mathematical modeling of the dynamics of a complex system is done by determining or identifying a dynamic system describing
its evolution. The study and analysis of the operation of the complex system is then done by studying the dynamic 
system modeling it. Likewise, optimizing the management of a complex system amounts to the optimization of one or more
associated criteria. This optimization can be static or dynamic (optimal control). Mathematical modeling of the
dynamics of a system can be classified according to several criteria: deterministic or stochastic, depending on whether 
the uncertainties are present and taken into account in the proposed model; discrete or continuous (in space and in time) depending
on whether the state and/or time variables are discrete or continuous; linear or nonlinear in standard or in other algebras; etc.
The understanding of the system dynamics can be entire, or it can be partial or unachievable, according to the 
nature and to the complexity of the dynamic system (examples: linear systems, ergodic systems, chaotic systems, etc.) In 
the case of a chaotic system, for example, whose state is unpredictable, although the system is deterministic, one is 
limited to look for the existence of stationary regimes, since the understanding of transient regimes is out of 
reach. Here we present both deterministic and stochastic models, although most of them are deterministic. However, we 
will see that some of our deterministic models admit interesting stochastic interpretations. For some of our models, we 
show the convergence of the dynamics towards a stationary regime and derive analytically the phases of the system.
In other cases, we show that the considered dynamic system is unstable. In more complicated cases, we will simply
simulate the dynamics of the system to analyze them.

One important point in the dynamic system modeling is the scale of modeling. The latter is determined by the
choice of the state variables for a dynamic system. In the case where one is interested in the fine dynamics of the mobile
units of the system (vehicle positions, their speeds, etc. for a transport system), one speaks of microscopic dynamic 
modeling. In the case where one is rather interested in the dynamics of variables representing aggregations of other finer 
variables (densities of vehicles, their flow rates, etc. for a transport system), one speaks of macroscopic modeling.
Advantages and disadvantages exist for each of the two modeling scales. The scale is usually chosen according to the 
modeling needs. In other words, depending on the phenomenon or the behavior we are interested in modeling or reproducing,
we determine the appropriate modeling scale. In general, we use the macroscopic scale to model the dynamics of large 
systems (large networks, etc.), and use the microscopic scale to model the dynamics of small systems (a part or an
axis of the system, etc.) One or more intermediate scales could also be considered if necessary. In these cases
we speak of mesoscopic modeling. The objective of this intermediate scale is to benefit from the advantages of both 
microscopic and macroscopic scales, without suffering their disadvantages.

Another parameter directly related to the modeling scale of dynamic systems, is the level of decision or management envisaged
for the system. As indicated above, one may model a system to understand and be able to reproduce and anticipate
its dynamics, or may model a system to be able to act on-, optimize, or even control in real time,
its dynamics. In this second case, the level at which one likes to act on the system could determine the required modeling 
scale. It may also be necessary to choose the decision level according to the possible or relevant modeling scale, in the
case where some modeling scales are not possible or relevant. Three main levels of decision are generally 
distinguished for the management of a system. A so-called \textit{strategic} level which concerns long-term decisions,
and which usually requires a macroscopic modeling scale for the dynamics of the system. A so-called \textit{tactical} 
level which concerns medium-term decisions, and whose strategies could be developed on the basis of macroscopic or 
microscopic (or even mesoscopic) modeling of the dynamics of the system. A so-called \textit{operational} level which
generally concerns decisions to be made in real time, and whose strategies could be developed on the basis of macroscopic,
microscopic or mesoscopic modeling.

The context in which the work presented here is carried out is the one in which the modeling and management of traffic 
in transport networks and of mobility in general are submitted to important changing, thanks to - the arrival of information and communication 
technologies (ICT) and digitalization, - the availability of big data and the development of new approaches and methods
for their analysis, and - the growing automation of vehicles. All or most of the components and
subsystems of a transport system, as well as their interactions, are affected by these developments. The mobility of goods
and people is in turn affected and is also changing. The modeling as well as the methods of mobility and traffic management
should therefore be adapted accordingly.

As well known, the state of traffic (vehicles and passengers) on a transportation network is a result of the interaction
between the travel demand of passengers/vehicles and the transport supply of the network. ICTs, digitalization, big data, and vehicle 
automation have very important effects on both the transport demand and supply. At the level of travel demand, all the
stages of the choices made by the users are affected, starting from the generation of the 
displacements, via their distribution and the modal choice, and up to the route choice. At the transport supply level,
infrastructures are (or would be) increasingly equipped to allow the exchange of information between vehicles/passengers
and the infrastructure in both directions. In addition, new strategies and regulation algorithms taking into account all 
these new technologies are currently in development. It should also be noted that the demand for mobility is affected by the 
effects of new technologies on supply, by interaction; and vice versa.

The works we present here concern the study of transport supply rather than its demand. Some works concern the modeling of 
road traffic (highway and urban). Other works concern public transport, and more specifically the modeling of train dynamics
on a metro line, taking into account the passenger travel demand. All models and control strategies presented here are 
developed in the spirit of this context of development of mobility and transport systems.

This dissertation summarizes my research since my PhD until early 2018. The confidential work done in
the framework of an industrial research project with the Metrolab company (4-year project) is excluded from this thesis.
Very recent research not yet published does not appear. Some research works are natural  developments of my previous
work done during my PhD thesis (generalization and extension of road modeling approaches). 
Other works are new approaches and models as well as new control strategies for road and mass transit traffic.

Thus, microscopic models of road traffic developed in my PhD thesis~\cite{Far08,Far12} are generalized to take
into account anticipation in driving~\cite{FHL13}. This generalization is important because anticipation in driving
becomes more and more relevant and efficient thanks to communications between vehicles and between vehicles and transport
infrastructure. On the other hand, by taking inspiration from macroscopic models of road traffic based on the min-plus 
algebra~\cite{Far08,FGQ11}, we have developed a dual approach for modeling the train dynamics on a metro 
line~\cite{FNHL16, FNHL17a, FNHL17b, SFCLG17, SFLG18_acc, Far18}, to describe the traffic phases and understand
its physics. Other models of the \textit{Network Calculus} theory are developed during my 
post-doc position at INRIA Rhones-Alpes and at ENS Paris~\cite{FG10, BFG11, BFG12}. I then adapted the Network Calculus theory
for the development of a new approach (theory of road traffic systems) allowing a modeling of road traffic systems in
the form of linear systems in the min-plus algebra~\cite{FHL14a, FHL14}, and then derive performance bounds on these 
systems, applying the Network Calculus theory. Other modeling and traffic control approaches are developed independent
of my PhD thesis work. I cite here the development of works on the optimal and robust guidance carried out in the framework
of the thesis of Ms. Farida Manseur, as well as the work on simulation of the traffic and of  the control strategies
for urban traffic, with taking into account possible communications between vehicles and infrastructure (works in 
collaboration with Cyril Nguyen Van Phu from the Grettia laboratory).

This thesis is entitled \textit{Dynamic programming systems for modeling and control of the traffic in transportation
networks}. Two parts are distinguished in this dissertation: 1) methods and approaches based on min-plus or max-plus algebra, 
where the dynamics are deterministic dynamic programming systems; 2) methods and approaches whose dynamic systems are non-linear 
but are interpreted as stochastic dynamic programming systems. Each of the two parts includes a chapter of necessary reviews,
two main chapters and a chapter summarizing other works related to the concerned part.

Part~1 includes a first chapter containing an introduction and some necessary reviews; two main chapters, one on the
max-plus algebra model for the train dynamics on a metro line, the other one on the network calculus approach for 
modeling and calculating performance bounds on road networks; and a final chapter summarizing my other contributions on the topic of this part. 
Part~2 includes a first chapter containing an introduction and some necessary reviews; two main chapters, one on the 
microscopic modeling of traffic taking into account anticipation in driving, the other one on the modeling of the train dynamics
on a metro line taking into account the passenger travel demand; and a final chapter summarizing my other contributions on the topic of this part.


\newpage
This page is intentionally left blank

\part{Deterministic dynamic programming based modeling}
\label{part-deter}

We present in this part our works related to the traffic modeling and control in transportation networks,
based on deterministic dynamic-programming systems.
As well known, the dynamic programming system associated to a deterministic discrete-time optimal control problem 
is a max-plus or min-plus linear system, depending on whether the optimal control problem maximizes
rewards or minimizes payoffs, respectively. We notice that the Max-plus and the Min-plus algebras are dual. 
All the results of the Max-plus algebra hold in the Min-plus algebra, and vice-versa.
One can easily pass from maximization to minimization by $\min(a,b) = - \max(-a, -b)$.
The traffic models we present here are mainly based on the max-plus or min-plus algebra theory.
This part is organized in four chapters.

Chapter~\ref{chap-deter} (the first chapter of Part~\ref{part-deter}) introduces the modeling and control 
approaches and gives the main theoretic reviews necessary for the
development of the models in the three other chapters.
The reviews consists in the definitions of the main notions of the max-plus and min-plus algebras,
the main theorems on the existence and uniqueness of stationary regimes of the linear max-plus systems,
and a review on the network calculus theory including the main notions and results.

Chapter~\ref{chap-maxtrain} presents the max-plus algebra modeling approach we developed for the
traffic modeling of the train dynamics in metro line systems. We describe the train dynamics in a metro line
system with constraints on the train run, dwell and safe-separation times. We show that the train dynamics
satisfying those constraints are written linearly in the max-plus algebra. Moreover, we show that the resulting
dynamic system admits a stationary regime with a unique asymptotic average growth rate, interpreted as the 
asymptotic average train time-headway of the metro line. 
Furthermore, we derive from that, phase diagrams for the train dynamics giving the average train frequency 
as a function of the number of trains running on the metro line. 
We distinguish three traffic phases which we describe  in detail. 
We then present briefly two extensions of the approach: 1) to the case of a
metro line with a junction, 2) to the case where the travel demand is taken into account.
We give the main results obtained with the two extensions.
We notice here that some other extensions have been performed, and presented in Chapter~7 of Part~II.
Those extensions introduce also the travel passenger demand in the train dynamics, and 
derive its effect on the traffic phases. However, the dynamic systems considered in the extensions of Chapter~7
are no more (max-plus) linear.

Chapter~\ref{chap-roadnetcal} presents the new approach we developed entitled ``the road network calculus''.
This approach applies the theory of the network calculus to the traffic modeling and control of road networks.
It proposes a kind of system theory that permits to build traffic systems of road networks basing on
predefined elementary traffic systems, such as a road section.
Network calculus results are then applied to derive performance bounds in the road networks, such as
upper bounds on the travel time through given itineraries, or upper bounds on the car-densities, etc.
Moreover, these performance bounds are dependent on the control strategies set on the road network.
Therefore, the approach can be used to assess and compare different control strategies in term of
various performance bounds. This is related to consideration of the reliability of the performance 
of road networks.

Chapter~\ref{chap-other-det} (the last chapter of Part~\ref{part-deter}) summarizes all my other contributions on the deterministic traffic modeling and control.
This chapter includes, first, some extensions of the network calculus where we proposed what we called \textit{packet curves}
in order to improve the modeling of the variety of packet sizes in arrival and departures flows of the network calculus.
Second, a semi-decentralized urban traffic control is presented, where basing on a centralized urban traffic control
approach, we proposed to add a decentralized level of control which is able to benefit from possible 
vehicle-to-infrastructure (v2i) communications at the level of a road intersection, in relatively short
time intervals. Third, a proactive-optimal-path selection model with coordinator agents assisted routing
for vehicular ad hoc networks is presented. We proposed a model for the routing in vehicular ad hoc networks,
where we take into account the vehicular traffic models, in particular the car-density variation in the road
network. This work is done in collaboration with Mrs Souaad Lahlah from University of Bejaia (Algeria). 
Finally, an algorithm for urban traffic control with vehicle to infrastructure (equipped traffic lights) 
communications is proposed. Both the vehicular and the communication traffics are simulated in this work,
in a closed loop, and with an interface permitting us to control in real-time the traffic lights, in function
of the information made available in particular by the communication possibilities.  

\newpage
This page is intentionally left blank

\chapter[Deterministic dyn. prog. - based systems]{Introduction to and reviews on deterministic dynamic programming - based systems}
\label{chap-deter}

We introduce in this part the readers to some theoretical tools we use in our models.
The reviews are organized in three sections. First, we give the main definitions and notions of the Max-plus
algebras on scalars, on square matrices, on polynomials, and on matrices of polynomials.
These reviews will be mainly used in Chapter~\ref{chap-maxtrain} where traffic models and control approaches are proposed
for the train dynamics in metro line systems.
Second, we give main reviews of the Min-plus algebra system theory.
Those reviews are necessary for Chapter~\ref{chap-roadnetcal} where the theory of network calculus, which is based on the Min-plus algebra
system theory, is applied on the road networks. 
Third and finally, we give some reviews on the network calculus theory which we need also for Chapter~\ref{chap-roadnetcal}.

\section{Deterministic dynamic programming and Min-plus algebra}
\label{sec-ddp}

Let us consider the following deterministic optimal control problem. 
\begin{align}
  \min_{s\in\mathcal S} & \sum_{k=0}^{K-1} c_{i(k)}^{u(k)} + \phi_{i(K)} \label{ddp1} \\
         & i(k+1) = f(i(k), u(k)), \;\; \forall \; k, 0\leq k \leq K-1, \label{ddp2}
\end{align}
where $\mathcal I := \{1,2, \ldots, n\} \ni i$ denotes the state space of the system, 
$\mathcal U \ni u$ denotes the control action space of the system,
$\mathcal S$ is the set of feedback strategies $s: \mathcal I \to \mathcal U$,
$c_{i(k)}^{u(k)}$ denotes the payoff at state $i(k)$ with control action $u(k)$ for $0\leq k < K$, and
$\phi_{i(K)}$ denotes the final payoff at state $i(K)$, which is independent of the control, because we
do not have (or need) a control action at the final time.
 
We have in~(\ref{ddp2}) a dynamic system where the state of the system $i(k+1)$ at time $(k+1)$ is given as a function
of the state $i(k)$ at time $k$, and of the control action $u(k)$ taken at time $k$.
Payoffs $c_{i(k)}^{u(k)}$ are associated to the dynamics, where at every time $k, 0\leq k < K$, the payoff $c_{i(k)}^{u(k)}$
depends on the state of the system $i(k)$ at time $k$, and on the control action $u(k)$ taken at time $k$.
Moreover, a final payoff $\phi_{i(K)}$ is associated to the final state of the system.
The optimal control problem consists then in minimizing the sum of all the payoffs from time zero to time $K$,
with respect to the control strategy, and under the dynamics~(\ref{ddp2}).
 
Both sets $\mathcal I$ and $\mathcal U$ are assumed to be discrete and finite.


In order to write the dynamic programming system associated to the deterministic 
optimal control problem~(\ref{ddp1})-(\ref{ddp2}), let us define the value function $x$ as follows.
\begin{align}
  x_j(t) := & \min_{s\in\mathcal S} \sum_{k=t}^{K-1} c_{i(k)}^{u(k)} + \phi_{i(K)} \label{ddp3} \\
            & i(k+1) = f(i(k), u(k)), \;\; \forall \; k, t\leq k \leq K-1, \label{ddp4} \\
            & i(t) = j. \label{ddp5}
\end{align}

It is then easy to check that the value function $x$ defined in~(\ref{ddp3})-(\ref{ddp5}) satisfies the following.
\begin{align}
  & x_j(K) = \phi_j, \;\; \forall \; j\in\mathcal I,  \label{ddp6} \\
  & x_j(t) = \min_{u\in \mathcal U} \left( c_{i(t)}^{u(t)} + x_{f(j,u(t))}(t+1)\right), \;\; \forall \; j\in\mathcal I, \forall \; t, 0\leq t\leq K-1. \label{ddp7}  
\end{align}
System~(\ref{ddp6})-(\ref{ddp7}) is called the dynamic programming system (written in backward, i.e. $x(t)$ is given as a function of $x(t+1)$)
associated to the optimal control problem~(\ref{ddp1})-(\ref{ddp2}).


Let us now denote by $I(t)\in \mathcal I^n$ the line vector giving the state of the dynamic system at time $t$, with
$$I_p(t) := \begin{cases}
	      1 & \text{ if } i(t) = p,\\
	      0 & \text{otherwise}.
            \end{cases}$$
We also denote by $M^u, u\in\mathcal U$ the matrices defined as follows.
$$M^u_{ij}\begin{cases}
	    1 & \text{ if } j = f(i,u), \\
	    0 & \text{ otherwise}.
          \end{cases}$$
We notice here that $M^u, u\in\mathcal U$ are permutation matrices.

Then the dynamics~(\ref{ddp2}) of the system can be written as follows.
\begin{equation}
   I(k+1) = I(k) M^{u(k)}, \;\; \forall k, 0\leq k\leq K-1.
\end{equation}


Moreover, if we denote by $X(t)$ the column vector giving the function value for each state, in such a way that
the $p$th component of $X(t)$ gives $x_p(t)$, denote by $C^{u(t)}$ the column vector giving the payoffs at each
state, in such a way that the $p$th component of $C^{u(t)}$ gives $c_p^{u(t)}$, and denote by $\Phi$ the column vector
giving the payoffs at the final states, in such a way that the $p$th component of $\Phi$ gives $\phi_p$, then the dynamic programming 
system~(\ref{ddp6})-(\ref{ddp7}) is written as follows.
\begin{align}
  & X(K) = \Phi, \label{ddp8} \\
  & \left( X(t) \right)_p = \min_{u\in\mathcal U} \left( \left( C^{u(t)}\right)_p + \left( M^{u(t)} X(t+1)\right)_p\right), \;\; \forall p, 1\leq p\leq n. \label{ddp9}
\end{align}


Let us now define the following matrix $A(t)$ with components in $\mathbb R \cup \{+\infty\}$.
$$A_{ij}(t)\begin{cases}
	      c_i^{u(t)} & \text{ if } M_{ij}^{u(t)} = 1, \\
	      +\infty & \text{ otherwise}.
           \end{cases}$$
Then we can easily check that the dynamic programming system~(\ref{ddp8})-(\ref{ddp9}) can be written as follows.
\begin{align}
  & X(K) = \Phi, \label{ddp10} \\
  & X(t) = A(t) \otimes X(t+1), \;\; \forall t, 0 \leq t \leq K-1, \label{ddp11}
\end{align}
where $X_i(t), 1\leq i \leq n, 0\leq t \leq K$ belong now to $\mathbb R \cup \{+\infty\}$, and where the operator $\otimes$ 
(which is the Min-plus algebra matrix product operator; see section~\ref{sec-mp} below) is defined as follows.
$$ \left( A \otimes B \right)_{ij} := \min_{1\leq p\leq l} \left( A_{ip} + B_{pj}\right), \;\; \forall i, 1\leq i\leq n, \forall j, 1\leq j\leq m.$$
$A$ and $B$ are $n\times l$ and $l\times m$ Min-plus matrices respectively; ; see section~\ref{sec-mp} below.

An important remark here is that the dynamic programming system associated to deterministic optimal control problem 
is written linearly in the Min-plus algebra.
One of the well known deterministic optimal control problems for which such formalism can be applied is the shortest path problem 
in a given network.

\section{Max-plus algebra}
\label{sec-mp}

As mentioned above, dynamic programming systems associated to deterministic optimal control problems can be written linearly in the Max-plus or Min-plus algebra.
We present here some reviews on the Max-plus algebra. 
The main traffic model we propose in Chapter~\ref{chap-maxtrain} is written in the Max-plus algebra of square matrices of polynomials.
We recall here the construction of this algebraic structure, and give some results which we used in the analysis
of our models.

\subsection{Max-plus algebra of scalars ($\mathbb R_{\max}$)}

Max-plus algebra~\cite{BCOQ92} is the idempotent semi-ring $(\mathbb R \cup \{-\infty\}, \oplus, \otimes)$
denoted by $\mathbb R_{\max}$, where the operations $\oplus$ and $\otimes$ are defined by:
$a \oplus b = \max\{a, b\}$ and $a \otimes b = a + b$. The zero element is $(-\infty)$ denoted by $\varepsilon$ and the
unity element is $0$ denoted by $e$. 

\subsection{Max-plus algebra of square matrices ($\mathbb R_{\max}^{n\times n}$)}

We have the same structure on the set of square matrices.
If $A$ and $B$ are two Max-plus matrices of size $n \times n$, the addition
$\oplus$ and the product $\otimes$ are defined by: $(A \oplus B)_{ij} = A_{ij} \oplus B_{ij} , \forall i, j$, and
$(A \otimes B)_{ij} = \bigoplus_k[A_{ik} \otimes B_{kj}]$. The zero and the unity matrices are still denoted
by $\varepsilon$ and $e$ respectively. 

A matrix $A$ is said to be \textit{reducible} if there exists a permutation matrix $P$ such that $P^T A P$ is lower block triangular.
A matrix that is not reducible is said to be \textit{irreducible}.

For a matrix $A\in \mathbb R_{\max}^{n\times n}$, a precedence graph $\mathcal G(A)$ is associated.
The set of nodes of $\mathcal G(A)$ is $\{1,2,\ldots,n\}$. There is an arc from node $i$ to node $j$
in $\mathcal G(A)$ if $A_{ji} \neq \varepsilon$.
A graph is said to be strongly connected if there exists a path from any node to any other node.
$A\in \mathbb R_{\max}^{n\times n}$ is irreducible if and only if $\mathcal G(A)$ is strongly connected~\cite{BCOQ92}.
That is, $A$ is irreducible if $\forall 1\leq i,j \leq n, \exists m\in\mathbb N, (A^m)_{ij} \neq \varepsilon$.

\subsection{Max-plus algebra of polynomials ($\mathbb R_{\max}[X]$)}

A (formal) polynomial in an indeterminate $X$ over $\mathbb R_{\max}$ is a finite sum $\bigoplus_{l=0}^p a_l X^l$ for some integer $p$
and coefficients $a_l\in \mathbb R_{\max}$. The set of formal polynomials in $\mathbb R_{\max}$ is denoted $\mathbb R_{\max}[X]$.
The support of a polynomial $\mathbf{f} = \bigoplus_{l=0}^p a_l X^l$ is $Supp(\mathbf{f}) = \{l, 0\leq l\leq p, a_l > \varepsilon\}$.
The degree of $\mathbf{f}$ is $Deg(\mathbf{f}) = \bigoplus_{l\in Supp(\mathbf{f})} l$.
We have the same algebraic structure of idempotent semi-ring on $\mathbb R_{\max}[X]$, where
the addition and the product of two polynomials $\mathbf{f}=\bigoplus_{l=0}^p a_l X^l$ and $\mathbf{g}=\bigoplus_{l=0}^q b_l X^l$ 
in $\mathbb R_{\max}[X]$ are defined as follows.
$$\mathbf{f}\oplus \mathbf{g} := \bigoplus_{l=0}^{deg(\mathbf{f})\oplus deg(\mathbf{g})} (a_l\oplus b_l) X^l,$$
$$\mathbf{f} \otimes \mathbf{g} := \bigoplus_{l=0}^{def(\mathbf{f})\otimes deg(\mathbf{g})} \left( \bigoplus_{i\otimes j=l} a_i\otimes b_j \right) X^l.$$
The zero element is $\varepsilon = \varepsilon X^0$ and the unity element is $e = eX^0$.
We notice that $\forall x\in \mathbb R_{\max}$, the valuation mapping $\varphi: \mathbf{f} \mapsto \mathbf{f}(x)$ is a homomorphism
from $\mathbb R_{\max}[X]$ into $\mathbb R_{\max}$.

\subsection{Max-plus algebra of polynomial square matrices $\left(\mathbb R_{\max}[X]\right)^{n\times n}$}

A polynomial matrix $A(X) \in \left(\mathbb R_{\max}[X]\right)^{n\times n}$ is a matrix with polynomial entries
$A_{ij}(X) = \bigoplus_{l=0}^p a_{ij}^{(l)} X^l$, where $a_{ij}^{(l)} \in \mathbb R_{\max}, \forall i,j, 0\leq i,j\leq n$ and
$\forall l, 1\leq l\leq p$.
$\left(\mathbb R_{\max}[X]\right)^{n\times n}$ is an idempotent semiring. The addition and the product are defined as follows.
$$(A(X) \oplus B(X))_{ij} = A_{ij}(X) \oplus B_{ij}(X), \forall i, j,$$
$$(A(X) \otimes B(X))_{ij} = \bigoplus_k[A_{ik}(X) \otimes B_{kj}(X)], \forall i, j.$$
The zero and the unity matrices are still denoted by $\varepsilon$ and $e$ respectively.
We also have, $\forall x\in \mathbb R_{\max}$, the valuation mapping $\varphi: A(X) \mapsto A(x)$ is a homomorphism 
from $\left(\mathbb R_{\max}[X]\right)^{n\times n}$ into $\mathbb R_{\max}^{n\times n}$.

A polynomial matrix $A (X) \in \left(\mathbb R_{\max}[X]\right)^{n\times n}$ is said to be \textit{reducible} if there exists a permutation
matrix $P(X)\in \left(\mathbb R_{\max}[X]\right)^{n\times n}$ such that $P(X)^T A(X) P(X)$ is lower block triangular.

Irreducibility of a matrix depends only on its support, that is the pattern of nonzero entries of the matrix.
As for finite values of $X$, the support of a matrix is preserved by the homomorphism valuation map $\varphi$,
then $A (X) \in \left(\mathbb R_{\max}[X]\right)^{n\times n}$ is irreducible if and only if $A(e)\in \left(\mathbb R_{\max}\right)^{n\times n}$ is so.
Therefore, $A(X) \in \left(\mathbb R_{\max}[X]\right)^{n\times n}$ is irreducible if and only if $\mathcal G(A(e))$ is strongly connected.

For a polynomial matrix $A(X) \in \left(\mathbb R_{\max}[X]\right)^{n\times n}$, a precedence graph $\mathcal G(A(X))$ is associated.
As for graphs associated to square matrices in $\mathbb R_{\max}^{n\times n}$, the set of nodes of $\mathcal G(A(X))$ is
$\{1,2,\ldots,n\}$. There is an arc $(i,j,l)$ from node $i$ to node $j$ in $\mathcal G(A(X))$ if $\exists l, 0\leq l\leq p, a^{(l)}_{ji} \neq \varepsilon$.
Moreover, a \textit{weight} $W(i,j,l)$ and a \textit{duration} $D(i,j,l)$ are associated to every arc $(i,j,l)$ in the graph,
with $W(i,j,l) = (A_l)_{ij} \neq \varepsilon$ and $D(i,j,l) = l$. Similarly, a weight, resp. duration of a cycle (directed cycle) in the graph is
the standard sum of the weights, resp. durations of all the arcs of the cycle. 
Finally, the \textit{cycle mean} of a cycle $c$ with a weight $W(c)$ and a duration $D(c)$ is $W(c)/D(c)$.

\subsection{Homogeneous linear Max-plus algebra systems}
\label{subsec-hlmp}

We are interested in the model for the first train dynamics we propose in Chapter~\ref{chap-maxtrain} in the dynamics of a homogeneous $p$-order max-plus system
\begin{equation}\label{eq-mp1-1}
  x(k) = \bigoplus_{l=0}^p A_l \otimes x(k-l),
\end{equation}
where $x(k), k\in \mathbb Z$ is a sequence of vectors in $\mathbb R_{\max}^n$, and $A_l, 0\leq l\leq p$ are matrices in $\mathbb R_{\max}^{n\times n}$.
If we define $\gamma$ as the back-shift operator applied on the sequences of vectors  in $\mathbb Z$, such that: $\gamma x(k) = x(k-1)$,
and then more generally $\gamma^l x(k) = x(k-l), \forall l\in\mathbb N$, then~(\ref{eq-mp1-1}) is written as follows.
\begin{equation}\label{eq-mp2}
  x^k = \bigoplus_{l=0}^p \gamma^l A_l x^k = A(\gamma) x^k,
\end{equation}
where $A(\gamma)=\bigoplus_{l=0}^p \gamma^l A_l \in \left(\mathbb R_{\max}[\gamma]\right)^{n\times n}$ is a polynomial matrix in the back-shift
operator $\gamma$; see~\cite{BCOQ92,Gov07} for more details.

\begin{definition}
    $\mu \in \mathbb R_{\max} \setminus \{\varepsilon\}$ is said to be a \textit{generalized} eigenvalue~\cite{CCGMQ98} of $A(\gamma)$, with associated
    \textit{generalized} eigenvector $v\in \mathbb R_{\max}^n \setminus \{\varepsilon\}$, if $A(\mu^{-1}) \otimes v = v$,
    where $A(\mu^{-1})$ is the matrix obtained by evaluating the polynomial matrix $A(\gamma)$ at $\mu^{-1}$.
\end{definition}    

\begin{theorem} \cite[Theorem 3.28]{BCOQ92} \cite[Theorems 7.4.1 and 7.4.7]{Gov05} \label{th-mpa}
  Let $A(\gamma) = \oplus_{l=0}^p A_l\gamma^l$
  be an irreducible polynomial matrix with acyclic
  sub-graph $\mathcal G(A_0)$. Then $A(\gamma)$ has a unique generalized eigenvalue $\mu > \varepsilon$ and finite eigenvectors $v > \varepsilon$
  such that $A(\mu^{-1}) \otimes v = v$, and $\mu$ is equal to the maximum cycle mean of $\mathcal G(A(\gamma))$, given as follows.
  $\mu = \max_{c\in\mathcal C} W(c) / D(c)$,
  where $\mathcal C$ is the set of all elementary cycles in $\mathcal G(A(\gamma))$.
  Moreover, the dynamic system $x^k = A(\gamma) x^k$ admits an asymptotic average growth vector (also called here cycle time vector) $\chi$ 
  whose components are all equal to $\mu$.
\end{theorem}

\section{Min-plus system theory}

We give in this section a short review of the main definitions and results of the min-plus system theory which we need
for the models proposed in Chapter~\ref{chap-roadnetcal}.
The main variables we use here are the cumulated traffic flows functions of time, which we denote with capital letters.
A dynamic system for the traffic is then seen as a system with input signals (car inflows) and output signals
(car outflows). The network calculus theory associates an arrival and a service curves to a such system, and
derives from those curves performance bounds like upper bounds of the delay of passing through the system.
An arrival curve upper-bounds the arrival inflows to the system, while a service curve
lower-bounds the guaranteed service and then the departure outflows from the system. In the following,
we review these two notions of arrival and service curves in the one-dimensional case (for systems with one
arrival inflow, and one departure outflow). 

\subsection{Min-plus algebra}

As for the max-plus algebra, the min-plus algebra of scalars, denoted $\mathbb R_{\min}$,
is the set $\mathbb R \cup \{+\infty\}$ endowed with the min-plus addition $\oplus$
such that $a \oplus b := \min (a,b)$ and the product $\otimes$ such that 
$a \times b := a + b$. The structure $\mathbb R_{\min}$ is a commutative dioid (idempotent semiring)~\cite{BCOQ92}.
The zero element is still denoted $\varepsilon = +\infty$~\footnote{$\varepsilon$ denotes $-\infty$ or $+\infty$ depending
on the structure considered, $\mathbb R_{\max}$ or $\mathbb R_{\min}$ respectively.}
and the identity element is $e=0$.

\subsubsection*{Min-plus algebra of functions.}

We are interested here in the min-plus algebra of functions in $\mathbb Z$, and then in the min-plus algebra
of matrices of functions in $\mathbb Z$.
We denote by $\mathcal F$ the set of functions $f$ indexed by $t\in\mathbb N$, such that $f(0)\geq 0$ and $f$ is increasing in $\mathbb N$.
Thus $f$ is non-negative.
We endow $\mathcal F$ with two intern operations: the addition $\oplus$ (element-wise minimum) and the multiplication $*$ (minimum convolution),
defined as follows.
\begin{itemize}
   \item addition: $(f\oplus g)(t) := \min(f(t),g(t))$.
   \item product: $(f*g)(t) := \min_{0\leq s\leq t} (f(s)+g(t-s))$.
\end{itemize}
The algebraic structure $(\mathcal F,\oplus,*)$ is a dioid, see~\cite{BCOQ92, Cha00, LT01}.
The zero element $\varepsilon$ and the identity element $e$ are given as follows.
\begin{itemize}
  \item $\varepsilon(t) = +\infty, \forall t\in \mathbb N$.
  \item $e(0) = 0$ and $e(t) = +\infty, \forall t > 0$.
\end{itemize} 
For, $f \in \mathcal F$, $f^k$ denotes the power operation with respect to the product $*$.
$$f^k = f*f^{k-1}, \text{ with } f^0 = e.$$
The sub-additive closure on $\mathcal F$ is then defined as follows.
$$f^*=\bigoplus_{k\geq 0} f^k.$$
We call deconvolution the operation $\oslash$ in $\mathcal F$ defined as follows.
$$(f \oslash g)(t) := \sup_{s\geq 0} \left( f(t+s) - g(s) \right).$$
We consider the two signals $\gamma^p$ (the gain signal) and $\delta^{\tau}$ (the time-shift signal) in $\mathcal F$.
\begin{itemize}
  \item $\gamma^p(t) = \begin{cases}
                          p & \text{ for } t = 0,\\
                          +\infty & \text{ for } t > 0.
                       \end{cases}$,
  \item $\delta^{\tau}(t) = \begin{cases}
                          0 & \text{ for } 0\leq t \leq T,\\
                          +\infty & \text{ for } t > T.
                       \end{cases}$.
\end{itemize}
\begin{figure}[h]
    \centering
    \includegraphics[scale=1]{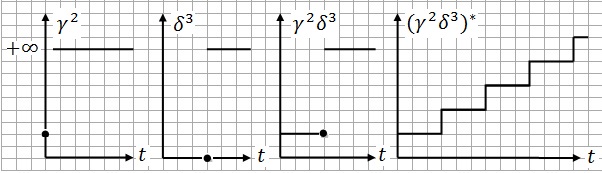}
    \caption{The signals $\gamma^2, \delta^3, \gamma^2\delta^3$ and $(\gamma^2\delta^3)^*$ respectively.}
    \label{gamma}
\end{figure} 
It is easy to check that (see Figure~\ref{gamma}).
\begin{align}
     \left(\gamma^p \delta^{\tau}\right)^* & \geq (p/ \tau) t, \label{res1} \\
     \gamma^{p_2}\delta_{\tau_2} \left(\gamma^{p_1} \delta^{\tau_1}\right)^* & \geq (p_1 / \tau_1)\left(t-(\tau_2 - \tau_1 p_2/p_1)\right)^+.  
        \label{res2}
\end{align}  
From Figure~\ref{gamma}, we can easily see that $(\gamma^2\delta^3)^*(t) \geq (2/3) t, \forall t\geq 0$.

We denote by $\mathcal F_0$ the set of functions $f$ indexed by $t\in\mathbb N$, such that $f(0) = 0$ and $f$ is increasing in $\mathbb N$.
The structure $(\mathcal F_0,\oplus, *)$ is also a dioid for which the zero $\varepsilon$ and the unity $e$ elements coincide: 
$\varepsilon(0)=e(0)=0$, and $\varepsilon(t) = e(t) = +\infty, \forall t>0$.

\subsubsection*{The dioid of matrices of functions.}

We denote by $\mathcal F^{n\times n}$ the set of $n\times n$ matrices with elements in $\mathcal F$.
The addition, still denoted by $\oplus$ is the element-wise minimum, and the product, still denoted
by $*$ is defined as follows.
$$(F*G)_{ij} = \bigoplus_{0\leq k\leq n}(F_{ik}*G_{kj}).$$

The zero element, still denoted by $\varepsilon$, is the $n\times n$ matrix with an $\varepsilon$
on all the entries. The unity element the $n\times n$ matrix with an $e$ on every diagonal entry,
and an $\varepsilon$ elsewhere. It is easy to check that we have again a dioid structure.

The sub-additive closure on $\mathcal F^{n\times n}$ is defined as follows.
$$F^*=\bigoplus_{k\geq 0} F^k,$$
where $F^k$ denotes the power operation with respect to the product $*$ on $\mathcal \mathcal F^{n\times n}$.
$$F^k = F*F^{k-1}, \text{ with } F^0 = e.$$

Let us consider the system $Y = f*Y \oplus U$ on the variable $Y \in \mathcal F^{n}$, where $f \in \mathcal F^{n\times n}$ and $U \in \mathcal F^{n}$.
A subsolution for that system is a $Y$ satisfying $Y \leq f*Y \oplus U$, where $F \leq G$
is defined here by $F \oplus G = F$.
Then, $Y'\in \mathcal F^{n}$ is called the maximum subsolution of the system if it is a subsolution of it, and if
for any subsolution $Y$ of that system, we have $Y \leq Y'$.

\begin{theorem}(Linear system with feedback){\cite{BCOQ92, Cha00}}\label{th-linsys}  
  The maximum subsolution, in the vector of signals $Y$, of the system $Y = f*Y \oplus U$ is given by $Y = f^* U$.
  If $f(0)>0$, then $Y = f^* U$ is unique. Under the same condition we have 
  $$Y \geq f*Y \oplus U \Rightarrow Y \geq f^* U.$$
\end{theorem} 

\section{Network calculus}

In this section we give a short review on the single-in-single-out (SISO) systems. 
Let us consider a system (seen as a server) with an arrival cumulated flow $U\in F$, and a departure
cumulated flow $Y\in F$.

\begin{figure}[h]
    \centering
    \includegraphics[scale=0.5]{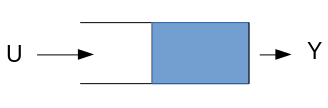}
    \caption{A server.}
    \label{server}
\end{figure} 

\begin{definition}(\cite{Cruz91a,Cruz91b,Cha00,LT01})
  \begin{itemize}
    \item The backlog $B(t)$ at time $t$ is defined $B(t) = U(t) - Y(t)$.
    \item The (virtual) delay $d(t)$ caused by the server at time $t$ is defined \\
       $d(t) = \inf\{h\geq 0, Y(t+h) \geq U(t)\}$.
    \item $\alpha$ is a maximum arrival curve for $U$, if $U \leq \alpha * U$.  
    \item $\underline{\alpha}$ is a minimum arrival curve for $U$, if $U \geq \underline{\alpha} \; \underline{*} \; U := \max_{0\leq s\leq t} (\underline{\alpha}(s) + U(t-s))$).
    \item $\beta$ is a service curve for the server, if $Y \geq \beta * U$.
    \item $\beta$ is a strict service curve for the server if for any busy time period $(s,t)$, we have \\
      $Y(t) - Y(s) \geq \beta(t-s)$.
  \end{itemize}
\end{definition}

Let us consider the following notations for affine and rate-latency curves.
\begin{align}
  & \Lambda(r,s)(t) := rt + s, \quad \text{ Affine curve}. \label{aff} \\
  & \lambda(R,T)(t) := R(t-T)^+, \quad \text{ Rate-latency curve}. \label{rl}
\end{align}

We notice that it is often used to upper-bound maximum arrival curves by affine curves~(\ref{aff}), and
lower-bound the service curves by rate-latency curves~(\ref{rl}).
The following result gives three bounds. 
\begin{theorem}(\cite{Cruz91a,Cruz91b,Cha00,LT01})
  The backlog $B$, the delay $d$ and the departure flow $Y$ are bounded as follows.
  \begin{itemize}  
     \item $B(t) \leq \sup_{s\geq 0} \left( \alpha(s) - \beta(s)\right) =: (\alpha\oslash\beta)(0), \forall t\geq 0$.
     \item $d(t) \leq \sup_{s\geq 0} \left\{ \inf \left\{ h \geq 0, \beta(s+h) \geq \alpha(s) \right\} \right\}, \forall t\geq 0$.
     \item $Y \leq (\alpha \oslash \beta) * Y$, ~~~~~i.e. $(\alpha\oslash\beta)$ is a maximum arrival curve for $Y$.
  \end{itemize}
\end{theorem} 
  
For the derivation of performance bounds from arrival and service curves in the one-dimension case, see~\cite{Cruz91a,Cruz91b,Cha00,LT01}.

\chapter[Max-plus model for the train dynamics]{Max-plus algebra model for the train dynamics in metro line systems}
\label{chap-maxtrain}

We present in this chapter a new approach for modeling and control of the train dynamics in metro line systems.
This approach has been firstly presented in~\cite{FNHL16} with the main ideas including a first model based on the Max-plus
algebra, for the modeling of the train dynamics in a linear metro line without taking into account the passenger 
travel demand. This first model is the basis of the approach. The references~\cite{FNHL16,FNHL18} include also two
other models (presented in Part~\ref{part-stoch} of this dissertation) for the modeling of the train dynamics taking
into account the passenger travel demand through average passenger arrival rates to platforms; see also~\cite{Far18}.
We present also in this chapter some extensions of the approach, realized with Master's and PhD degree students from
Ecole des Ponts ParisTech (Enpc).

The first extension~\cite{SFCLG17} models the train dynamics in a metro line with a junction. It has permitted
the derivation of analytic formulas giving the phase diagrams of the train dynamics in such a system.
This allowed us to understand wholly the physics of the traffic in this case, and, in particular, to 
derive the effect of the junction on the traffic phases and on the traffic control. An application of the 
the model on the metro line~13 of Paris has confirmed our analytic results.

The second extension~\cite{SFLG18_acc} introduces the running times (or train speeds) as control variables in addition to
the train dwell times at platforms. This extension permits in particular to respond easily to increases
of the passenger travel demand, by extending train dwell times, and then compensate by shortening train run times
on inter-stations. We provided here thresholds on the margins on the train run time, from which it is possible
to compensate an extension of train dwell times at platforms corresponding to a given increase of the passenger travel demand.

\section{Introduction}
\label{introduction}

Mass transit by metro is one of the most efficient, safe, and comfortable urban passenger transport modes.
However, it is also known to be naturally unstable in case where it is exploited at high frequencies~\cite{BCB91}.
Indeed, in the latter case, the capacity margins are reduced, and therefore train delays can be amplified and propagated
through the whole metro line.
In order to avoid such scenarios, the development of innovative approaches and methods for real-time railway traffic 
control is needed.
We present in this chapter discrete-event traffic models and control strategies for the train dynamics.
We provide here guarantees for the train dynamics stability, and derive analytically phase diagrams of the dynamics
with an interpretation of the traffic phases.

One of the main control parameters of the train dynamics is the train dwell
times at platforms. 
The passenger dynamics at platforms and inside the trains do have an effect on train dwell times at platforms,
and by that, on the whole train dynamics.
Indeed, accumulation of passengers at platforms and inside the trains induce additional constraints for the train dwell times.
This situation can be due to high level of passenger travel demand, or to delayed trains in case of incidents.
These additional constraints induce direct extensions of the train dwell times at platforms, which induces train delays.
The latter may then propagate in space and in time inducing secondary delays, more incidents, and so on. 

Several approaches have been developed (mathematical, simulation-based, expert system, etc.) for optimization and control of the train dynamics.
We cite here~\cite{Cur80,BCB91,Lee97,Ass04,EK04}.
Breusegem et al. (1991)~\cite{BCB91} developed interesting discrete event traffic and control models, pointing out the nature of traffic instability in metro lines.
The approach of~\cite{BCB91} is based on linear quadratic (LQ) control.
Another interesting approach is the one of Goverd (2007)~\cite{Gov07} who developed an algebraic model for the analysis of timetable stability and robustness against delays.
Other approaches studying and taking into account passenger dynamics and their effect on the train dynamics can be found in~\cite{Cao09,And,SZDZ13}.
For a recent overview on recovery models and algorithms for real-time railway disturbance and disruption management; see~\cite{Cac14}.

The modeling approach we propose here as well as the results are new. 
The first model~\cite{FNHL16,FNHL17a} of the train dynamics for a metro line without junction and without taking into account the
passenger travel demand, assumes bounds for the train dwell times at platforms, bounds on the train safe separation times on the line,
and nominal train run times on inter-stations. The train dynamics are then described with a discrete event system.
This description is new. Moreover, we characterize here the stationary regime of the train dynamics
by giving the conditions for its existence, and by deriving analytically the asymptotic average train-frequency 
as a function of the train density in the metro line, in the cases where a stationary regime exists.
This derivation is new and original.
It permits to distinguish the traffic phases of the train dynamics, and by that, determine 
the capacity of the metro line and the optimal number of running trains (which realizes the line capacity),
as functions of all the parameters of the line (bounds on the train dwell times at platforms,
nominal train run times at inter-stations and bounds on the train safe separation times).

The extensions presented here and developed in~\cite{SFCLG17} and~\cite{SFLG18_acc}, are also new.
The model of the train dynamics in a metro line with a junction~\cite{SFCLG17} derives the asymptotic
average train frequency as a function of the number of running trains on the line, and of the difference
between the number of running trains on the two branches of the line (the metro line has a central part and 
two branches crossing at a junction).
The traffic phases are then derived. This derivation is new and original. To our best knowledge, it is
the first description and analytic derivation of the physics of traffic on a metro line with a junction.
The second extension~\cite{SFLG18_acc} proposes a model of the train dynamics in a metro line without 
junction but where both train dwell times at platforms 
and train run times on inter-stations are controlled depending on the passenger arrival demand.
Similarly, the stationary regime of the train dynamics is characterized, and the asymptotic 
average train frequency is derived as a function of the number of running trains on the line, and of the level
of passenger arrival demand.
The traffic phases are then derived and interpreted. This derivation is also new and original.

This chapter is organized as follows.
In section~\ref{sec-mpm} we present the first model which describes the trains dynamics without taking into account the passenger travel demand.
The model considers lower bounds for train run times on inter-stations, train dwell times at platforms, and safe separation times between
successive trains. We show that the model is linear in the Max-plus algebra, and derive the stationary regime of the train dynamics.
In section~\ref{sec-ftd}, we define what we called \textit{fundamental traffic diagrams} of the train dynamics, by similarity to the road traffic.
In section~\ref{sec-tph}, we describe the traffic phases of the train dynamics derived for the model.
Finally and briefly, we present two extensions of the model which have been developed with Florain Schanzenbächer
(a PhD student at RATP and University of Paris Est).
Section~\ref{sec-ext1} presents an extension of the model to the case of metro lines with a junction; see~\cite{SFCLG17}.
Section~\ref{sec-ext2} presents an extension of the model to take into account the travel passenger demand; see~\cite{SFLG18_acc}.

\section{The Max-plus algebra model}
\label{sec-mpm}

The main references of this section are~\cite{FNHL16,FNHL17a,FNHL18}.
As mentioned above, this is a model for the train dynamics in a metro line.
It is the basis of all other extensions presented in the next sections.
This first model, written in the Max-plus algebra, takes into account minimum run, dwell and
safe separation time constraints, without any control of the train dwell times at platforms, and
without consideration of the passenger travel demand.

We show that the dynamics admit a stationary regime with a unique asymptotic average growth rate,
interpreted here as the asymptotic average train time-headway.
Indeed, the latter is derived as a Max-plus eigenvalue of the Max-plus linear system corresponding to the train dynamics.
Moreover, the asymptotic average train time-headway, dwell time, as well as safe separation time, are derived analytically, as
functions of the number of running trains on the metro line.
By that, three traffic phases of the train dynamics are clearly distinguished.

\begin{figure}[thpb]
      \centering
	  \includegraphics[scale = 0.6]{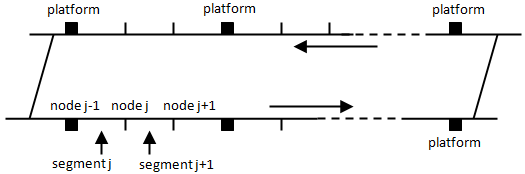}
      \caption{Representation of a linear metro line.}
      \label{fig-loop}
\end{figure}

We consider a linear metro line of $N$ platforms as shown in Figure~\ref{fig-loop}.
In order to model the train dynamics on the whole line, including the dynamics
on inter-stations, we discretize the inter-stations space, and thus the whole line, in segments
(or sections, or blocks). The length of every segment must be larger than the length of
one train. We then consider the following notations. \\~~

\noindent
\begin{tabular}{ll}
  $N$ & number of platforms.\\
  $n$ & number of all segments of the line.\\ 
  $m$ & number of running trains.\\
  $L$ & the length of the whole line. \\ 
  $b_j$ & $\in \{0,1\}$: boolean number of trains being on segment $j$ at time zero.\\
  $\bar{b}_j$ & $= 1 - b_j \in \{0,1\}$.
\end{tabular}~~\\
\begin{tabular}{ll}  
  $d^k_j$ & instant of the $k$th departure from node $j$. Notice that $k$ do not index trains, \\
          & but count the number of departures from segment $j$. \\
  $a^k_j$ & instant of the $k$th arrival to node $j$. Notice that $k$ do not index trains, but \\
          & count the number of arrivals to segment $j$.
\end{tabular}~~\\
\begin{tabular}{ll}  	  
  $r_j$ & the running time of a train on segment $j$, i.e. from node $j-1$ to node $j$.\\
  $w_j^k$ & $ = d_j^k - a_j^k$: train dwell time corresponding to the $k$th arrival to- and departure \\
          & from node $j$.
\end{tabular}~~\\
\begin{tabular}{ll}            
  $t_j^k$ & $ = r_{j} + w_j^k$: train travel time from node $j-1$ to node $j$, corresponding to the \\
          & $k$th arrival to- and departure from node $j$.
\end{tabular}~~\\
\begin{tabular}{ll}  		
  $g^k_j$ & $ = a_{j}^{k} - d^{k-1}_{j}$: node- (or station-) safe separation time (also known as close-in \\
          & time), corresponding to the $k$th arrival to- and $(k-1)$st departure from node $j$.\\
  $h_j^k$ & $ = d^k_j - d^{k-1}_j = g_j^k + w_j^k$~: departure time headway at node $j$, associated to the \\
          & $(k-1)$st and $k$th departures from node $j$.\\
  $s_j^k$ & $ = g_j^{k+b_j} - r_j$: a kind of segment safe separation time, not taking into \\
          & account the running time. $s_j^k$ is associated to segment $j$.
\end{tabular}~~\\

We also use underlined and over-lined notations to denote the maximum and minimum bounds of 
the corresponding variables respectively. Then
\begin{itemize}
    \item $\bar{r}_j, \bar{t}_j, \bar{w}_j, \bar{g}, \bar{h}_j$ and $\bar{s}_j$ denote maximum run, travel, dwell, safe separation,
             headway and $s$ times, respectively.
    \item $\underline{r}_j, \underline{t}_j, \underline{w}_j, \underline{g}, \underline{h}_j$ and $\underline{s}_j$
             denote minimum run, travel, dwell, safe separation, headway and $s$ times, respectively.
\end{itemize}
The average on $j$ and on $k$ (asymptotic) of those variables are denoted without any subscript or superscript.
Then 
\begin{itemize}
    \item $r, t, w, g, h$ and $s$ denote the average run, travel, dwell, safe separation, headway and $s$ times, respectively.
\end{itemize}

\begin{proposition}
  We have the following relationships.  
  \begin{align}
    & g = r + s, \label{form1} \\
    & t = r + w, \label{form2} \\
    & h = g + w = t + s = \frac{n}{m} t = \frac{n}{n-m} s. \label{form3}
  \end{align}    
\end{proposition}
\proof 
Indeed, (\ref{form1}) comes from the definition of $s_j^k$ and~(\ref{form2}) comes from
the definition of $t_j^k$. For~(\ref{form3}),
\begin{itemize}
  \item $h=g+w$ comes from the definition of $h_j^k$,
  \item $h=t+s$ comes from the definition of $t_j^k$ and $s_j^k$ and from $h=g+w$,
  \item $h=nt/m$ average train time-headway is given by the travel time of the whole line ($nt$) divided by
     the number of trains.
  \item $h=ns/(n-m)$ can be derived from $h=t+s$ and $h=nt/m$.
\end{itemize}  
\endproof

The running times $r_j$ of trains on every segment $j$, are considered to be constant (do not change with~$k$).
They can be calculated from the running times on inter-stations, and by means of given inter-station speed profiles
(in the free train flow case), depending on the characteristics of the line and of the trains running on it.
We then have 
\begin{align}
  & \underline{t}_j = r_j + \underline{w}_j, \forall j \\
  & \bar{t}_j = r_j + \bar{w}_j, \forall j.
\end{align}

Let us notice here that the variable $w_j^k$ denote dwell times at all nodes $j\in\{1,\ldots,n\}$
including non-platform nodes. The lower bounds $\underline{w}_j$ should be zero for the non-platform nodes $j$, and they
should be strictly positive for platform nodes. Therefore, the asymptotic average dwell time $w$ on all the nodes
is lower than the asymptotic dwell time on the platform nodes.
We also define $g^*, t^*$ and $s^*$ corresponding to platforms and distinguish them from $g, t$ and $s$.
\begin{itemize}
  \item $w^*$: asymptotic average dwell time at platform nodes.
  \item $g^*$: asymptotic average safe separation time on platform nodes.
  \item $t^* = r + w^*$: asymptotic travel time on segment $j$ upstream of platform node $j$.
  \item $s^* = g^* - r$: asymptotic average safe separation time on segment $j$ upstream of platform node $j$.
\end{itemize}

Another important remark related to the one above is that, as we consider in this first model constant running times on segments and on inter-stations,
then every train deceleration or stopping at the level of an inter-station, generally caused by an interaction with the train ahead,
is modeled here by a dwell time extension at one of the nodes at the considered inter-station.
In one of the extensions we present in the next sections, inter-station train running times are considered as control variables,
in addition to train dwell times at platforms.

The model we propose here writes the train dynamics basing on two time constraints:
\begin{itemize}
  \item A constraint on the travel time on every segment $j$.
     \begin{equation}\label{const1}
       d^k_j \geq d^{k-b_{j}}_{j-1} + \underline{t}_j.
     \end{equation}
     Constraint~(\ref{const1}) tells first that the $k$th departure from node $j$ corresponds to the same train as
     the $k$th departure from node $(j-1)$ in case where there is no train at segment $j$
     at time zero ($b_j=0$), and corresponds to the same train as the $(k-1)$st departure from node $(j-1)$
     in case where there is a train at segment $j$ at time zero.
     Constraint~(\ref{const1}) tells in addition that the departure from node $j$
     cannot be realized before the corresponding departure from node $(j-1)$ plus
     the minimum travel time $\underline{t}_j = \underline{r}_j + \underline{w}_j$ from node $j-1$ to node $j$.
  \item A constraint on the safe separation time at every segment~$j$.
     $$\begin{array}{ll}
         d_j^k - d^{k-\bar{b}_{j+1}}_{j+1} & = a_{j+1}^{k+b_{j+1}} - r_{j+1} - d^{k-\bar{b}_{j+1}}_{j+1} \\~~\\
                                           & = g_{j+1}^{k+b_{j+1}} - r_{j+1} \\~~\\
                                           & \geq \underline{g}_{j+1} - r_{j+1} = \underline{s}_{j+1}.
       \end{array}$$
     That is
     \begin{equation}\label{const2}
        d^k_j \geq d^{k-\bar{b}_{j+1}}_{j+1} + \underline{s}_{j+1}.
     \end{equation}
     Constraint~(\ref{const2}) tells first that, in term of safety, the $k$th departure from node $j$
     is constrained by the $(k-1)$st departure from node $(j+1)$ in case where there is no train at segment
     $(j+1)$ at time zero, and it is constrained by the $k$th departure from node $(j+1)$ in case where there is
     a train at segment $(j+1)$ at time zero. Constraint~(\ref{const2}) tells in addition that the $k$th
     departure from node $j$ cannot be realized before the departure constraining it from node $(j+1)$
     plus the minimum safety time at node $(j+1)$.
\end{itemize}

The model then assumes that the $k$th train departure from segment $j$ is realized as soon as both constraints~(\ref{const1}) and~(\ref{const2})
are satisfied. Therefore, the $k$th train departure time from segment $j$ is given as follows.
\begin{equation}\label{eq-d1}
  d^k_j = \max \{d^{k-b_{j}}_{j-1} + \underline{t}_j, d^{k-\bar{b}_{j+1}}_{j+1} + \underline{s}_{j+1} \}, \; k\geq 1, 1\leq j\leq n,
\end{equation}
where the index $j$ is taken with modulo $n$. That is to say that, for the two particular cases of $j=1$ and $j=n$, the dynamics are written as follows.
\begin{align}
  & d^k_1 = \max \{d^{k-b_{1}}_n + \underline{t}_1, d^{k-\bar{b}_2}_{2} + \underline{s}_2 \}, \quad k\in \mathbb N, \nonumber \\
  & d^k_n = \max \{d^{k-b_{n}}_{n-1} + \underline{t}_n, d^{k-\bar{b}_1}_{1} + \underline{s}_1 \}, \quad k\in \mathbb N. \nonumber
\end{align}

With Max-plus notations, and using the back-shift operator~$\gamma$, defined in 
section~\ref{subsec-hlmp}, the dynamics~(\ref{eq-d1})
are written as follows.
\begin{equation}\label{eq-d11}
  d_j = \underline{t}_j \gamma^{b_{j}} d_{j-1} \oplus \underline{s}_{j+1} \gamma^{\bar{b}_{j+1}} d_{j+1}, \quad 1\leq j\leq n.
\end{equation}

We denote by $d^k$ the vector with components $d^k_j$ for $j=1,\ldots, n$. The dynamics~(\ref{eq-d11}) are then
written as follows.
\begin{equation}
   d^k = A(\gamma) \otimes d^k,
\end{equation}
where $A(\gamma)$ is the following Max-plus polynomial matrix.
\small
$$A(\gamma) = \begin{pmatrix}
                 \varepsilon & \gamma^{\bar{b}_2}\underline{s}_2 & \varepsilon & \cdots & \varepsilon & \gamma^{b_1}\underline{t}_1 \\
                 \gamma^{b_2}\underline{t}_2 & \varepsilon & \gamma^{\bar{b}_3}\underline{s}_3 & \varepsilon & \cdots & \varepsilon \\
                  & \ddots & \varepsilon & \ddots & & \\
                  \varepsilon & \cdots & \gamma^{b_{j}}\underline{t}_j & \varepsilon & \gamma^{\bar{b}_{j+1}}\underline{s}_{j+1} & \varepsilon \\
                   & & & \ddots & \varepsilon & \\
                  \gamma^{\bar{b}_1}\underline{s}_1 & \varepsilon & \cdots & \varepsilon & \gamma^{b_{n}}\underline{t}_n & \varepsilon
              \end{pmatrix}$$
\normalsize

\begin{theorem}\label{th-mpm}
  The train dynamics~(\ref{eq-d1}) converges to a stable stationary regime with a unique average asymptotic growth vector,
  whose components are all equal and are
  interpreted here as the asymptotic average train time-headway $h$ given as follows.
  $$h = \max \left\{ \frac{\sum_j \underline{t}_j}{m}, \max_j (\underline{t}_j+\underline{s}_j), \frac{\sum_j \underline{s}_j}{n-m}\right\}.$$
\end{theorem}~\\

\begin{figure}[thpb]
      \centering
	  \includegraphics[scale = 0.42]{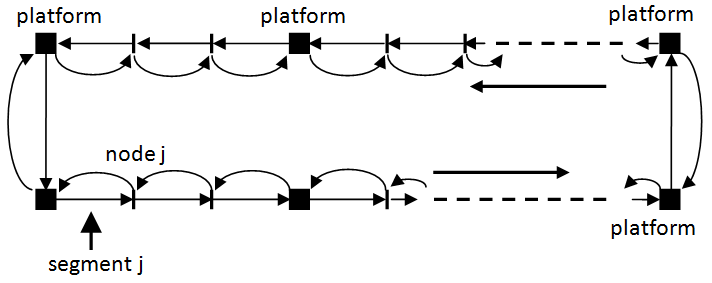}
      \caption{The graph $\mathcal G(A(\gamma))$.}
      \label{graph1}
\end{figure}

\proof
The graph $\mathcal G(A(\gamma))$ associated to the matrix $A(\gamma)$ is strongly 
connected; see Figure~\ref{graph1}.
Therefore, by Theorem~\ref{th-mpa}, we know that the asymptotic average growth 
vector of the dynamic system~(\ref{eq-d1}), whose components $h_j$
are interpreted here as the asymptotic average time-headway of the trains on segment $j$,
exists, and that all its components are the same
$h = h_j = \lim_{k\to +\infty} d^k_j/k, \forall j=1,\ldots,n$.
Moreover, $h$ coincides with the unique generalized eigenvalue of $A(\gamma)$,
given by Theorem~\ref{th-mpa} as the maximum cycle mean of the graph $\mathcal G(A(\gamma))$.
Three different elementary cycles are distinguished in $\mathcal 
G(A(\gamma))$; see Figure~\ref{graph1}.

\begin{itemize}
  \item The Hamiltonian cycle $c$ in the direction of the train running, with cycle mean
     $$\frac{W(c)}{D(c)} = \frac{\sum_j \underline{t}_j}{\sum_j b_j} = \frac{\sum_j \underline{t}_j}{m}.$$
  \item All the cycles $c_j$ of two links relying nodes $j-1$ and $j$, with cycle means
     $$\frac{W(c_j)}{D(c_j)} = \frac{(\underline{t}_j+\underline{s}_j)}{(b_j+\bar{b}_j)} = \underline{t}_j+\underline{s}_j, \quad \forall j.$$
  \item The Hamiltonian cycle $\bar{c}$ in the reverse direction of the train running, with cycle mean
     $$\frac{W(\bar{c})}{D(\bar{c})} = \frac{\sum_j \underline{s}_j}{\sum_j \bar{b}_j} = \frac{\sum_j \underline{s}_j}{n-m}.$$
\end{itemize}
\endproof

An important remark on Theorem~\ref{th-mpm} is that the asymptotic average 
train time-headway depends
on the average number of trains running on the metro line, without depending on 
the initial departure times of the trains (initial condition of the dynamic system).

\subsection{Fundamental traffic diagram for train dynamics}
\label{sec-ftd}

By similarity to the road traffic, one can define what is called \textit{fundamental traffic diagram} for the train dynamics.
In road traffic, such diagrams give relationships between car-flow and car-density on a road section; see for example~\cite{Far12, FHL13};
also extended to network (or macroscopic) fundamental diagrams; see for example~\cite{FGQ05, FGQ11, FGQ07, Far09, FGQ07b,FGQ11b}.

\begin{figure}[thpb]
      \begin{center}    
      \includegraphics[scale=0.6]{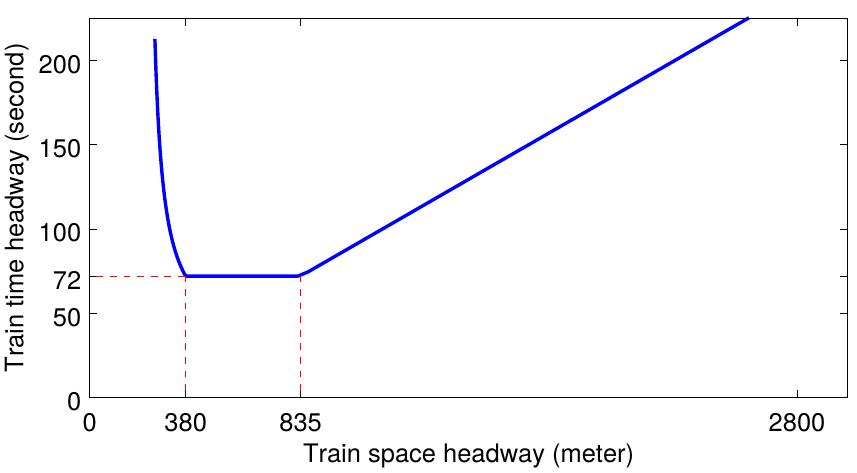} \\ ~~ \\
      \includegraphics[scale=0.6]{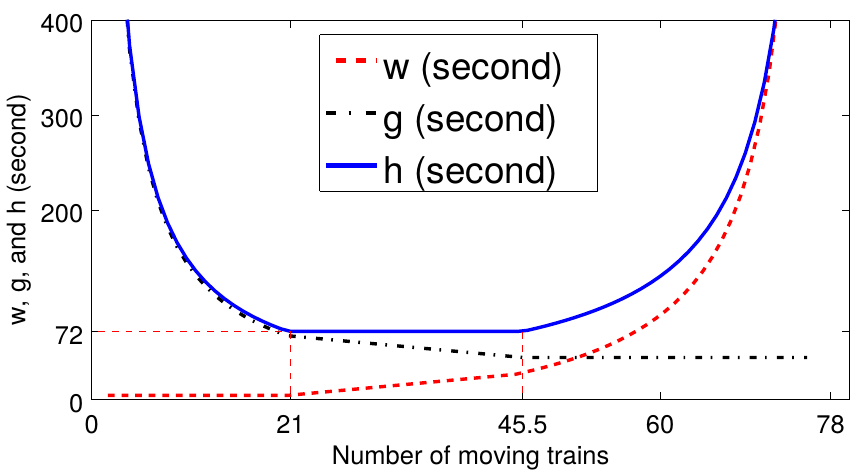} \\ ~~ \\
      \includegraphics[scale=0.6]{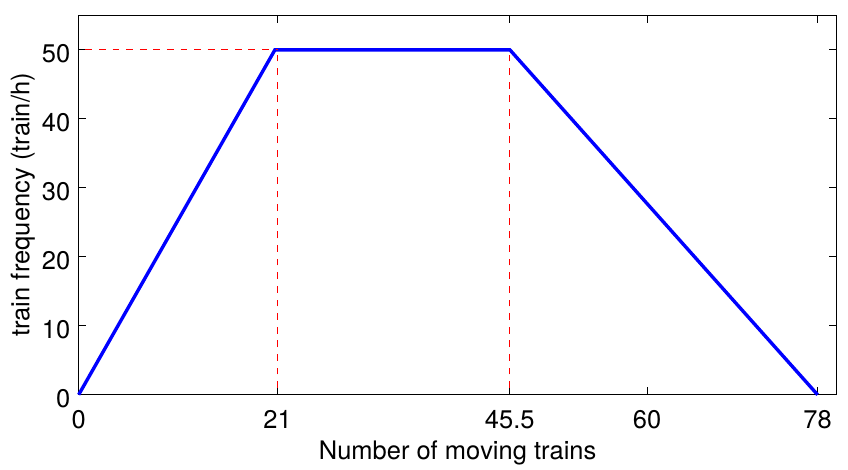}
      \caption{Analytical Phase diagrams for the train dynamics in a linear metro line.}
      \label{fig-diag1}
      \end{center}
\end{figure}

Let us first notice that the result given in Theorem~\ref{th-mpm} can be written as follows.
\begin{equation}\label{diag1}
  h(\sigma) = \max\left\{ \tau \sigma, h_{\min}, \frac{\omega}{\frac{1}{\underline{\sigma}} - \frac{1}{\sigma}}\right\},
\end{equation}
where $h$ is the asymptotic average train time-headway, and 
\begin{itemize}
    \item $\sigma := L/m$ is the average train space-headway,
    \item $\tau := \sum_j \underline{t}_j / L = 1/v$ is the inverse of the maximum train speed $v$,
    \item $h_{\min} := \max_j h_j = \max_j (\underline{t}_j + \underline{s}_j)$,
    \item $\omega := \sum_j \underline{s}_j / L$,
    \item $\underline{\sigma} := L/n$ is the minimum train space-headway.
\end{itemize}
Relationship~(\ref{diag1}) gives the asymptotic average train time-headway as a function of the average train space-headway;
see~Figure~\ref{fig-diag1}.

One can also write a relationship giving the average train time-headway as a function of the average train density $\rho := m/L = 1/\sigma$;
see~Figure~\ref{fig-diag1}.
\begin{equation}\label{diag2}
  h(\rho) = \max\left\{ \frac{\tau}{\rho}, h_{\min}, \frac{\omega}{\bar{\rho} - \rho}\right\},
\end{equation}
where 
\begin{itemize}
    \item $\bar{\rho} := n/L = 1/\underline{\sigma}$ is the maximum train density on the metro line.
\end{itemize}    
Let us know denote
\begin{itemize}
    \item $f=1/h$: the average train frequency (or flow) on the metro line.
\end{itemize}    
Then, from~(\ref{diag2}), we obtain a trapezoidal fundamental traffic diagram (well known in the road traffic) for the metro line;
see~Figure~\ref{fig-diag1}.
\begin{equation}\label{diag3}
  f(\rho) = \min \left\{v \rho,  f_{\max}, w' (\bar{\rho} - \rho)\right\},
\end{equation}
where 
\begin{itemize}
    \item $f_{\max} = 1/h_{\min}$ is the maximum train frequency over the metro line segments,
    \item $v = 1/\tau$ is the free (or maximum) train-speed on the metro line,
    \item $w' = 1/\omega$ is the backward wave-speed for the train dynamics.~\footnote{We use the
     notation $w'$, with a prime, for the backward wave speed, in order to distinguish 
     it with dwell time notation $w$.}
\end{itemize}
Relationships~(\ref{diag1}),~(\ref{diag2}) and~(\ref{diag3}) show how the asymptotic average train time-headway, 
and the asymptotic average train frequency change in function of the number of trains running on the metro line.
Moreover, they give the (maximum) train capacity of the metro line (expressed by the average train time-headway or by the average train frequency),
as well as the corresponding optimal number of trains. Furthermore, those relationships describe wholly the traffic phases of the train dynamics.

\begin{theorem}\label{cor-wg}
  The asymptotic average dwell time $w$ and safe separation time $g$ are given as follows.
  \begin{align}
    & w(\rho) = \max\left\{\underline{w}, \frac{h_{\min}}{\bar{\rho}}\;\rho - r, \frac{\omega}{\bar{\rho} - \rho} - \underline{g} \right\}.\label{diag-w}\\
    & g(\rho) = \max\left\{ \frac{\tau}{\rho} - \underline{w}, (r+h_{\min}) - \frac{h_{\min}}{\bar{\rho}}\;\rho, \underline{g} \right\}. \label{diag-g}
  \end{align}  
  where $\underline{w} = \sum_j \underline{w}_j /n, r = \sum_j r_j /n$ and $\underline{g} = \sum_j \underline{g_j} /n$.
\end{theorem}
\proof
We have
\begin{itemize}
  \item By (\ref{form2}) and (\ref{form3}), $w = t - r = (m/n) h - r$, then we replace $h$ using~(\ref{diag2}).
  \item By (\ref{form1}) and (\ref{form3}), $g = r + s = r + ((n-m)/n) h$, then we replace $h$ using~(\ref{diag2}). 
    Or directly form~(\ref{form3}), we have $g = h - w$, then we replace $h$ using~(\ref{diag2}) and $w$ using~(\ref{diag-w}). \endproof
\end{itemize}

Figure \ref{fig-diag1} illustrates the relationships~(\ref{diag1}),
(\ref{diag2}), (\ref{diag3}), (\ref{diag-w}) and (\ref{diag-g}) for a linear metro line of 9 stations (18 platforms),
inspired from the automated metro line~14 of Paris~\cite{ratp99}.
The parameters considered for the line are given in Table~\ref{tab-param}.

\begin{table}
\centering
\caption{Parameters of the metro line considered.}
\begin{tabular}{|l||l|}
  \hline
  Number of stations & 9 ($\Rightarrow $ 18 platforms)\\
  \hline
  Segment length & about 200 meters (m.) \\
  \hline
  Free train speed $v_{\text{run}}$ & 22 m/s (about 80 km/h) \\
  \hline 
  Train speed on terminus & 11 m/s (about 40 km/h) \\
  \hline
  Min. dwell time $\underline{w}$ & 20 seconds \\
  \hline
  Min. safety time $\underline{s}$ & 30 seconds \\
  \hline  
  Inter-station length (in meters) & \begin{tabular}{lr}
			       S.-Laz. $\to$ Mad. & 618 m. \\
			       Mad. $\to$ Pyr. & 712 m. \\
                               Pyr. $\to$ Cha. & 1359 m. \\
			       Cha. $\to$ G.-Lyo & 2499 m.\\
			       G.-Lyo $\to$ Ber. & 624 m.\\
			       Ber. $\to$ C. S. Emi. & 970 m. \\
			       C. S. Emi. $\to$ Bib. & 947 m.\\
			       Bib. $\to$ Oly. & 713 m.
			     \end{tabular} \\
  \hline                          
\end{tabular}
\label{tab-param}
\end{table}

According to Figure~\ref{fig-diag1}, the maximum average train frequency for the considered metro line is about 50 trains/hour,
corresponding to an average time-headway of 72 seconds. The optimal number of running trains to reach the
maximum frequency is 21 trains. We note that time-margins for robustness are not considered here. 

Formulas~(\ref{diag-w}) and~(\ref{diag-g}) are important for the control model we present 
in section~\ref{sec-stable} of chapter~\ref{chap-stochtrain}, where we consider $w$ as the control vector and $g$ as the
traffic state vector, and where the formulas~(\ref{diag2}),~(\ref{diag-w}), and~(\ref{diag-g})
of the max-plus traffic model are used.

\subsection{The traffic phases of the train dynamics}
\label{sec-tph}

Theorems~\ref{th-mpm} and~\ref{cor-wg}, and formulas~(\ref{diag1}), (\ref{diag2}) and~(\ref{diag3}) allow the description of the
traffic phases of the train dynamics~(\ref{eq-d1}). Three traffic phases are distinguished.

\textit{\textbf{Train free flow traffic phase}}. ($0 \leq \rho \leq f_{\max}/v$). During this phase, trains move freely on the line, which
    operates under capacity, with high average train time-headways. The average time-headway is a
    sum of the average minimum train dwell time with the average train safe separation time.
    The average train dwell time is independent of the number of
    running trains, while the average train time-headway as well as the average train safe separation time
    decrease rapidly with the number of running trains. We notice that the average train frequency increases linearly with respect to 
    the number of running trains. Similarly, the average train time-headway increases linearly with respect to the average space-headway.

\textit{\textbf{Maximum train-frequency traffic phase}}. ($f_{\max}/v \leq \rho \leq \bar{\rho} - f_{\max}/w'$).
    During this phase, the metro line operates at its maximum train-capacity. The latter is constant, i.e. independent of the
    number of running trains. The average train dwell time $w$ increases linearly with the number of the running trains.
    The average train safe separation time $g$ decreases linearly with the number of running trains. 
    The average train time-headway $h = g + w$ remains constant and independent of the number of running trains.
    The optimum number of running trains on the line is attained at the beginning of the this traffic phase.
    This optimum number is $L f_{\max}/v$, corresponding to train density $f_{\max}/v$.
    However, in the case where passenger arrivals are taken into account, it can be interesting to increase the number of running
    trains on the line. Indeed, although the average train time-headway remains constant during this phase, the average
    train dwell time increases, while the average train safe separation time decreases with the increasing of the number of trains
    running on the line. This induces less average train safe separation time, so less time for the accumulation of passengers on platforms,
    in one side, and more time to passengers to go onto the trains, on the other side,
    without affecting the average train time-headway; see Figure~\ref{fig-diag1}.

\textit{\textbf{Train congestion traffic phase}}. ($\bar{\rho} - f_{\max}/w' \leq \rho \leq \bar{\rho}$).
    During this phase, trains interact with each other and the metro line operates under capacity with high average train time-headways.
    The latter is given as the sum of the average train dwell time with the average minimum safe separation time.
    The latter is independent of the number of running trains, while the average train time-headways, as well as the average train dwell times
    increase rapidly with the number of running trains on the metro line. The average train frequency
    decreases linearly with the number of trains running on the metro line.

\section{Extension 1 - Train dynamics in a metro line with a junction}
\label{sec-ext1}

\begin{figure}[h]
  \centering
  \includegraphics[width=0.8\textwidth]{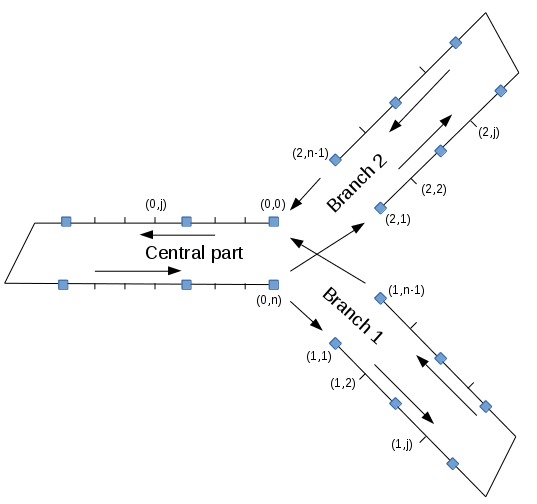} 
  \caption {A metro line with one junction.}
  \label{fig_line}
\end{figure}

The main reference of this section is~\cite{SFCLG17}.
We presented above a Max-plus model for the train dynamics on a linear metro line without junction.
We present here an extension of that model to the train dynamics on a metro line with a junction.
We assume that the junction is symmetrically operating (i.e. one train by two goes to the left/right at the divergent,
and one train by two enters from the left/right at the merge); see Figure~\ref{fig_line}.

We extend the train dynamics by adding the description of the train movements at the divergent and at the merge of the junction.
We show that, by means of a changing of variable, the model can be still written linearly in the Max-plus algebra.
We derive analytically the traffic phases of the train dynamics, where eight traffic phases are distinguished in this case.

We consider a metro line with one junction as in Figure~\ref{fig_line}.
The junction includes a divergent where trains go from the central part of the line to 
the branches, and a merge where trains go from the two branches to the central part of the line.
As above, and as in~\cite{FNHL16,FNHL17a,FNHL17b}, the line is discretized in a number of segments.
Segments and nodes are indexed as in Figure~\ref{fig_line}.
We adapted the notations of section~\ref{sec-mpm} as follows.~~\\~~

\noindent
\begin{tabular}{ll}
  $u$ & $\in \mathcal U = \{0,1,2\}$ indexes the central part if $u=0$, 
      the branch 1 if $u=1$, \\
      & and the branch 2 if $u=2$. \\
  $N_u$ & the number of platforms on part $u$ of the line. \\
  $n_u$ & the number of segments on part $u$ of the line. \\
  $J(u)$ & the set of indexes on part $u$ of the line.
         $J(u) = \{1,2, \ldots, n_u\}$.\\        
  $m_u$ & the number of trains being on the part $u$ of the  
        line, at time zero. \\
  $b_{(u,j)}$ & $\in \{0,1\}$. It is $0$ (resp. $1$) if there is no
              train (resp. one train) at segment \\
              & $j$ of part $u$.\\
  $\bar{b}_{(u,j)}$ & $= 1 - b_{(u,j)}$.
\end{tabular} \\~~\\

\noindent
\begin{tabular}{ll}
  $d^k_{(u,j)}$ & the $k^{\text{th}}$ departure time from node $j$, on part $u$
      of the line. Notice that $k$ \\
      & do not index trains, but count the number of train departures. \\  
  $a^k_{(u,j)}$ & the $k^{\text{th}}$ arrival time to node $j$, on part $u$
      of the line.
\end{tabular} \\~~\\
\begin{tabular}{ll}      
  $r{(u,j)}$ & the average running time of trains on segment $j$ 
      (between nodes $j-1$ \\
      & and $j$) of part $u$.\\
  $w^k_{(u,j)}$ & $=d^k_{(u,j)}-a^k_{(u,j)}$ the $k^{\text{th}}$ dwell time on node $j$ 
                on part $u$ of the line.\\
  $t^k_{(u,j)}$ & $=r_{(u,j)}+w^k_{(u,j)}$ the $k^{\text{th}}$ travel time from node
		$j-1$ to node $j$ on part $u$ of \\
		& the line.
\end{tabular} \\~~\\
\begin{tabular}{ll}		
  $g^k_{(u,j)}$ & $=a^k_{(u,j)}-d^{k-1}_{(u,j)}$ the $k^{\text{th}}$ safe separation time 
                (or close-in time) at node $j$ on\\
                & part $u$.\\
  $h^k_{(u,j)}$ & $=d^k_{(u,j)}-d^{k-1}_{(u,j)} = g^k_{(u,j)}+w^{k-1}_{(u,j)}$ the $k^{\text{th}}$
      departure time-headway at node $j$ on \\
      & part $u$.\\
  $s^k_{(u,j)}$ & $=g^{k+b_{(u,j)}}_{(u,j)}-r_{(u,j)}$.    
\end{tabular}

\subsection{The train dynamics besides the junction}

The train dynamics besides the junction are similar to the train dynamics in the case of 
a linear metro line without junction (section~\ref{sec-mpm}), and remain linear in the Max-plus algebra.
\begin{align}
   & d^k_{(u,j)} \geq d^{k-b_{(u,j)}}_{(u,j-1)} + \underline{t}_{(u,j)}, \; \forall k\geq 0, u\in\mathcal U, j \neq n_u. \\
   & d^k_{(u,j)} \geq d^{k-\bar{b}_{(u,j+1)}} _{(u,j+1)} + \underline{s}_{(u,j+1)}, \; \forall k\geq 0, u\in\mathcal U, j \neq n_u.
\end{align}

\subsection{Train dynamics at the divergence}

We assume that odd departures go to branch~1 while even ones go to branch~2.

The $k^{\text{th}}$ departures from the central part:
\begin{equation}\label{eq-d1_b}
  d^k_{(0,n)} \geq d^{k-b_{(0,n)}}_{(0,n-1)} + \underline{t}_{(0,n)}, \; \forall k\geq 0, \\
\end{equation}

\begin{equation}\label{eq-d2}
  d^k_{(0,n)} \geq \begin{cases}
                      d^{(k+1)/2-\bar{b}_{(1,1)}}_{(1,1)} + \underline{s}_{(1,1)} & \text{if } k \text{ is odd} \\ ~~ \\
                      d^{k/2-\bar{b}_{(2,1)}}_{(2,1)} + \underline{s}_{(2,1)} & \text{if } k \text{ is even} \\
                   \end{cases}
\end{equation}

The $k^{\text{th}}$ departures from the entry of branch 1:
\begin{eqnarray}
   d^k_{(1,1)} \geq d^{(2k-1)-b_{(1,1)}}_{(0,n)} + \underline{t}_{(1,1)}, \; \forall k\geq 0, \label{eq-d3}\\
   d^k_{(1,1)} \geq d^{k-\bar{b}_{(1,2)}} _{(1,2)} + \underline{s}_{(1,2)}, \; \forall k\geq 0. \label{eq-d4}
\end{eqnarray}

The $k^{\text{th}}$ departures from the entry of branch 2:
\begin{eqnarray}
   d^k_{(2,1)} \geq d^{2k-b_{(2,1)}}_{(0,n)} + \underline{t}_{(2,1)}, \; \forall k\geq 0, \label{eq-d5}\\
   d^k_{(2,1)} \geq d^{k-\bar{b}_{(2,2)}}_{(2,2)} + \underline{s}_{(2,2)}, \; \forall k\geq 0. \label{eq-d6}
\end{eqnarray}

\subsection{Train dynamics at the merge}

We assume that odd departures at node $(0,0)$ towards the central part correspond to trains coming from
branch~1 while even ones correspond to trains coming from branch~2.

The $k^{\text{th}}$ departures from the central part:
\begin{equation}\label{eq-m1}
  d^k_{(0,0)} \geq \begin{cases}
                      d^{(k+1)/2-b_{(1,n)}}_{(1,n-1)} + \underline{t}_{(1,n)} & \text{if } k \text{ is odd} \\ ~~ \\
                      d^{k/2-b_{(2,n)}}_{(2,n-1)} + \underline{t}_{(2,n)} & \text{if } k \text{ is even} \\
                   \end{cases}
\end{equation}
\begin{equation}\label{eq-m2}
  d^k_{(0,0)} \geq d^{k-\bar{b}_{(0,1)}}_{(0,1)} + \underline{s}_{(0,1)}, \; \forall k\geq 0, \\
\end{equation}

The $k^{\text{th}}$ departures from the entry of branch 1:
\begin{eqnarray}
   d^k_{(1,n-1)} \geq d^{k-b_{(1,n-1)}}_{(1,n-2)} + \underline{t}_{(1,n-1)}, \; \forall k\geq 0, \label{eq-m3}\\
   d^k_{(1,n-1)} \geq d^{(2k-1)-\bar{b}_{(1,n)}}_{(0,0)} + \underline{s}_{(1,n)}, \; \forall k\geq 0. \label{eq-m4}
\end{eqnarray}

The $k^{\text{th}}$ departures from the entry of branch 2:
\begin{eqnarray}
   d^k_{(2,n-1)} \geq d^{k-b_{(2,n-1)}}_{(2,n-2)} + \underline{t}_{(2,n-1)}, \; \forall k\geq 0, \label{eq-m5} \\
   d^k_{(2,n-1)} \geq d^{2k-\bar{b}_{(2,n)}}_{(0,0)} + \underline{s}_{(2,n)}, \; \forall k\geq 0. \label{eq-m6}
\end{eqnarray}

We have shown that the train dynamics at the divergence as well as the ones at the merge can be rewritten
linearly in the Max-plus algebra, by means of a changing of variables; see~\cite{SFCLG17} for the details. 
The whole dynamic system is then a Max-plus linear system.
As in section~\ref{sec-mpm}, we use the graphic approach to derive analytically the asymptotic average growth rate of the train dynamics,
which gives the average train time-headway at the central part, and at the two branches.
We have then obtained the traffic phases of the train dynamics on the metro line with a junction,
and by that we derive the effect of the junction on the physics of traffic; see~\cite{SFCLG17} for more details.

\section{The traffic phases}

The asymptotic average 
train time-headway and frequency depend on the total number of trains and on the difference between the number
of trains on the branches. Both parameters are time-invariant (with a time unit of two steps of the train dynamics).
Indeed, the difference in the number of trains between branches~2 and~1 is invariant in two steps of the train dynamics
because we assume here that the junction is symmetrically operating (i.e. one train by two goes to the left/right at the divergent,
and one train by two enters from the left/right at the merge).
Let us consider the following additional notations.~~\\~~ 

\begin{tabular}{ll}
  $m_u$ & the number of trains on part $u$ of the line at time zero. \\  
  $m$ & $= m_0+m_1+m_2$ the total number of trains on the line.\\
  $\Delta m$ & $= m_2 - m_1$ the difference in the number of trains between branches 2 and 1.\\
  $\bar{m}_u$ & $= n_u - m_u, \forall u\in\{0,1,2\}$.\\      
  $\bar{m}$ & $= \bar{m}_0 + \bar{m}_1 + \bar{m}_2$.\\
  $\Delta \bar{m}$ & $= \bar{m}_2 - \bar{m}_1$.\\
  $\underline{T}_u$ & $= \sum_j \underline{t}_{(u,j)}, \forall u\in\{0,1,2\}$. \\
  $\underline{S}_u$ & $= \sum_j \underline{s}_{(u,j)}, \forall u\in\{0,1,2\}$.
\end{tabular}~~\\~~

\noindent
The main results gives then the traffic phases of the train dynamics as follows.

\begin{theorem}(\cite{SFCLG17})\label{thm-1}
  The train dynamics admit a unique average growth rate
  which represents the average train time-headway $h_0$ at the central part of the metro line.
  Moreover we have
  $$ h_0 = h_1/2 = h_2/2 = \max \{ h_{fw}, h_{\min} ,h_{bw}, h_{br}\}, $$
  with\footnote{fw: forward, bw: backward, min: minimum, br: branches.}
  $$ h_{fw} = \max\left\{ \frac{\underline{T}_0 + \underline{T}_1}{m - \Delta m},
                          \frac{\underline{T}_0 + \underline{T}_2}{m + \Delta m} \right\}, $$
	$$ h_{\min} = \max \begin{cases}
		  \max_{u,j} (t_{(u,j)} + s_{(u,j)}) & \forall u \in \{0\},\\
		  \max_{u,j} (t_{(u,j)} + s_{(u,j)})/2 & \forall u \in \{1,2\},
	\end{cases}$$
  $$ h_{bw} = \max\left\{ \frac{\underline{S}_0 + \underline{S}_1}{\bar{m} - \Delta \bar{m}},
                          \frac{\underline{S}_0 + \underline{S}_2}{\bar{m} + \Delta \bar{m}} \right\}, $$
  $$ h_{br} =  \max\left\{\frac{\underline{T}_1 + \underline{S}_2}{2(n_2 - \Delta m)},
                  \frac{\underline{S}_1 + \underline{T}_2}{2(n_1 + \Delta m)}\right\}.$$
\end{theorem}
\proof See~\cite{SFCLG17}. \endproof

\begin{corollary}\label{cor1}
  The average train frequency $f_0$ on the central part and $f_1=f_2$ on the branches
  are given as follows:
  $$f_0 = 2 f_1 = 2 f_2 = \max \left\{ 0, \min \left\{ \frac{1}{h_{fw}}, \frac{1}{h_{\min}} ,\frac{1}{h_{bw}}, \frac{1}{h_{br}} \right\} \right\}.$$
\end{corollary}
\proof Directly from Theorem~\ref{thm-1} with $f_0 = 1/h_0$ and $f_0 \geq 0$. \endproof

Fig.~\ref{fig-2D} illustrates the results of Theorem~\ref{thm-1} and Corollary~\ref{cor1}.
We can clearly distinguish the eight traffic phases of the train dynamics.
For more details on the interpretation of those phase diagrams, please see~\cite{SFCLG17} and~\cite{SFLG18_acc}.

\begin{table}
\centering
\begin{tabular}{|c|}
  \hline \\
  \includegraphics[width=0.7\textwidth]{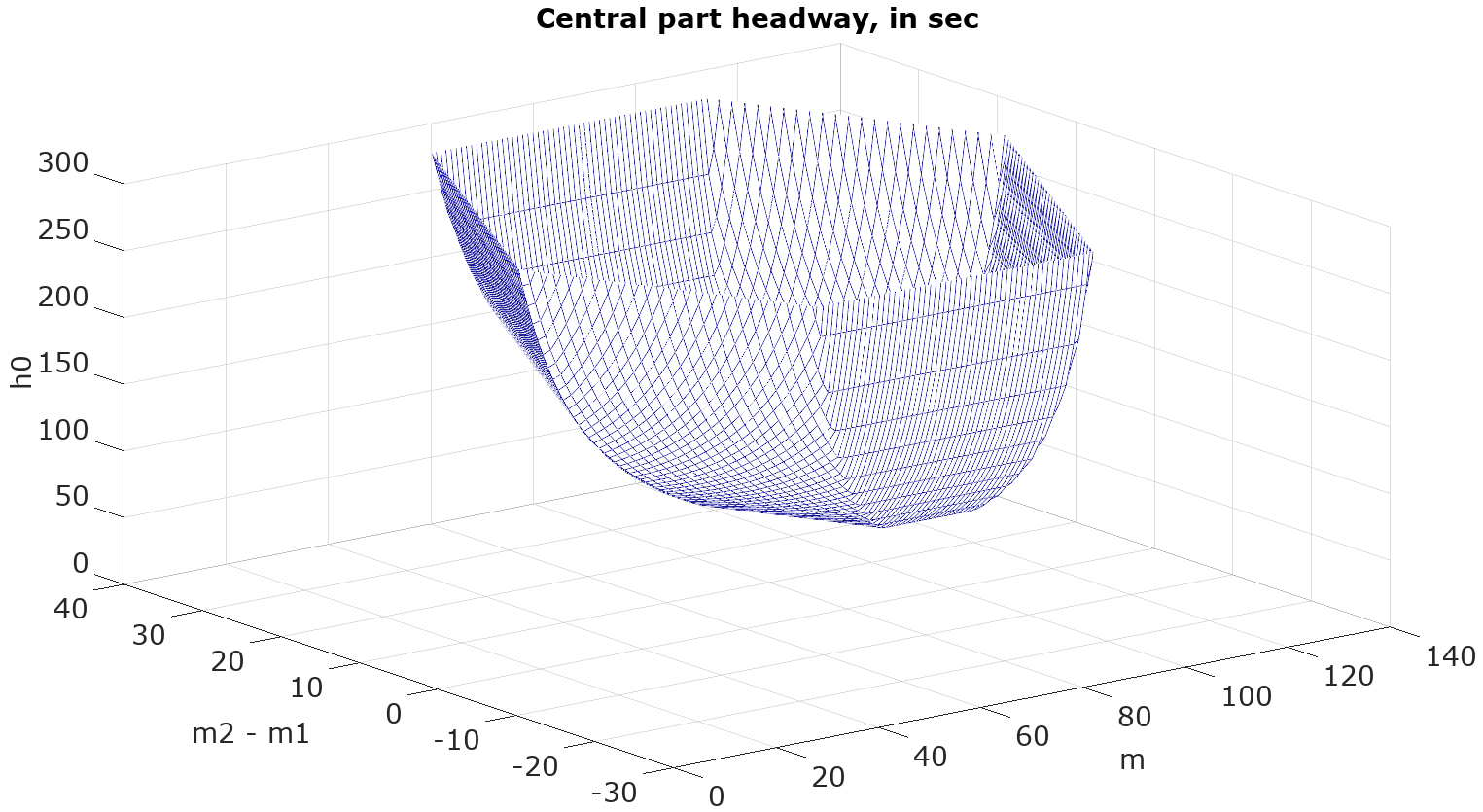}\\
  \hline \\
  \includegraphics[width=0.7\textwidth]{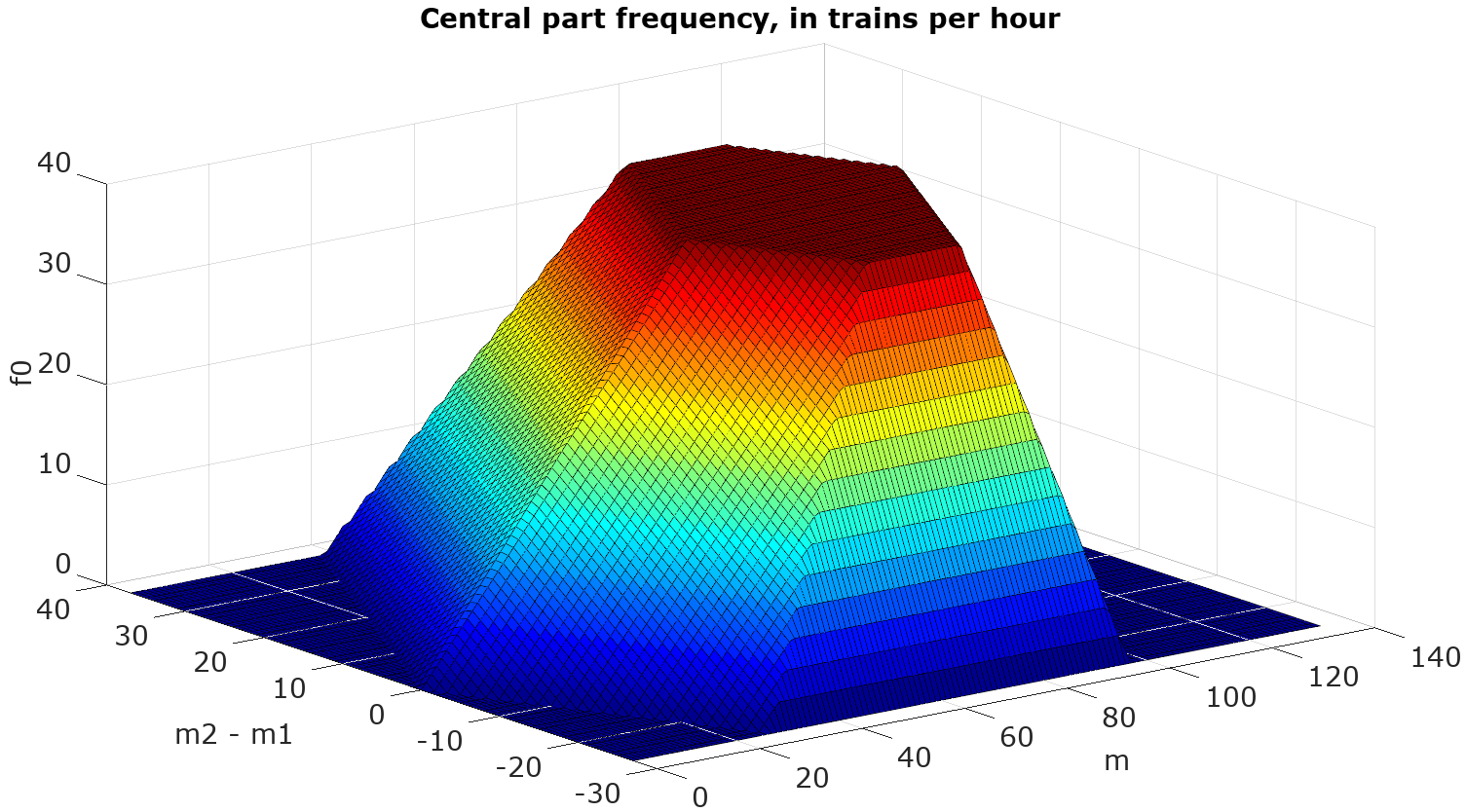}\\
  \hline \\
  \includegraphics[width=0.7\textwidth]{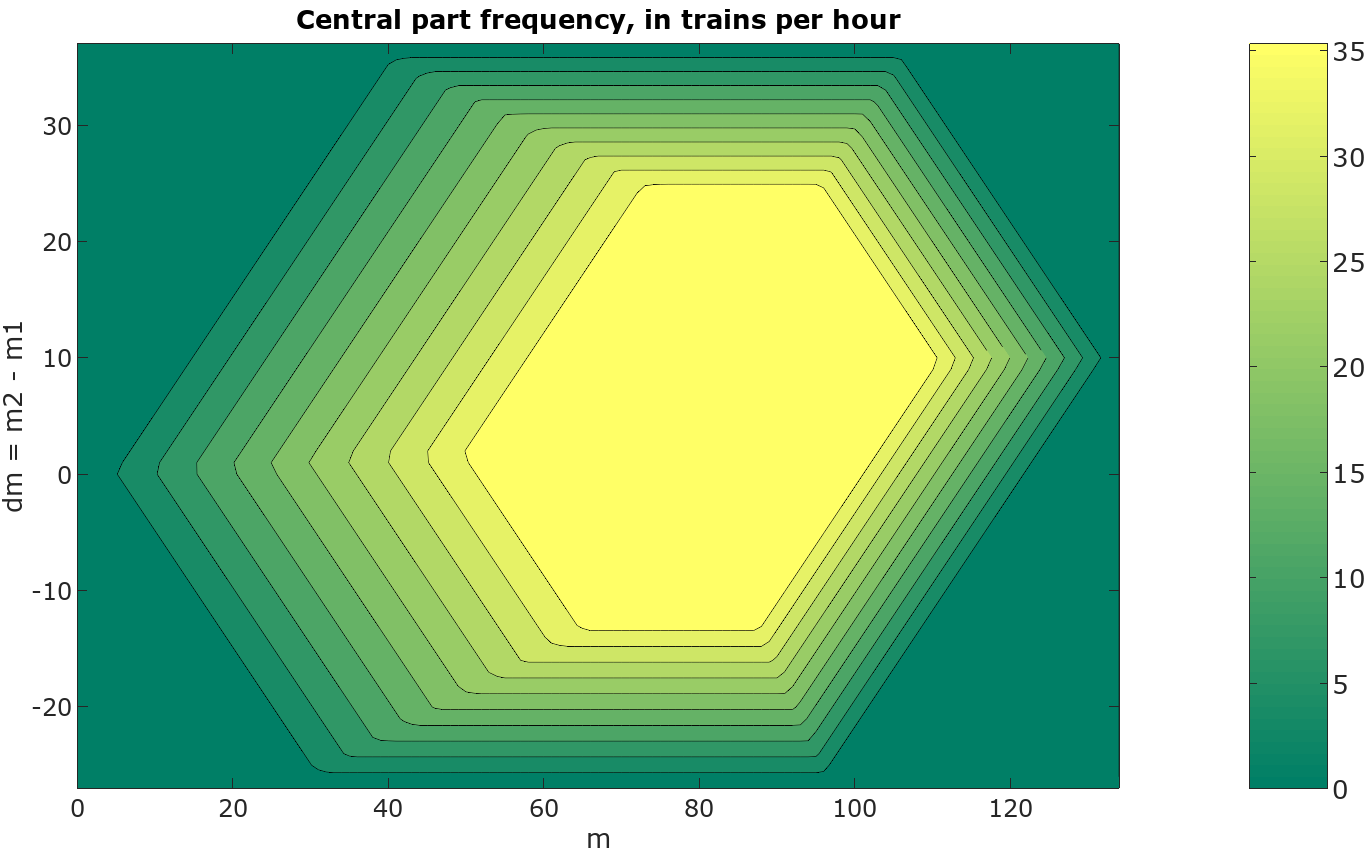}\\
  \hline 
\end{tabular}    
  \caption{The asymptotic average train time headway $h_0$ and frequency $f_0$, as functions of the total number $m$ of trains on the line 
          and of the difference $\Delta m$ of the number of trains on the two branches.
          We recognize the eight traffic phases of Theorem~\ref{thm-1} and Corollary~\ref{cor1}.} 
  \label{fig-2D}
\end{table}

\newpage
\section{Extension 2 - Train dynamics taking into account the travel demand}
\label{sec-ext2}

The main reference of this part is~\cite{SFLG18_acc}.
We propose in this part a second extension of the model~\ref{eq-d1} for the train dynamic modeling and control.
We present the model on a linear metro line without junction, but it is still extensible to a metro line with a junction. 
We consider the same train dynamics as in~(\ref{eq-d1}), but here the train dwell and run times are both controlled.
Two control laws are associated to the train dwell and run times, depending on the real-time train dynamics and on the
passenger arrival flow rates.

The idea here is to be able to extend train dwell times at platforms to respond to high levels of passenger travel demand,
and in counterpart to shorten the train run time on inter-stations in order to compensate, and avoid 
delays and their propagation.
We show that under given run time margins, the train dynamics remain stable against small perturbations relatively to the given
run time margins.
We show that, in this case, the train dynamics are still linear in the Max-plus algebra, and the 
asymptotic average train time-headway as derived analytically as a function of the number of running trains on the metro line
and of the level of the passenger travel demand.

We use here the same notations as in section~\ref{sec-mpm} for the train dynamics.
For the passenger travel demand, we use the following additional notations.\\~\\

\noindent
\begin{tabular}{ll}  
  $\lambda_{ij}$ & passenger travel demand from platform~$i$ to \\
      & platform~$j$, when $i$ and $j$ denote platforms, and\\
      & $\lambda_{ij} = 0$ if $i$ or $j$ is not a platform node.\\
  $\lambda^{\text{in}}_i$ & $=\sum_j \lambda_{ij}$ the average passenger arrival rate on \\
      & origin platform $i$ to any destination platform. \\      
  $\lambda^{\text{out}}_j$ & $=\sum_i \lambda_{ij}$ the average passenger arrival rate on \\
      & any origin platform to destination platform $j$. \\
  $\alpha^{\text{in}}_j$ & average passenger upload rate on platform $j$. \\~~\\      
  $\alpha^{\text{out}}_j$ & average passenger download rate on platform $j$. \\~~\\
  $\sum_i \lambda_{ij} h_i / \alpha_j^{out}$ & time needed for passenger download at platform $j$. \\~~\\
  $\lambda^{\text{in}}_j\; h_j / \alpha^{\text{in}}_j$ & time needed for passenger upload at platform $j$.
\end{tabular}  
\\~\\

The time needed for passenger upload at platform $j$ is $\lambda^{\text{in}}_j\; h_j / \alpha^{\text{in}}_j$.
The time needed for passenger download at platform $j$ is $\sum_{i} \lambda_{ij} h_{i} / \alpha_{j}^{out}$.
As in~\cite{SFLG18_acc}, we approximate it as follows.

$$\lambda^{\text{out}}_j\; h_j / \alpha^{\text{out}}_j \approx \sum_i \lambda_{i,j} h_i / \alpha_j^{out}.$$

Two passenger demand parameters $x_j$ and $X_j$ have been introduced in~\cite{SFLG18_acc}.
\begin{equation}
x_j = \left(\frac{\lambda_j^{out}}{\alpha_j^{out}} + \frac{\lambda_j^{in}}{\alpha_j^{in}}\right),
\label{eq-6}
\end{equation}
in such a way that $x_j h_j$ gives the time needed for passenger download and upload at platform~$j$,
\begin{equation}
  X_j = \frac{x_j}{1 - x_j},
\end{equation}
in such a way that $x_j h_j = X_j g_j, \forall \underline{h}_j \leq h_j \leq \bar{h}_j$ and 
$\underline{g}_j \leq g_j \leq \bar{g}_j$, and $\forall j$.

\subsection{The dwell time control law}

The dwell time model is the following.
\begin{align}
  & w_j^k(h_j^k,x_j) = \min (x_j h_j^k, \bar{w}_j ), \label{eq-dwell}
\end{align}
where the $k^{th}$ dwell time on platform $j$ is the minimum of the needed time for passenger upload and download
on the same platform, and a maximum dwell time $\bar{w}_j$ to avoid congestion behind a train stopping too long at a platform. 

\subsection{The running time control law}

We propose the following running time law.
\begin{align}
  & r_j^k(h_j^k,x_j) = \max \left\{ \underline{r}_j, \tilde{r}_j - x_j \left( h_j^k - \underline{h}_j \right) \right\}, \label{eq-run}
\end{align}
where $\tilde{r}_j$ is the average (nominal) running time of trains on section $j$.

The model~(\ref{eq-run}) gives the running time as the maximum between a given minimum running time $\underline{r}_j$
and a term that subtracts $x_j \left( h_j^k - \underline{h}_j \right)$ from the nominal running time.
The term $x_j \left( h_j^k - \underline{h}_j \right)$ expresses a deviation of the upload and download time, 
due to a deviation of the train time-headway.
We notice here that the term $x_j h^k_j$ appearing in the dwell time law~(\ref{eq-dwell}) with a sign ``$+$'',
appears in the running time law~(\ref{eq-run}) with a sign ``$-$''.

The train travel time law is then obtained by summing the dwell time law~(\ref{eq-dwell}) with the running time law~(\ref{eq-run}).
\begin{equation}\label{eq-travel}
   t_j^k(x_j) = r_j^k(h_j^k,x_j) + w_j^k(h_j^k,x_j).
\end{equation}

Let us use the notations.
\begin{align}
    & \Delta h_j := \bar{h}_j - \underline{h}_j, \quad \Delta g_j := \bar{g}_j - \underline{g}_j, \quad 
                   \Delta w_j := \bar{w}_j - \underline{w}_j, \quad \Delta r_j := \tilde{r}_j - \underline{r}_j. \nonumber
\end{align}    
It is then easy to check the following. $\Delta w_j = x_j \Delta h_j = X_j \Delta g_j, \forall j$.
Then we have the following result.
\begin{theorem}\cite{SFLG18_acc}\label{thm-mp}
  If $h^1_j \leq \bar{h}_j=1/(1 - x_j) \; \bar{g}_j, \forall j$ and $\Delta r_j \geq \Delta w_j = X_j \Delta g_j, \forall j$,
  then the dynamic system~(\ref{eq-d1}), where the train dwell and run times are controlled as in~(\ref{eq-dwell}) and~(\ref{eq-run}) respectively,
  is a Max-plus linear system, and is equivalent to
  \begin{equation}
     d^k_j = \max\left\{ d^{k-b_{j}}_{j-1} + \tilde{r}_j + X_j \underline{g}_j \;,\;
                           d^{k-\bar{b}_{j+1}}_{j+1} + \underline{s}_{j+1} \right\}. \label{eq-thm-mp}
  \end{equation}
\end{theorem}
\proof See~\cite{SFLG18_acc}. \endproof

\noindent
Let us interpret the two conditions of Theorem~\ref{thm-mp}.
\begin{itemize}
 \item Condition $h^1_j \leq \bar{h}_j=1/(1 - x_j) \; \bar{g}_j, \forall j,k$ limits the initial headway $h^1_j$
   (i.e. the initial condition) to 
   its upper bound $\bar{h}_j$, which is given by the level of the passenger travel demand $x_j$ at platform $j$,
   and by the upper bound $\bar{g}_j$ on the safe separation time at the same platform. This condition tells that 
   the perturbation of train dynamics at time zero is bounded.
 \item Condition $\Delta r_j \geq \Delta w_j = X_j \Delta g_j, \forall j$ limits the margin on the train dwell times
   to the margin on the train running times. This condition tells that every extension of the train dwell time at any platform
   due to the level of the passenger travel demand at this platform can be compensated by the available margin in
   train run time at this platform.
\end{itemize}

\subsection{The traffic phases}

Under the conditions of Theorem~\ref{thm-mp}, we have a Max-plus linear system, from which we derive analytically the 
traffic phases of the train dynamics, giving the average train time headway and frequency as functions of the number $m$
of running trains and of the passenger travel demand level $x$ (or equivalently $X$).
This derivation is similar to the one done in section~\ref{sec-mpm} of Chapter~\ref{chap-maxtrain} (Theorem~\ref{th-mpm}).

\begin{theorem}\label{thm-1b}
The average asymptotic train time-headway of the linear Max-plus system with dynamic dwell and run
times depending on the demand is given by the average asymptotic growth rate
of the system.
The average headway depends on the number of trains $m$ on the line and the passenger travel demand for every station $X_j$.
   $$  h(m,X) = \max \left\{ \frac{\sum_j (g^{\min}_j X_j + \tilde{r}_j)}{m} \;,\;
                             \max_j ((g^{\min}_j X_j + \tilde{r}_j)+\underline{s}_j) \;,\;
                               \frac{\sum_j \underline{s}_j}{n-m} \right\}.$$                               
\end{theorem}
\proof See~\cite{SFLG18_acc}. \endproof

\begin{corollary}\label{cor-1b}
The average frequency of the linear Max-plus system with demand-dependent dwell and run times is a function of the number of trains and the travel demand.
	$$f (m,X) = \max \left\{ 0 \;,\; \min \left\{\frac{m}{\sum_j (g^{\min}_j X_j + \tilde{r}_j)} \;,\;
					             \frac{1}{\max_j ((g^{\min}_j X_j + \tilde{r}_j)+\underline{s}_j)} \;,\;
					             \frac{n-m}{\sum_j \underline{s}_j} \right\} \right\}.$$
\end{corollary}
\proof Directly from Theorem~\ref{thm-1b} with $f = 1/h$ and $f \geq 0$. \endproof

Figures of Table~\ref{tab} illustrate the results of Theorem~\ref{thm-1b}) and Corollary~\ref{cor-1b}.
The parameters of the metro line are from RATP Paris (18 stations, with the same passenger demand at all platforms).
For the interpretation of the traffic phases, please see\cite{SFLG18_acc}.

\begin{table*}[thbp]
\centering
\caption{Average asymptotic train time-headway $h$ and frequency $f$ as functions of the total number $m$ of running trains 
     and of the passenger travel demand level $x$ (or equivalently $X$).
     Three traffic phases are distinguished. Linear metro line. Parameters from RATP Paris.}~~\\~~ 
\begin{tabular}{|c|c|}
  \hline
  $h(m,x)$ & $h(m,X)$ \\
  \hline
   & \\
  \includegraphics[scale=0.18]{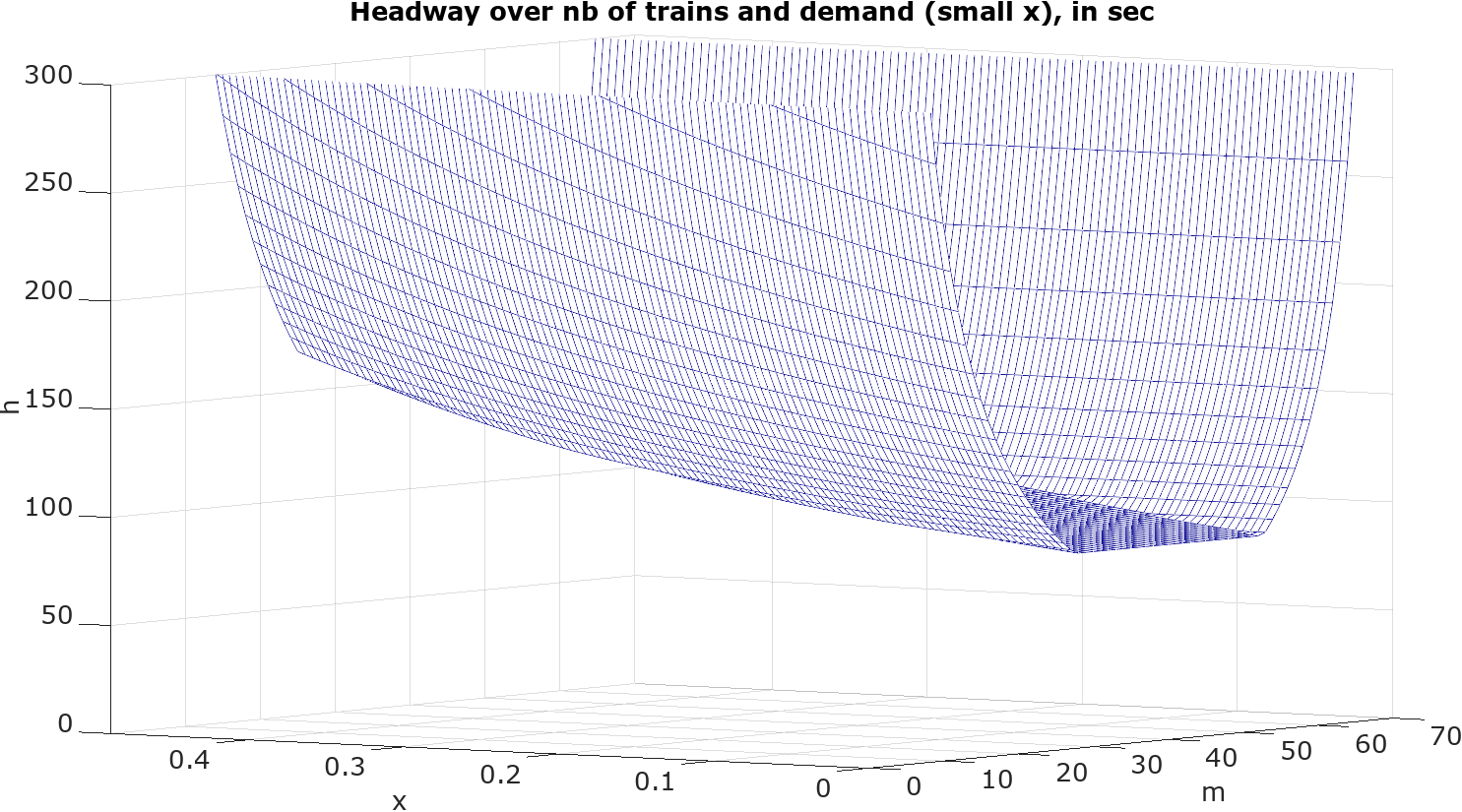} & \includegraphics[scale=0.18]{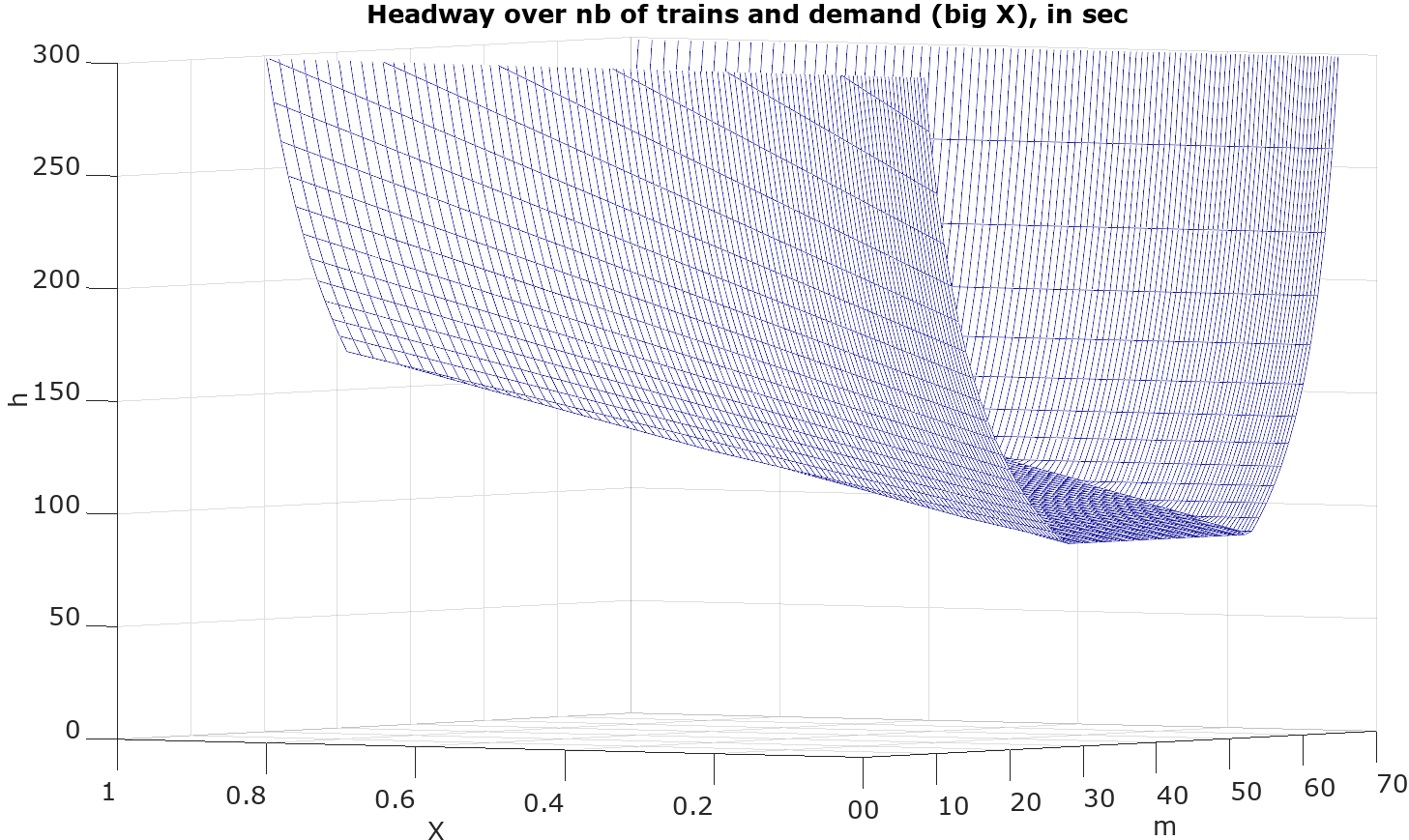} \\
  \hline 
  $f(m,x)$ & $f(m,X)$ \\
   \hline
   & \\
  \includegraphics[scale=0.18]{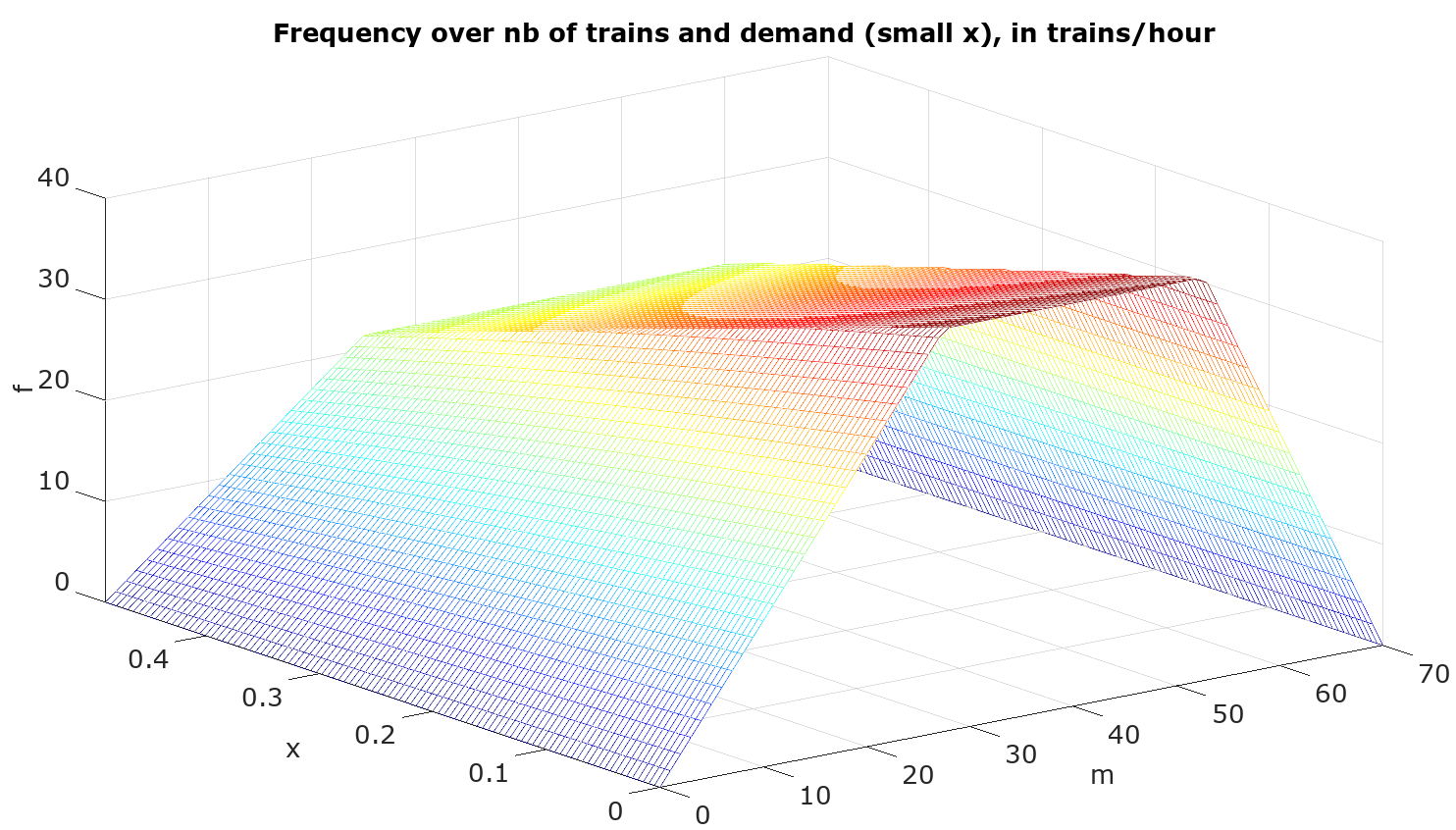} & \includegraphics[scale=0.18]{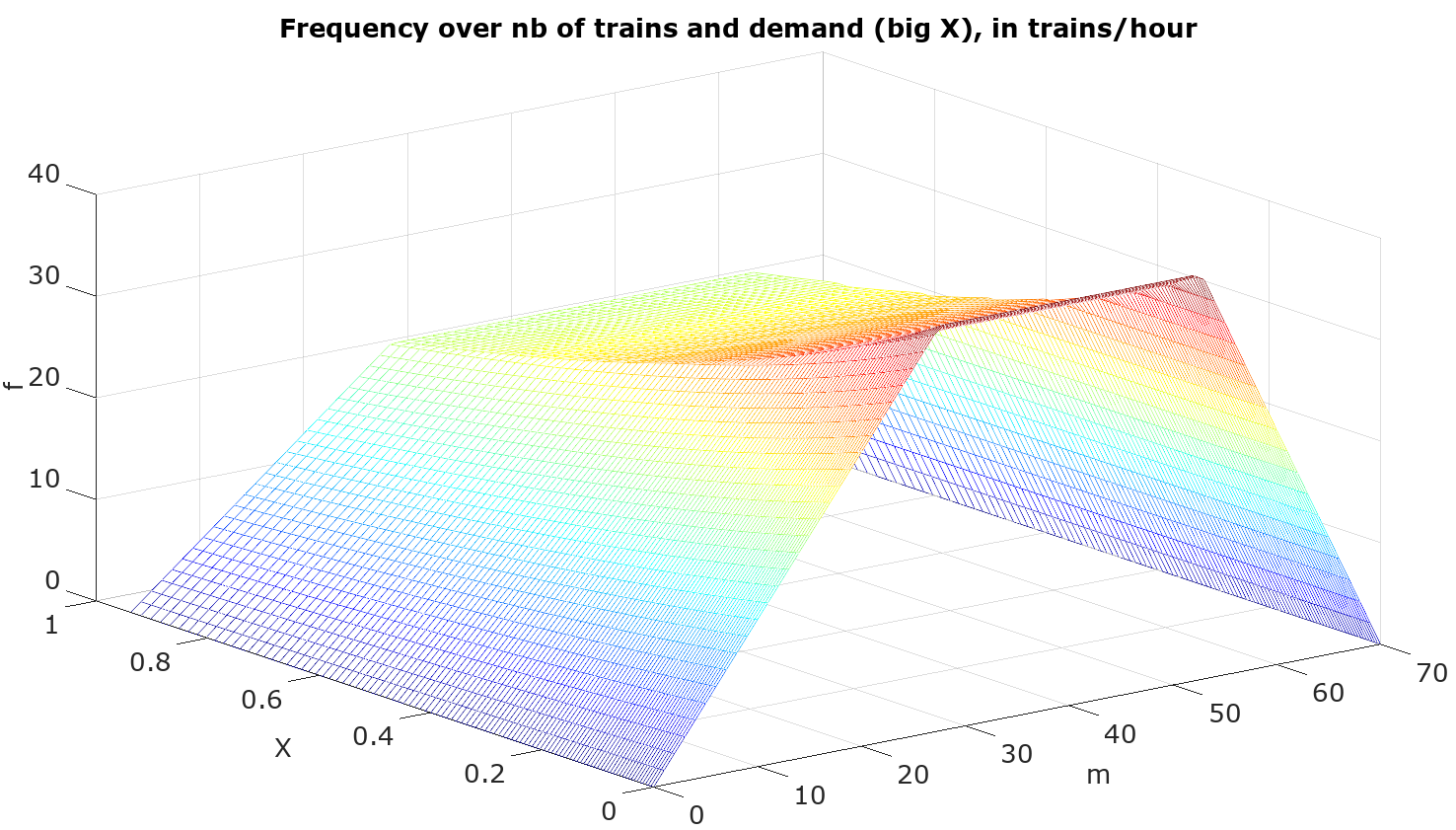} \\
  \hline 
  & \\
  \includegraphics[scale=0.18]{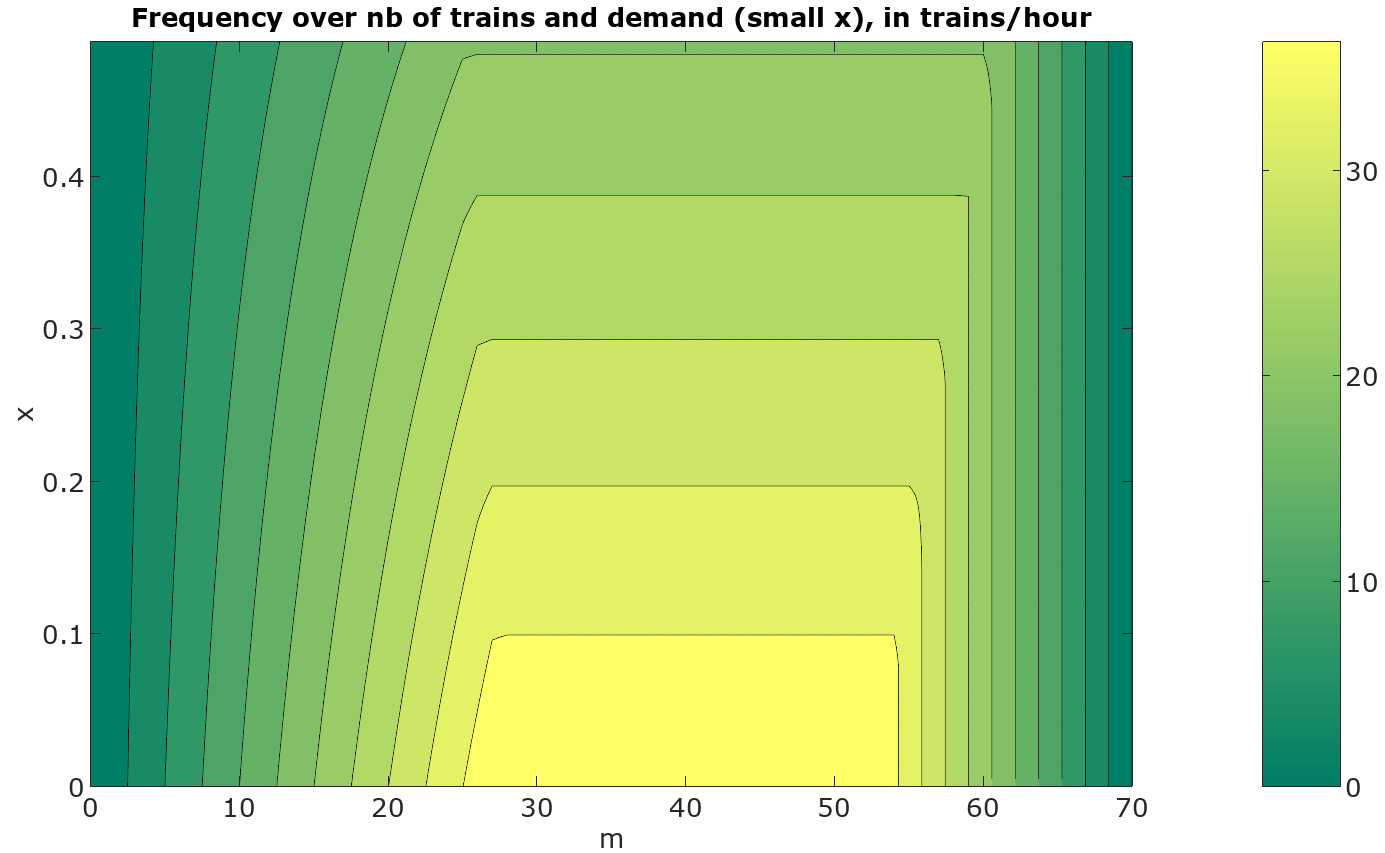} & \includegraphics[scale=0.18]{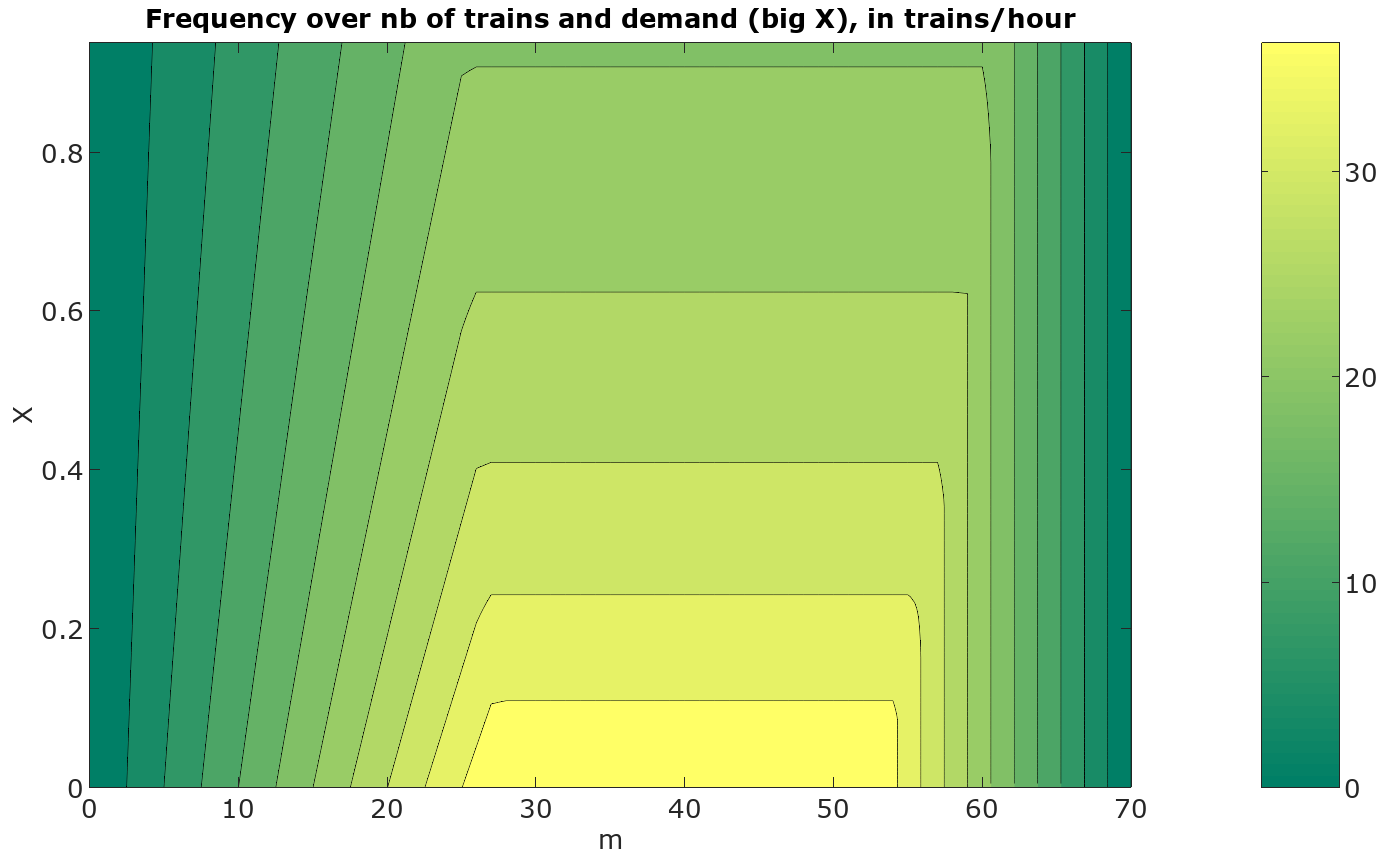} \\
  \hline
\end{tabular}
\label{tab}
\end{table*}


\newpage
This page is intentionally left blank

\chapter{Road network calculus}
\label{chap-roadnetcal}

The main references of this chapter are~\cite{FHL14,FHL14a}.
We propose a traffic system theory for the calculus of performance
bounds in road traffic networks.
The system theory is made in the min-plus algebra of 2$\times$2 matrices of functions~\cite{BCOQ92}.
The performance bounds are derived basing on basic and extended results of the
network calculus theory~\cite{Cruz91a,Cruz91b,Cha00,LT01}.
The road traffic is modeled with the well known 1st order model LWR (Lighthil Whiteham and Richards),
and following the well known cell transmission model scheme.

We first present the road section model, which will be the basis of our modeling here.
We write the traffic dynamics on one road section under the cell transmission
model, with given traffic demand and supply. 
We show that the dynamics can be written linearly in the Min-plus algebra of 2$\times$2 matrices of functions.
We derive the impulse response of that system, and show that the response can be interpreted
as a service guarantee of the traffic system seen as a server.

Second, we define a concatenation operator of two road section systems, and show that
a road link of many sections can be built algebraically by concatenating elementary
road section systems. We also define a feedback operator for this traffic system theory,
and use it to model closed traffic systems.

Finally, we show how this theory can be extended to two dimensional road traffic.
We present the model of an urban junction controlled with a traffic light, give
some insights on other two dimensional traffic control models, and on the highway traffic modeling,
and some directions for further extensions.

\begin{table}
\caption{Nomenclature.}~~\\~~
\begin{tabular}{ll}
$u(t)$	    & inflow at time t.\\
$U(t)$	    & cumulated inflow from time zero to time t.\\
$y(t)$	    & outflow at time t.\\
$Y(t)$	    & cumulated outflow from time zero to time t.\\
$\alpha$    & maximum arrival curve (time function).\\
$\beta$	    & minimum service curve (time function).
\end{tabular}\\~~\\~~\\
\begin{tabular}{ll}
$\mathcal F$	    & the set of non-decreasing time functions $f$ satisfying $f(t)=0,\forall t<0$.\\
$\oplus$    & element-wise operation (min-plus addition in $\mathcal F$), $(f\oplus g)(t)= \min (f(t),g(t))$.\\
$*$	    & min-plus convolution in $\mathcal F$, $(f*g)(t)= \inf_{0\leq s\leq t} (f(s)+g(t-s))$.\\
$\oslash$  & min-plus deconvolution in $\mathcal F$, $(f\oslash g)(t)= \sup_{s\geq 0} (f(t+s)-f(s))$.\\
$\underline{\oslash}$ & max-plus deconvolution in $\mathcal F$, $(f\underline{\oslash} g)(t)= \inf_{s\geq 0} (f(t+s)-f(s))$.\\
$\varepsilon$	& the zero element of the dioid $(\mathcal F,\oplus,*)$, $\varepsilon(t)=+\infty, \forall t\geq 0$.\\
$e$         & the unity element of the dioid $(\mathcal F,\oplus,*)$, $e(t)=+ \infty, \forall t>0$, and $e(0)=0$.\\
$f^k$	    & convolution power.  $f^k=f*f*…*f \;\;$   ($k$ times).
\end{tabular}\\~~\\~~\\
\begin{tabular}{ll}
$B(t)$      & the backlog at time $t$. $B(t)=U(t)-Y(t)$. \\
$d(t)$	    & the virtual delay at time $t$, $d(t)=\inf \{h\geq 0, Y(t+h) \geq U(t) \}$.\\
$\gamma^p$	& a particular function in $\mathcal F$, $\gamma^p(t)= + \infty\; \forall t>0$, and  $\gamma^p(0)=p$.\\
$\delta^T$	& a particular function in $\mathcal F$. $\delta^T (t)=0 \;\forall t\leq T$, and $\delta^T(t)=+\infty \forall t>T$.\\
$[f]^+$	    & $:= \max(0,f)$. \\
$\Lambda(r,s)$ & affine maximum arrival curve $\Lambda(r,s)(t) = rt+s$. \\
$\lambda(R,T)$ &  rate-latency minimum service curve $\lambda(R,T)(t) = R(t-T)^+$.
\end{tabular}\\~~\\~~\\
\begin{tabular}{ll}
$q$	        & car-flow. \\
$q_{\max}$	& maximum car-flow.\\
$q_i(t)$	& car outflow from the ith section at time t.\\
$Q_i(t)$	& cumulated car outflow from the $i$th section from time zero to time t.
\end{tabular}\\~~\\~~\\
\begin{tabular}{ll}                
$n_i$	    & number of cars on the $i$th section at time zero. \\
$n_{max}$	& maximum number of cars that a section can contain. \\
$\rho$      & average car-density. \\
$\rho_c$    & critical car-density. \\
$\rho_j$	& jam (or maximum) car-density. \\
$v$	        & free car speed. \\
$w$	        & backward wave speed. \\
$L$	& road section length. \\
$q(\rho)$	& the car-flow function of the car-density (fundamental traffic diagram). 
\end{tabular}\\~~\\~~\\
\begin{tabular}{ll}
$M$	    	& the set of $n\times m$ matrices with entries in $\mathcal F$. \\
$\varepsilon$ & the zero element for $(M_{(n\times n)},\oplus,*)$, \\
$e$		& the unity element for $(M_{(n\times n)},\oplus,*)$, \\
$A^k$		& power operation in $M_{(n\times n)}$. \\ 
$A^*$		& sub-additive closure in $M_{(n\times n)}$, such that $A^*=\bigoplus_(k\geq 0)(A^k)$.
\end{tabular}
\label{tab-rnc}
\end{table}

\section{Complements - Network calculus for MIMO systems}
\label{sec-cnc}

We give in this section some complements on the Multiple-Input Multiple-Output (MIMO) systems.
We are concerned here by the multi-dimensional case, 
where multiple inflows arrive to- and depart from the system.
The particular signals $\varepsilon, e , \gamma^p$ and $\delta^T$,
defined in Table~\ref{tab-rnc} will be  used here.

Let us consider a MIMO server with $n$ (cumulated) arrival flows $U_i, i=1,2,\ldots,n$, and $n$ (cumulated) departure flows $Y_i, i=1,2,\ldots,n$.
We recall here that maximum and minimum arrival curves $\alpha_i$ and $\underline{\alpha}_i$ resp. can be estimated from
the flow $U_i$ as follows.
\begin{align}
  & U_i \oslash U_i \leq \alpha_i \\
  & U_i \underline{\oslash} U_i \geq \underline{\alpha_i}
\end{align}

In a system where the service depend on multiple arrival flows, it is not sufficient to
upper bound each arrival separately, but one needs also to bound the dependencies of the arrivals on each other.
In the following we show that the extension of the definition of arrival curves in the case of one arrival,
to arrival matrices of curves for the case of multiple arrivals, includes naturally the dependencies
of the arrivals on each other. Let us first introduce a new notion of \textit{Time-shift arrival matrix},
which we will use to define an \textit{arrival matrix}.

\subsection*{Time-shift arrival matrix}

For a system with $n$ input and $n$ output flows, we need to estimate a matrix arrival $\alpha$, which is
a matrix of arrival curves such that $\alpha_{ii}, 1\leq i\leq n$ are nothing but maximum arrival curves
for arrival flows $U_i, 1\leq i\leq n$, while $\alpha_{ij}, i\neq j$ are curves that bound the deviations
between each pair of arrival flows $(i,j), i\neq j$.

For $i \neq j$, the difference with respect to the case $i=j$ is that, it is possible to have
$U_i(t)-U_j(s)>0$, even for $t < s$. Indeed, if we assume that $U_i(t)-U_j(s) \leq 0, \forall t<s$, then
we get $U_i(s)-U_j(s) \leq 0,\forall s \geq 0$, and similarly $U_j(s)-U_i(s) \leq 0, \forall s \geq 0$.
Therefore, $U_i(s)-U_j(s) = 0, \forall s \geq 0$. It is trivial that such an assumption is very restrictive.
Therefore, if we like to upper bound $U_i(t)-U_j(s)$ for all $s,t \geq 0$, then we need to work with
negative times for the arrival curves. In order to keep working with non-negative times, we back-shift
the curves with negative times to zero.
To obtain such back-shifted arrival matrices, we first define what we call here \textit{time-shift matrix} $T$
(of non negative entries).

\begin{definition}(Time-shift matrix)\label{shift}
  The time-shift matrix $T \in \mathbb R_+^{n\times n}$ for arrival flows $U_i, i=1,2,\ldots,n$ is defined as follows
  $$T_{ij} = \sup_{t\geq 0} \inf \{s\geq 0, U_i(t+s)-U_j(t) \geq 0\}, \quad \forall i,j.$$
\end{definition}

It is easy to see that $T_{ii}=0, \forall i=1,2,\ldots,n$.

For two arrival flows $U_i$ and $U_j$, if we  see $U_i$ as an output flow of $U_j$, then $T_{ij}$ can be seen as the maximum
virtual delay (see definition of virtual delay in Table~\ref{tab-rnc}). From that remark, an easy way to estimate the matrix $T$ is given b Proposition~\ref{prop0} below.

\begin{proposition}\label{prop0}
  If $\alpha_i$ and $\underline{\alpha}_i$ are resp. maximum and minimum arrival curves for $U_i, i=1,2,\ldots,n$, and if $T$ is
  a time-shift matrix for $U_i, i=1,2,\ldots,n$, then
  $$T_{ij} \leq \sup_{t\geq 0} \inf \{h\geq 0, \underline{\alpha}_i(t+h) - \alpha_j(t) \geq 0  \}.$$
\end{proposition}

\proof Follows directly from the definitions of $T, \alpha$ and $\underline{\alpha}$. \endproof

In order that both elements $T_{ij}$ and $T_{ji}$ of the shift matrix $T$ be finite, we need to guarantee that
$$\lim_{t\to\infty} U_i(t)/t \leq \lim_{t\to\infty} U_j(t)/t, \quad \text{ and } \quad \lim_{t\to\infty} U_j(t)/t \leq \lim_{t\to\infty} U_i(t)/t.$$
Hence, we need that
$$\lim_{t\to\infty} U_i(t)/t = \lim_{t\to\infty} U_j(t)/t,$$
which means that the two flows have the same asymptotic average rate (they are comparable in average).

\begin{definition}(Arrival matrix)\label{def-arrival}
   For a given $n\times 1$ vector $U$ of cumulated arrival flows $U_i, i=1,\ldots, n$, a matrix
   $\alpha \in \mathcal F^{n\times n}$ is said to be a maximum (resp. minimum) arrival matrix for $U$ if    
   there exists a time-shift matrix $T$ for $U_i, i=1,2,\ldots, n$, such that
   $$ \forall i,j=1,2,\ldots ,n, \forall s,t \in \mathbb N, U_i(t)-U_j(s) \leq \alpha_{ij}(T_{ij}+t-s).$$
   $$(resp. 
            U_i(t)-U_j(s) \geq \alpha_{ij}(T_{ij}+t-s).)$$
   which can also be written
   $$\forall i,j = 1,2,\ldots ,n, U_i \leq \delta^{-T_{ij}} (\alpha_{ij}*U_j).$$
   $$(resp. 
            U_i \geq \delta^{-T_{ij}} (\alpha_{ij}*U_j)).$$
\end{definition}

Let us notice that Definition~\ref{def-arrival} is different form Definition~4.2.1 given in~\cite{Cha00}.
Definition~\ref{def-arrival} is illustrated in Example~\ref{ex1} below.
A procedure for estimating an arrival matrix for a cumulated inflow vector $U$, from the data of $U$ itself, is the following.
\begin{enumerate}
 \item Estimate maximum and minimum arrival curves for each arrival flow, individually, by the deconvolution.
   $$\begin{array}{ll}
        \alpha_i = \alpha_{ii} = U_i \oslash U_i, \forall i.\\
        \underline{\alpha}_i = \underline{\alpha}_{ii} = U_i \underline{ \oslash } U_i, \forall i.
      \end{array}$$
 \item Estimate a shift-time matrix from those curves, by maximum delay calculus, as explained above in Proposition~\ref{prop0}.
   $$T_{ij} \leq \max_{t\geq 0} \min \{s\geq 0, \underline{\alpha}_i(t+s) - \alpha_j(t) \geq 0  \}.$$
 \item Estimate non-diagonal shift arrival curves, by shift deconvolution. 
    From Definition~\ref{def-arrival}, $\alpha_{ij}$ satisfy
    $$\alpha_{ij} \geq \delta^{T_{ij}} (U_i \oslash U_j).$$
    It is easy to check that for $i=j$, we have $T_{ii}=0$, and then $\alpha_{ii}$ is a one-dimensional arrival curve
\end{enumerate}

\begin{figure*}[h]
    \centering
    \includegraphics[width=1\textwidth]{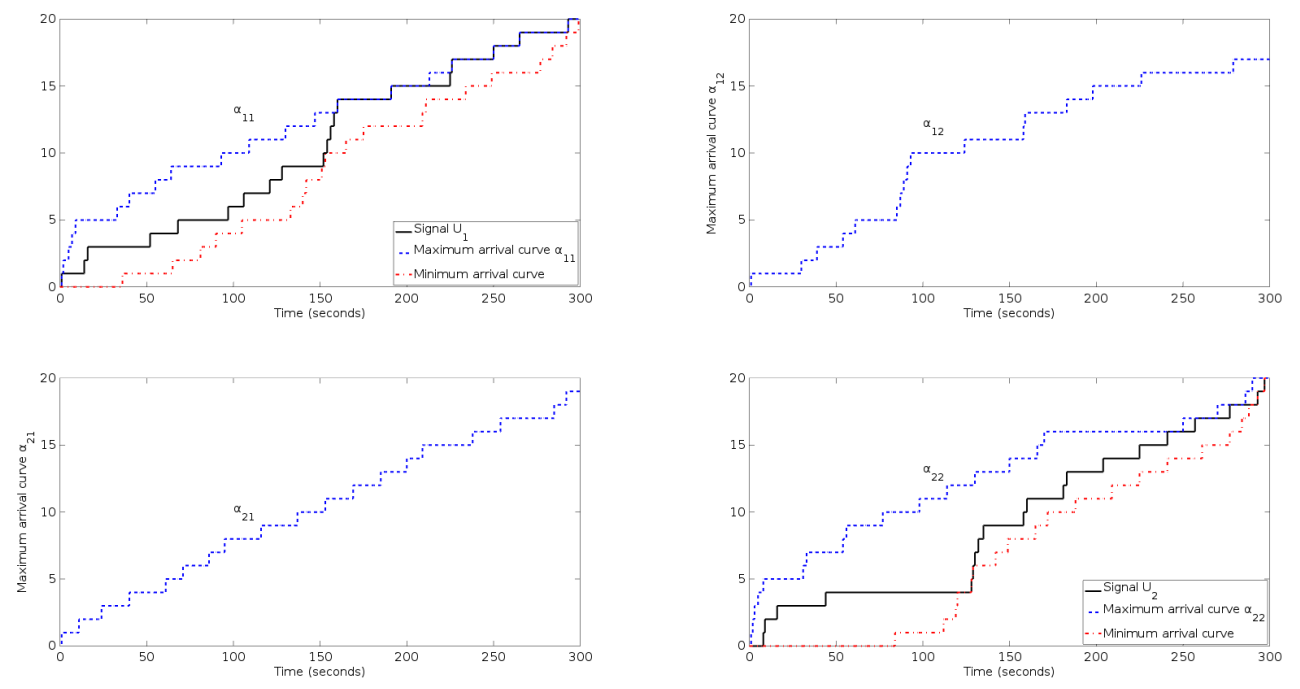}
    \caption{The four curves of the arrival matrix of a vector of two signals.}
    \label{arrival_matrix}
\end{figure*}

\begin{example}\label{ex1}
In Figure~\ref{arrival_matrix} we show the four curves of an arrival matrix (estimated as explained above) for a vector of two signals $U_1$ and $U_2$.
Those signal are supposed given from time zero to time 300. The calculus of those curves is done as follows.
\begin{itemize}
  \item $\alpha_{11} = \bar{\alpha}_1 = U_1 \oslash U_1$.
  \item $\alpha_{22} = \bar{\alpha}_2 = U_2 \oslash U_2$.
  \item $T_{12} = \max_{t\geq 0} \min\{ s\geq 0, U_1(t+s) - U_2(t) \geq 0\} = 60$.
  \item $T_{21} = \max_{t\geq 0} \min\{ s\geq 0, U_2(t+s) - U_1(t) \geq 0\} = 8$.  
  \item $\alpha_{12} = \delta^{T_{12}} (U_1\oslash U_2)$.
  \item $\alpha_{21} = \delta^{T_{21}} (U_2\oslash U_1)$.
\end{itemize}
\end{example}

Let us now define the service for a MIMO server.
We base here on an existing definition of service guarantee in multi-dimensional case (multiple input flows and multiple output flows)
given in~\cite{Cha00}. Then by combining that definition with Definition~\ref{def-arrival} given above,
we recalculate upper bounds for the virtual delays as well as other performance bounds.

\begin{definition} (Service matrix)~\cite{Cha00}\label{def-sm}
   For a given server with input vector $U$ and output vector $Y$, a $n\times n$ matrix $\beta$ is said to be
   a service matrix for the server, if $Y \geq \beta * U$.
\end{definition}
   
\begin{definition} (Virtual delay)~\cite{Cha00}
   For a given server with input vector $U$ and output vector $Y$, the virtual delay of the last quantity arrived
   at time $t$ from the $i$th input to depart from the $i$th output, denoted $d_i(t)$ is defined:
   $$d_i(t)= \inf\{d\geq 0, Y_i(t+d) \geq U_i(t)\}.$$
\end{definition}
   
The following result improves Theorem~4.3.6 of~\cite{Cha00}, which derives upper bounds for the virtual delays
through input-output pairs.
\begin{theorem}\label{th-vec-del}
  For a given server with input vector $U$ and output vector $Y$, if $\alpha$ is a shift-arrival matrix with $T$ a shift matrix, for $U$,
  and if $\beta$ is a service matrix for the server, then $\forall i=1,2,\ldots ,n, \forall t\in \mathbb N $, 
  $$ d_i(t) \leq \inf\{d\geq 0, \alpha_{ij}(T_{ij}+s)
     \leq \beta_{ij}(s+d), -T_{ij} \leq s \leq t, \forall j\}.$$
  and then the virtual delays $d_i, i=1,\ldots, n$ are bounded as follows. $\forall i=1,2,\ldots,n, \forall t\in \mathbb N$,
  $$d_i(t) \leq \max_{1\leq j\leq n} \sup_{-T_{ij} \leq s \leq t} 
     \inf \{ d\geq 0, \alpha_{ij}(T_{ij}+s) \leq \beta_{ij}(s+d)\}.$$
  or equivalently
  $$d_i(t) \leq \max_{1\leq j\leq n} \left\{ T_{ij} + \sup_{s\geq 0}
     \inf \{ d\geq 0, \alpha_{ij}(s) \leq \beta_{ij}(s+d)\} \right\}.$$
\end{theorem}
  
\proof. The proof is an adaptation of the proof of Theorem~4.3.6 in~\cite{Cha00}.
\endproof

In order to build the traffic system corresponding to a whole traffic network, we base in elementary traffic
systems which will be used as unit
systems in the composition. We consider here two elementary systems (an uncontrolled and a controlled road
sections). The composition (concatenation) we use here are inspired from~\cite{Far08};
see also~\cite{Far09, FGQ11, FGQ05, FGQ07, FGQ11b, Far12, FHL14}.
It consists in two operators: a concatenation one and a feedback one. More details are given in
Section~\ref{sec-rsm} below.

\section{The road section model}
\label{sec-rsm}

We give in this section the derivation of a service matrix for a road section traffic system.
The model gives the traffic dynamics on which this approach is based (cell transmission model~\cite{Dag94}),
and by that explain why the network calculus theory can be applied to road traffic networks.

\begin{figure}[h]
  \begin{center}
    \includegraphics[width=5cm]{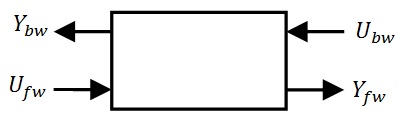}    
    \caption{The traffic system of a road section.}
  \end{center}
  \label{fig-section}
\end{figure}

In a stretch of road subdivided in $n$ sections $1,2,\ldots,i-1,i,i+1,\ldots,n$, we consider 
here one road section, section $i$, as the basic traffic system.
We denote by $L$ the length of the road section. 
Cars arrive from the left side of the road section (input flow), pass through it, and departs from
the right side of it (output flow). The inputs $U_{fw}$ and $U_{bw}$ represent respectively the traffic demand from the upstream
section~$i-1$ to section~$i$, and the traffic supply of the downstream section~$i+1$ to section~$i$. 
The outputs $Y_{fw}$ and $Y_{bw}$ represent respectively the traffic outflow from section~$i$ to the downstream
section~$i+1$ (which is also the traffic demand of section $i+1$), and the traffic supply of section~$i$ to the
upstream section~$i-1$; see Figure~\ref{fig-section}.

Let us clarify the notations $U_{fw}, Y_{fw}, U_{bw}$ and $Y_{bw}$:
\begin{itemize}
	\item $U_{fw}(t)$ : cumulated forward inflow of cars from time zero to time $t$.
	\item $Y_{fw}(t)$ : cumulated forward outflow of cars from time zero to time $t$.
	\item $U_{bw}(t)$ : cumulated backward supply of section $i+1$ from time zero to time $t$.
	\item $Y_{bw}(t)$ : cumulated backward supply of section $i$ from time zero to time $t$.
\end{itemize}

A number of assumptions is considered here.
First, we assume given an initial number $n$ of cars on the road section a time zero.
Second, we suppose that cars pass through the road section under the FIFO (first in first out) rule.
Third, we assume that the road section cannot contain more than $n_{\max}$ cars in a given time.
We denote by $\bar{n} = n_{\max} - n$ the additional number of cars that the road section could contain at time zero.
The number $\bar{n}$ corresponds to the free spaces in the section at time zero.

In order to simplify the presentation of the ideas, and without loss of generality, we assume that all the cumulated flows are initialized to
zero~\footnote{Indeed, we only need that $U_{fw}=Y_{fw}$ and $U_{bw}=Y_{bw}$.}.
\begin{equation}\label{eq-causality}
  U_{fw}(0) = Y_{fw}(0) = U_{bw}(0) = Y_{bw}(0) = 0.		
\end{equation}	

\begin{figure}[h]
  \begin{center}
    \includegraphics[width=10cm]{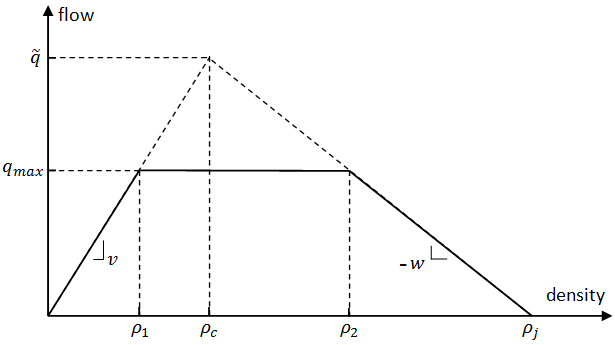}    
    \caption{A trapezoidal fundamental diagram; with $\rho_j = n_{\max}/L, \rho_1 = q_{\max}/v, \rho_2 = \rho_j - q_{\max}/w, 
          \tilde{q} = \rho_j / (1/v + 1/w), \rho_c = \tilde{q}/v = \rho_j - \tilde{q}/w$.}
  \end{center}
  \label{fig-trapez}
\end{figure} 

Let us now write the traffic dynamics on the road section. As in~\cite{FHL13}, we base on the cell-transmission model~\cite{Dag94} with a
trapezoidal fundamental diagram; see Figure~\ref{fig-trapez}. $v, q_{\max}$ and $w$ denote here the free flow speed, the capacity (maximum)
flow, and the backward wave speed for the road section.
We obtain the following dynamics, where we introduce an intermediate variable $Q$, which is simply the
cumulated forward outflow $Y_{fw}$.

\begin{align}
  & Q(t) = \min \left\{ \begin{array}{l}
			    U_{fw}\left(t-L/v \right)+n, \\
			    Q\left(t-L/v\right)+q_{max}L/v, \\
			    U_{bw}(t)
                        \end{array} \right., \label{dyn1a}\\
  & Y_{fw}(t) = Q(t), \label{dyn2a} \\
  & Y_{bw}(t) = Q\left(t-L/w+\bar{n}\right). \label{dyn3a}
\end{align}

We will see here the variables $U, Y$, and $Q$ as time signals in the dioïd $(\mathcal F,\oplus,*)$ defined above, where the addition of two signals
is the point-wise minimum, and where the multiplication of two signals is the minimum convolution of the signals.
By using the min-plus algebra notations (see~\cite{FHL13},~\cite{BCOQ92}), we get:
\begin{align}
  & Q = f_3 Q \oplus f_1 U_{fw} \oplus U_{bw}, \label{sys1} \\
  & Y_{fw} = Q \oplus e, \label{sys2}\\
  & Y_{bw} = f_2 Q \oplus e, \label{sys3}
\end{align}
where $f_1 = \gamma^{n}\delta^{L/v}, f_2 = \gamma^{\bar{n}}\delta^{L/w}$ and $f_3 = \gamma^{q_{\max}L/v}\delta^{L/v}$.

\begin{figure}[h]
  \centering
  \includegraphics[width=11cm]{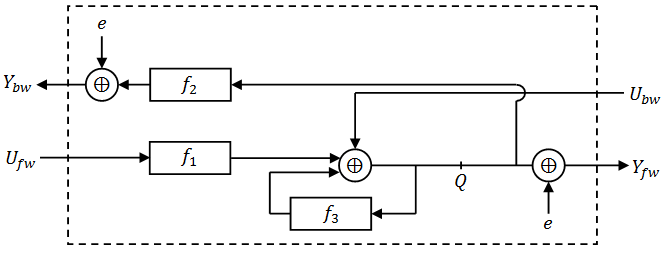}
  \caption{Diagram for the system~(\ref{sys1})-(\ref{sys3})}.
  \label{diag-sys}
\end{figure} 

where we added (min-plus addition) the unity vector $e$ to $Y_{fw}$ and to $Y_{bw}$ in order to satisfy condition~\ref{eq-causality}.
Then, by denoting by $U$ and $Y$ the two column vectors $U = (U_{fw} \; , \; U_{bw})$, and $Y = (Y_{fw} \; , \; Y_{bw})$, we can write
\begin{equation}\label{eq-matrix-dyn}
  \begin{array}{l}
     Q = A Q \oplus B U \\
     Y = C Q \oplus e
   \end{array}
\end{equation}
with $A = \gamma^{q_{\max} L/v} \delta^{L/v}$, $B$ denotes the line vector B$ = (\gamma^n \delta^{L/v} \; , \; e)$, 
$C$ denotes the column vector $C = (e \; , \; \gamma^{\bar{n}}\delta^{L/w})$, and with $e$ denoting in~(\ref{eq-matrix-dyn}) the
column vector $(e \; , \; e)$.

Therefore, the traffic dynamics on a road section may be represented with three matrices: a signal $A$, a line matrix $B$ of two signals,
and a column matrix $C$ of two signals. 
In a such configuration, $B$ represents the action on the traffic demand of the section, $C$ gives action on the traffic
supply that the section offers for an eventual upstream section, while $A$ models the outflow limit $q_{\max}$ imposed by the trapezoidal
fundamental diagram.

From Theorem\ref{th-linsys}, the maximum sub-solution of~(\ref{eq-matrix-dyn}) on the variable $Q$ is $Q = A^* * B * U$. Moreover, we have
$Q \geq A^* * B * U$. Therefore, we obtain $Y \geq (C*A^**B)*U \oplus e$.

Then, since $U(0)=0$, we have $e \geq U$. Hence
\begin{equation}\label{eq-impulse}
  Y \geq (C*A^**B)*U\oplus U = (e\oplus C*A^**B)*U.  
\end{equation}
From~(\ref{eq-impulse}) we obtain a bound for the impulse response of the min-plus linear system~(\ref{eq-matrix-dyn}).
Moreover, the formula~(\ref{eq-impulse}) tells that $e\oplus C*A^* * B$ is 
a service matrix of the traffic system of one road section, seen as a server; see Definition~\ref{def-sm}.

\begin{theorem}\label{thm-impulse}
  A service matrix of the traffic system~(\ref{dyn1a})-(\ref{dyn3a}) is
  $\beta = e\oplus C*A^**B$, also written
  $$\beta = e \oplus (\gamma^{q_{\max} L/v} \delta^{L/v})^* 
      * \begin{pmatrix}
			\gamma^n \delta^{L/v} & e\\
			\gamma^{n+\bar{n}} \delta^{L/v + L/w} & \gamma^{\bar{n}} \delta^{L/w}
        \end{pmatrix}.$$
\end{theorem}

\proof It follows directly from~(\ref{eq-impulse}). \endproof

\begin{corollary}
   The service matrix $\beta$ for the traffic system~(\ref{eq-matrix-dyn}) satisfies
   \begin{align}
      & \beta_{11} \geq  \Lambda(q_{\max}, L/v) + n, \nonumber \\
      & \beta_{12} \geq \Lambda(q_{\max}, L/v), \nonumber \\
      & \beta_{21} \geq \Lambda(q_{\max}, L/v+L/w) + n + \bar{n}, \nonumber \\
      & \beta_{22} \geq \Lambda(q_{\max}, L/w) + \bar{n}. \nonumber
   \end{align}
   \label{cor-serv}
\end{corollary}

\proof It follows directly from Theorem~\ref{thm-impulse}. \endproof

We notice here that the $n$ vehicles being on the road section at time zero are counted in the
output flow $Y_{fw}$. Similarly, the $\bar{n}$ free spaces available in the road section at time zero 
are counted in the backward output flow (traffic supply flow) $Y_{bw}$.
By consequent, for the calculus of upper bounds (travel times (end-to-end delays), unserved demand (backlogs), etc.) for the road section,
we have to consider augmented arrival flows $U_{fw} + n$ and $U_{bw} + \bar{n}$ rather than $U_{fw}$
and $U_{bw}$ respectively. Therefore, we need to calculate a matrix arrival for the vector 
$(U_{fw}+n,U_{bw}+\bar{n})$ of arrival flows.

We notice also that for the calculus of an upper bound for the delay (travel time in case of a road traffic system)
for $U$ without taking into account the $n$ particles being in the system at time zero, we need to take the maximum
of the horizontal distances between the curves $\alpha_{ij}$ and $\beta_{ij}$ by considering only the horizontal segments
above $n$ and above $\bar{n}$. Moreover, the upper bound on the backlog (unsatisfied demand) includes
the initial $n$ particles.

\begin{example}\label{ex-lambda}
  If we consider a road section of length $L = 200$ m. (meters), with fundamental diagram parameters $v = 28$ m/s (about 100 km/h),
  $w = 7$ m/s (about 25 km/h), $\rho_j = 1/10$ veh/m (100 veh/km), and $q_{\max} = 0.5$ veh/s (1800 veh/h), 
  and if we consider a time step $dt = 5$ s, and an initial number of vehicles in the 
  road $n = 10$ vehicles (i.e. $\rho = 1/20 $ veh/m)  
  (and thus $\bar{n} = 10$, since $n_{\max} = 20$), then 
  $\rho_1=1/50$ veh/m and $\rho_2 = 1/35$ veh/m),
  and the service matrix $\beta$ is bounded by Corollary~\ref{cor-serv}
  as follows (the time unity being a second).
  \begin{itemize}
     \item $\beta_{11}(t) \geq  0.5(t - 7.14)^+ + 10$,
     \item $\beta_{12}(t) \geq 0.5(t - 7.14)^+$,
     \item $\beta_{21}(t) \geq 0.5(t - 35.71)^+ + 20$,
     \item $\beta_{22}(t) \geq 0.5(t - 28.57)^+ + 10$,
  \end{itemize}
  More precisely, the curves of the service matrix of the road section are 
  given by Theorem~\ref{thm-impulse}. 
\end{example}

\section{Composition of 1D traffic systems}

We propose here two operators for the connection of one-dimensional road traffic systems (traffic systems without junctions).
The first operator is for \textit{concatenating} two road traffic systems of two input flows and two output flows.
The resulting traffic system has also two input and two output flows.
The concatenation operator permits to build more complicated systems than a road section, such as 
a whole road of many sections.
The second operator is the \textit{feedback} operator. It permits to model closed traffic systems, such as
a circular (or ring) road.

\subsection{Concatenation}

\begin{figure}[h]
  \begin{center}
    \includegraphics[width=0.8\textwidth]{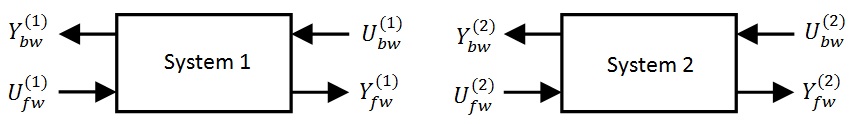}
    \caption{Composition of two min-plus linear traffic systems.}  
    \label{fig2}
  \end{center} 
\end{figure}

The composition of traffic systems is done in two dimensions, since each system has two inputs and two outputs.
The connection is not in series, in the sense that connecting two systems does not mean connect the outputs of
one system to the inputs of the other. Indeed, the connection is made here in the two directions, by connecting
an output of system~1 to an input of system~2, and an output of system~2 to an input of system~1. 
In Figure~\ref{fig2}, we illustrate the connection of two elementary systems (road sections).

Let us consider two min-plus linear traffic systems~1 and~2, with service curve matrices $\beta^{(1)}$ and $\beta^{(2)}$. We then have:
\begin{equation}\label{}
   \begin{pmatrix} 
     Y^{(i)}_{fw} \\
     Y^{(i)}_{bw}
   \end{pmatrix}
   \geq
   \begin{pmatrix}
     (\beta^{(i)})_{11} & (\beta^{(i)})_{12} \\
     (\beta^{(i)})_{21} & (\beta^{(i)})_{22}
   \end{pmatrix}     
   \begin{pmatrix}
     U^{(i)}_{fw} \\
     U^{(i)}_{bw}     
   \end{pmatrix}.
\end{equation} 
The resulting system from the concatenation has two input signals $U^{(1)}_{fw}$ and $U^{(2)}_{bw}$, and two output signals 
$Y^{(2)}_{fw}$ and $Y^{(1)}_{bw}$; see Figure~\ref{fig2}.
We have the following result.

\begin{theorem}\label{thm-connect}
  	A service matrix $\beta$ for the whole system is given by:
  	$$\begin{array}{l}
  	  \beta_{11} = \beta^{(2)}_{11} \beta^{(1)}_{11} \oplus 
  	                   \beta^{(2)}_{11} \beta^{(1)}_{12} \left( \beta^{(2)}_{21} \beta^{(1)}_{12} \right)^*
  	                     \beta^{(2)}_{21} \beta^{(1)}_{11} \\
	  \beta_{12} = \beta^{(2)}_{11} \beta^{(1)}_{12} \left( \beta^{(2)}_{21} \beta^{(1)}_{12} \right)^*\beta^{(2)}_{22} 
	                   \oplus \beta^{(2)}_{12} \\
	  \beta_{21} = \beta^{(1)}_{21} 
	                  \oplus \beta^{(1)}_{22} \left( \beta^{(2)}_{21} \beta^{(1)}_{12} \right)^* \beta^{(2)}_{21}\beta^{(1)}_{11} \\
	  \beta_{22} = \beta^{(1)}_{22}\left( \beta^{(2)}_{21} \beta^{(1)}_{12} \right)^* \beta^{(2)}_{22}.
   	\end{array}$$
    such that
  $$
   \begin{pmatrix} 
     Y^{(2)}_{fw} \\
     Y^{(1)}_{bw}
   \end{pmatrix}
   \geq 
   \begin{pmatrix}
     \beta_{11} & \beta_{12} \\
     \beta_{21} & \beta_{22}
   \end{pmatrix}     
   \begin{pmatrix}
     U^{(1)}_{fw} \\
     U^{(2)}_{bw}     
   \end{pmatrix}.
   $$    	
\end{theorem}
\proof The proof is available in~\cite{FHL14}. \endproof

As mentioned above, the immediate example for the concatenation of road traffic systems is
the construction of a whole road with a given number of sections, eventually with different lengths and characteristics (fundamental diagrams).
Service matrices of the road sections can be fixed based on the characteristics of the sections and by means of
Theorem~\ref{thm-impulse} and Corollary~\ref{cor-serv}.
The concatenation permits us to have a service matrix for the whole road by means of Theorem~\ref{thm-connect}.
It suffices then to upper bound the arrival flows of the system (forward demand $U_{fw}$ and backward supply $U_{bw}$) by
calculating an arrival matrix as explained above in section~\ref{sec-cnc}, to be able to derive deterministic upper bounds
for the travel time through the road, the unsatisfied traffic demand, the maximum car density, etc.

\begin{remark}\label{rem-road}
  An alternative way to build a road of many sections is to consider the whole dynamics of the road, as done in
  section~\ref{sec-rsm} (dynamics~(\ref{dyn1a}-\ref{dyn3a})),
  and derive the service matrix directly from the dynamics (as done by Theorem~\ref{thm-impulse}).
\end{remark}

\subsection{Feedback connection}

The objective here is to build closed traffic systems. We define a \textit{feedback} operator that permits
to link the input flows (or some of them) of a given traffic system with its output flows (or some of them).
One then obtain a closed traffic system. Such systems can be interesting for academic example analyses, and for
the derivation of qualitative characteristics.

\begin{figure}[h]
    \centering
    \includegraphics[width=0.5\textwidth]{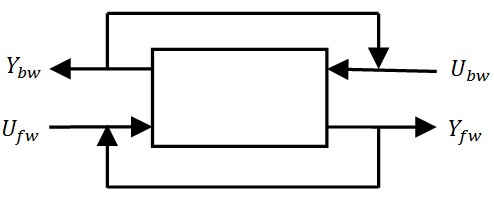}    
    \caption{A feedback on a traffic system with two inputs and two outputs.}
    \label{fig-feed}
\end{figure}

Without loss of generality, we consider a traffic system with two inputs $U_{fw}$ and $U_{bw}$ and two outputs $Y_{fw}$
and $Y_{bw}$. We then consider the system where the outputs are backed to the inputs in such a way that in the new system
the inputs are $(U_{fw}\oplus Y_{fw})$ and $(U_{bw}\oplus Y_{ww})$ rather than $U_{fw}$ and $U_{bw}$ respectively;
see Figure~\ref{fig-feed}.

The traffic system is set in feedback as shown in Figure~\ref{fig-feed}. We assume given a service matrix 
$\beta$ for the initial open system. We then have $Y\geq \beta * U$.
For the system set in feedback, we have
\begin{align}
  Y_{fw} & \geq \beta_{11}(U_{fw}\oplus Y_{fw}) \oplus \beta_{12}(U_{bw}\oplus Y_{bw}), \nonumber \\
  Y_{bw} & \geq \beta_{21}(U_{fw}\oplus Y_{fw}) \oplus \beta_{22}(U_{bw}\oplus Y_{bw}), \nonumber
\end{align}
which we can simply write $Y\geq \beta Y \oplus \beta U$. Then from Theorem~\ref{th-linsys}, we immediately get
$Y\geq (\beta^* \beta)*U$. Hence $\beta^*\beta$ is a service matrix for the system set in feedback.
Therefore we have the following result.
\begin{theorem}\label{thm-feed}
  If $\beta$ is a service matrix for an open traffic system, then $\beta^*\beta$ is a service matrix 
  for the system set in feedback.
\end{theorem}
\proof Follows directly from the arguments above. \endproof

\section{Controlled 2D traffic systems}

The calculus of guaranteed service on routes passing through intersections need to take into account the control set on every intersection.
In the communication network framework, one calculates residual services for respective flows passing through data routers.
We focus here on the case of an intersection controlled with a traffic light set in an open loop. 
In this case, the control depends only on time, and does not depend on the state of the traffic on the road entering or exiting the intersection.
One can calculate a service guarantee for every flow passing through the intersection, independent of the other flows. 

\subsection{The system of a road section controlled with a traffic light}

We consider here an urban road section controlled with a traffic light. We denote by $c$ the cycle time of the traffic light,
and by $g$ and $r$ the green and red times, with $c=g+r$. We only consider this road section, but it is implicitly supposed
that other road sections are controlled with the same traffic light, in such a way that when the light is red for the
considered road section, it is green for the other traffic systems.

The traffic dynamics in the controlled road section is the following (the same as~(\ref{dyn1a})-(\ref{dyn3a}), except the first equation~(\ref{dyn1a})).
\begin{align}
  & Q(t) \geq \min \left\{ \begin{array}{l}
			    U_{fw}\left(t-L/v - r \right)+n, \\
			    Q\left(t-L/v\right)+(g/c) q_{max}L/v, \\
			    U_{bw}(t)
                        \end{array} \right\}, \label{dyn1b}\\
  & Y_{fw}(t) = Q(t), \label{dyn2b} \\
  & Y_{bw}(t) = Q\left(t-L/w+\bar{n}\right). \label{dyn3b}
\end{align}

\begin{itemize}
  \item The term $U_{fw}(t-L/v)+n$ in~(\ref{dyn1a}) is time-shifted by $r$ in~(\ref{dyn1b}).
  \item The second term in~(\ref{dyn1b}) replaces the term $Q(t-L/v)+q_{max} L/v$ in~(\ref{dyn1a}) with $g/c<1$.
  \item The last term in~(\ref{dyn1b}) remains as in~(\ref{dyn1a}).
\end{itemize}  
Cars arriving to the traffic light may have an additional delay upper-bounded by $r$.
The dynamics~(\ref{dyn1b}) tell that the inflow to the traffic light passes through the light with a maximum time delay of $r$ time units,
under the supply constraint downstream of the light, and with a maximum flow of $(g/c) q_{\max}$.

\begin{theorem}\label{thm_rsc}
   A service matrix $\beta$ for the controlled road section is $C*A^* * B$, with 
   $A=\gamma^{(g/c)q_{\max}L/v} \delta^{L/v}, B=(\gamma^n \delta^{L/v+r} \quad e)$.
\end{theorem}   
\proof. Same proof as the one of Theorem~\ref{thm-impulse}. \endproof

\section{Roads and itineraries}

In order to build a road of $m$ sections, we need to compose $m$ elementary traffic systems of road sections.
The service matrix of each road section can be obtained by Theorem~\ref{thm-impulse}, giving fundamental traffic diagrams
on each section. Then the service matrix of the whole road is obtained by the composition of the road section 
systems and by applying Theorem~\ref{thm-connect}. A controlled road of $m$ sections is obtained similarly by 
composing $(m-1)$ uncontrolled road sections with one controlled road section.

\begin{figure}[h]
  \begin{center}
    \includegraphics[width=6cm]{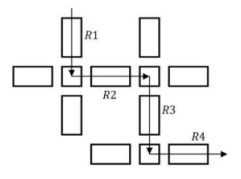}
    \caption{Route traffic system.}
  \end{center} 
  \label{fig-itin}
\end{figure}

A route (or an itinerary) in a controlled road network is build by composing a number of controlled roads.
In Figure~\ref{fig-itin}, we illustrate the composition of controlled roads to obtain a traffic system associated
to a whole road network. The procedure of computing a service matrix for a traffic flow passing respectively
through roads R1, R2,R3, and R4 is the following.
\begin{itemize}
    \item Determine service matrices for all the uncontrolled sections of the itinerary, by Theorem~\ref{thm-impulse}.
    \item Determine service matrices for all the controlled sections of the itinerary, by Theorem~\ref{thm_rsc}.
    \item Determine a service matrix for the itinerary by connecting the systems R1, R2, R3, R4, by Theorem~\ref{thm-connect}.
\end{itemize}	

Once a service matrix is determined for a given traffic system, and having an arrival matrix expressing the traffic
demand arriving to the system and the traffic supply that backwards to it, it suffices to apply Theorem~\ref{th-vec-del} to
obtain upper bounds for the travel time for an input - output couple of the traffic system.

\begin{table}\label{tab-itiner}
  \centering
  \begin{tabular}{c||c|c|c|c}
	  & R1 & R2 & R3 & R4 \\
	  \hline
	  \hline
      Length $L$ (meter) & 150 & 150 & 100 & 100 \\
      \hline
      Maximum flow $q_{\max}$ (veh/sec) & 0.32 & 0.35 & 0.4 & 0.38 \\
      \hline
      Initial car-density $n/L$ (veh/meter) & 5/150 & 10/150 & 3/100 & 7/100 \\
      \hline
      Cycle time $c$ (sec.) & 60 & 90 & 80 & - \\
      \hline
      Green time $g$ (sec.) & 30 & 50 & 45 & - 
   \end{tabular}
   \caption{Parameters of the numerical example.}
\end{table}

\begin{figure}[h]
    \includegraphics[width=1\textwidth]{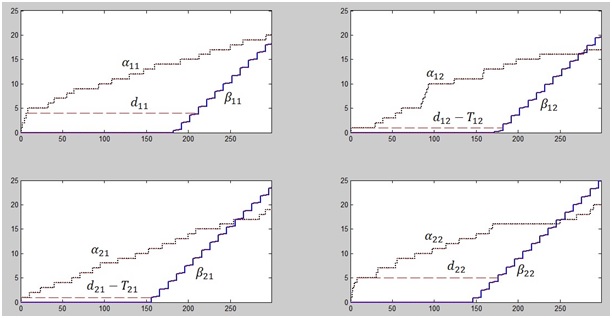}
    \caption{Arrival curves of the arrival matrix, service curves of the service matrix, and time delays.}
    \label{figg4}
\end{figure}

The results of this example are shown in Figure~\ref{figg4}. The input signals $U_{fw}$ arriving to
road~1, and $U_{bw}$ backing from road~4 are taken such that the arrival flows do not exceed the service 
offered by the whole route. The arrival curves of the arrival matrix $\alpha$ are computed following 
Definition~\ref{def-arrival}. First the shift times $T_{12}=60$ s., and $T_{21}=8$ s. are computed. Then the
curves are deduced by Definition~\ref{def-arrival}. The service curves are computed following the steps cited
above in this section. An upper bound for the travel time through the route is then calculated according to 
Theorem~\ref{th-vec-del}. We are concerned here by the delay $d_1$ corresponding to the forward travel time
(the delay $d_2$ corresponds to the backward travel time of the backward waves). We obtained for 
this example the following result.
$d_1 = \max (d_{11},d_{12}) = \max (205,241) = 241$ seconds.

\section{Possible extensions and limits}

Numerous extensions of this approach are possible. We enumerate below some of them.
\begin{itemize}
  \item By extending the approach to intersections controlled by traffic lights where the decision depends on the feedback of the state of the traffic,
     one may derive performance bounds corresponding to given control strategies at intersections. Different control strategies for urban traffic
     can be studied. Extension to controlled network is then straightforward.
  \item Extension to the case of intersections managed by priority rules can also be considered.
  \item Extension to roundabouts.
  \item Extension to highway merge and divergent.
  \item Extension to highway traffic control (ramp metering).
  \item Extension to stochastic guarantees, by applying stochastic network calculus.
\end{itemize}

Finally, let us mention one of the main limits of this approach, and of the applications of network calculus in general.
In case where cyclic dependencies of the different considered flows, the unique existing way to derive the residuated service curves is
to solve implicit systems deduced from all the dependencies of the different flows.
It is known that this way is not optimal, and the obtained bounds are in general not tight; see~\cite{Cha00} for example.


\newpage
This page is intentionally left blank

\chapter[Summary of my other contributions]{Summary of my other contributions on the deterministic traffic modeling and control}
\label{chap-other-det}

This chapter summarizes some of my other contributions on the deterministic traffic modeling and control.
Four contributions are summarized. 

First, we summarize a work on the network calculus theory, where we proposed a new concept of \textit{packet curves}, which permits to
take into account the variance in packet sizes of arrival flows to given servers.
The existing classical approach takes into account only the maximum packet size for each arrival flows, which is very pessimistic.
The approach we propose here consists in associating to arrival flows curves giving for an amount of data the minimum number of packets.
This approach follows the same philosophy as the one used for upper-bounding arrival flows, and adapts it to the number of packets, ot equivalently
to their sizes. The main references of this part are~\cite{FG10,BFG11,BFG12}.

Second, we present briefly a new semi-decentralized approach for urban traffic control, which we proposed in 2015.
We have based on a the centralized control strategy TUC (Traffic Urban Control)~\cite{Dia99} and
introduced a new time interval within the time cycle, 
whose length is controlled by central level, and during which the traffic is locally controlled.
We have shown that the traffic management in urban networks can be improved by this kind of semi-decentralized controls.
Indeed, in free flow traffic conditions, where the intersections operate under capacity, it is possible to set off the traffic lights, and 
alternate a priority rule at the intersection. Our results are confirmed by means of numeric simulation using the traffic simulator SUMO.
The main reference of this part is~\cite{FNHL15}.

Third, we summarize a work done in collaboration with Mrs Souaad Lahlah from University of Bejaja (Algeria), where we proposed a new approach
for optimal path selection for VANETs. It is new proactive protocol where the selection of optimal paths is based on two criteria: the distance and the car-density.
The latter is estimated basing on estimations of the car-speed on the links of the network, and using given fundamental diagrams for the links.
The performance of the proposed protocol are assessed by numerical simulation using the vehicular traffic simulator SUMO combined with
the communication simulator ns2. The new proposed protocol has been compared to existing protocols.
The main reference of this part is~\cite{LSBF18}.

Fourth, we give some details on a first algorithm we developed for urban traffic control including v2i (vehicle to infrastructure) communications.
We have assessed the algorithm with the open-source tool \textit{Veins} which combines the vehicular traffic simulator SUMO with the communication
simulator Omnet++. 
We present here briefly the algorithm and evaluate its performances in term of vehicular traffic and in term of communications. 
The main reference of this part is~\cite{NFHL17}.

\newpage
\section{Packet curves in network calculus}

The main references of this part are~\cite{FG10,BFG11,BFG12}.
We proposed in  this work a new notion of \textit{packet curve} for data packetization in the theory
of \textit{Network Calculus}. Packet curves permit to optimize the utilization of all available information
of the packet characteristics in communication networks. Moreover, the definition of packet curves we made
here allows us to efficiently improve the derivation and the tightness of performance bounds in such networks.
This is useful in particular in the case where we have aggregation of flows, packet-based service policies and shared
buffers. We illustrate this new notion by a model of \textit{wormhole switch} and show how our results can be used
to get efficient delay bounds.
We give here the main definitions and results without any proof. All the details and proofs can be found in~\cite{BFG11, BFG12}.

Let us consider an arrival flow $A$ for a given server $S$, where
$A(t)$ denotes the cumulative arrival data from time zero to time $t$.
We denote by $B$ the departure flow, in such a way that $B(t)$ gives
the cumulative departure data from the server from time zero to time $t$.
We assume that the data arrive in packets of different sizes.

\begin{definition}
  A packet operator $\mathcal P: \mathbb R_+ \to \mathbb N$ gives for any amount $x\in\mathbb R_+$ of data 
  the number of entire packets in $x$.
\end{definition}

For example, the application of the packet operator $P$ to an arrival flow data $A(t)$ at time $t$
gives the number of entire packets of $A$ arrived from time zero to time $t$.
Therefore, $P$ transforms an arrival flow $A$ to a packet flow $P$, where $P(t) = \mathcal P \circ A(t)$.

The packet flow $P$ may not be perfectly known, but some information about it may be available, more precise than
only the minimum and maximum packet length respectively denoted by $l_{\min}$ and $l_{\max}$.

\begin{example}\label{examp1}
Let us for example consider a flow with packets of sizes~1 and~2.
We may have more information than only the minimum and maximum length sizes~1 and~2.
For example, we may know that in three successive packets, there are
at least one packet of size~1 and at least one packet on size~2.
\end{example}

In order to take into account such additional information, we introduce here the notion of
\textit{packet curve} of a packet operator.

\begin{definition}
  A curve $\pi$ (resp. $\Pi$) is a minimum (resp. maximum) packet curve for $\mathcal P$ if
  $\forall 0 \leq x \leq y, \mathcal P(y) - \mathcal P(x) \geq \pi(y - x)$ (resp. $\mathcal P(y) - \mathcal P(x) \leq \Pi(y - x)$).
\end{definition}

For the case considered in Example~\ref{examp1}, one can take as a minimum packet curve, the curve
$$\pi : x \mapsto (3/5(x - 2/3))^+$$
and as a maximum packet curve, the curve
$$ \Pi : x \mapsto 3/4 x + 3/2,$$
where $(x)^+ := \max(0, x)$.

Stair-case functions can be more accurate, but rate-latency functions and affine functions are easier to handle from an computational viewpoint.
Note that if $\pi$ is defined as $\pi : x \mapsto \mu (x - \nu)^+$ , then $\nu \geq l_{\max}$ and $\mu \geq 1 / l_{\max}$
and if $\Pi$ is defined as $\Pi : x \mapsto V + U x$, then $V \leq 1/l_{\min}$ and $V \geq 1$.

\subsection*{Main properties}

The following result gives some properties derived directly from the definitions above.
\begin{proposition}
   Let us consider arrival flows $A_1$ and $A_2$ to a given server, with associated departure flows $B_1$ and $B_2$ respectively,
   and with packet operators $\mathcal P_1$ and $\mathcal P_2$ respectively. We denote $A = A_1 + A_2$ and $B = B_1 + B_2$.
   We consider also a packet operator $\mathcal P$ for $A$.
   If $\alpha, \alpha_1$ and $\alpha_2$ are maximum arrival curves for $A, A_1$ and $A_2$ respectively, if $\pi, \pi_1$ and $\pi_2$
   are minimum, and $\Pi, \Pi_1, \Pi-2$ are maximum packet curves for $\mathcal P, \mathcal P_1$ and $\mathcal Pi_2$ respectively, and if
   $\beta$ is a minimum strict service curve for the server, then we have the following.
   \begin{enumerate}
      \item $\Pi$  and $\pi$ are maximal and minimal packet curves for $B$.
      \item $\Pi \circ \alpha$ is an arrival curve for the packet flow $P := \mathcal P \circ A$.
      \item $\pi \circ \beta$ is a minimum strict service curve for $P$ .
      \item If $\beta'$ is a minimum simple (resp. strict) service curve for packet flow $P$, then, $(\Pi)^{-1} \circ \lfloor \beta' \rfloor$
         (resp. $(\Pi)^{-1} \circ \lceil \beta' \rceil$) is a minimum simple (resp. strict) service curve for $\mathcal P^{-1} \circ \mathcal P \circ A$.
      \item $\pi_1 * \pi_2$ is a minimum packet curve for the (blind) aggregation of $A_1$ and $A_2$.
   \end{enumerate}
\end{proposition}

\proof See~\cite{BFG11} \endproof

\subsection*{Examples of modeling network elements}

\subsubsection*{a) Superposition of periodic flows}

\begin{theorem}(Superposition of periodic flows)
  For periodic arrival flows $A_1, A_2, \ldots, A_n$ such that flow $A_i$ includes packets of size $S_i$ arriving with period $T_i$ for $1\leq i\leq n$,
  the total flow $A = \sum_i A_i$ where $A_i$ are aggregated in FIFO admits minimum and maximum packet curves $\pi$ and $\Pi$ given as follows.
  $$\pi(x) = \left(\sum_{i=1}^n \frac{x}{T_i \sum_j \rho_j} - \sum_{i=1}^n \frac{\sum_j T_j \rho_j}{T_i \sum_j \rho_j}  \right)^+$$
  $$\Pi(x) = \sum_{i=1}^n \frac{x}{T_i \sum_j \rho_j} + \sum_{i=1}^n \frac{\sum_j T_j \rho_j}{T_i \sum_j \rho_j}$$
  where $\rho_i := S_i/T_i, \forall i=1,2, \ldots,n$.
  The rates of the two curves $\pi(x)$ and $\Pi(x)$ are equal. Therefore, $\pi$ and $\Pi$ are optimal. By consequent, $\pi$ is better than
  $\pi_1 * \pi_2 * \ldots * \pi_n$.
\end{theorem}

\proof See~\cite{BFG11} \endproof

\subsubsection*{b) Non-preemptive service curves}

\begin{figure}[htbp]
   \centering
   \includegraphics[scale = 0.7]{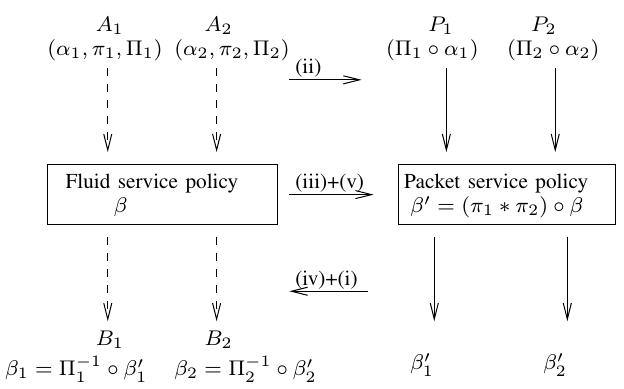}
   \caption{Non-preemptive service calculus scheme.}
   \label{packet1}
\end{figure}

Let us now show how one can derive non-preemptive service curves by using this new notion of packet curve.
Figure~\ref{packet1} gives the general scheme of computation, using the basic properties of packet curves.
We notice that this scheme is more efficient than the classical one taking into account the minimum
and maximum packet sizes, only when the service policy is based on counting packets.
We illustrate here with an example of the Round-Robin routing policy.
We consider two flows arriving to a server operating under the round-robin policy. 
We are interested here in the calculus of the residual service curve for flow~1.

We present here the derivation of a service curve for flow~1 with three different methods.
\begin{itemize}
  \item Classical method:
    $$\beta_c(t) = \left( \frac{l_1^{\min}}{n l_1^{\max}} \beta - l_2^{\max} \right)^+$$
  \item Scheme method:
    $$\beta_s(t) = \left( \Pi_1\right)^{-1} \left( \frac{1}{n} (\pi_1 * \pi_2 \circ \beta) - 1\right)^+$$
  \item Ad-hoc method:
    $$\beta_{ah}(t) = \left( Id + \pi_2^{-1}\circ \Pi_1 + 1 \right)^{-1} \circ \left(\beta - l_2^{\max} \right)^+$$
\end{itemize} 

\begin{example}
  With $l^{\min} = 1, l^{\max} = 2, \pi_1 = \pi_2 = \pi, \Pi_1 = \Pi_2 = \Pi$, and $\beta(t) = 10 t$, we have
  \begin{itemize}
    \item $\beta_c(t) = 2.5(t - 0.8)^+$,
    \item $\beta_s(t) = 4(t - 0.97)^+$,
    \item $\beta_{ah}(t) = 4.44(t - 0.68)^+$.
  \end{itemize}
  The residual service rate using the scheme method is much better than the one obtained by the classical method.
  The ad-hoc method is the best one.
\end{example}

\subsubsection*{c) Shared queues}

Let $A_1$ and $A_2$ be two arrival packeted flows to a given server.
The latter serves the flows packet by packet picking a packet from flow $A_1$
or flow $A_2$ respecting FIFO service, without preemption.
We assume that when the service of a packet from a given flow finishes, then
if the next one is from a different flow, then the service is re-initiated,
as if the switching time is the beginning of a backlogged period.
By that we can consider that the servers provide a different service for 
flows $A_1$ and $A_2$: packets of flow $A_1$ have a strict minimum service curve $\beta_1$
and packets of flow $A_2$ have a strict minimum service curve $\beta_2$. 

\begin{figure}[htbp]
   \centering
   \includegraphics[scale = 0.6]{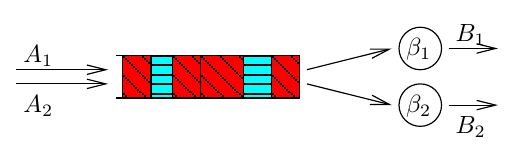}
   \caption{Shared queues: the two servers cannot be active at the same time.}
\end{figure}

Then we have the following result.
\begin{theorem}(Shared queues)\label{thm-packet2}
  Considering only the periods of time when server~1 is active (and server~2 is idle), $\tilde{\beta}_1$ is a service curve for flow~1, with
  $$\tilde{\beta}_1^{-1} (x) = \left( \beta_1^{\lfloor \Pi_1(x)\rfloor}\right)^{-1} \left( \Pi_1^{-1}\left( \lfloor \Pi_1(x)\rfloor\right)\right) + \beta_1^{-1}\left( \pi_1^{-1}(0^+)\right).$$
\end{theorem}
\proof See~\cite{BFG11} \endproof

Then an overall service curve for the shared queue is $\tilde{\beta}_1 * \tilde{\beta}_2$.
Even though the curve given by Theorem~\ref{thm-packet2} is rather pessimistic, it gives the key idea to
find better service curves; that is, for each flow compressing time when the server of interest is idle.
Then by bounding the number of idle periods and the lengths of those idle periods, we improve the curve.
This is possible because minimal and maximal packet curves are known.

\subsection*{Application}

\begin{figure}[htbp]
   \centering
   \includegraphics[scale = 0.8]{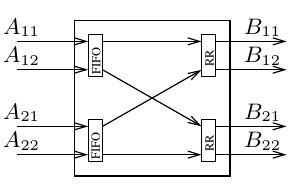}
   \caption{A $2 \times 2$ switch.}
   \label{fig-packet3}
\end{figure}

To illustrate the results presented above, we give here an application to performance evaluation on a $2 \times 2$ switch; see Fig.~\ref{fig-packet3}.
We assume that the input ports of the switch serve the arrival flows under FIFO policy, while the output ports serve them under packetized round-robin (RR) policy
(with fixed service rate).

We propose the following steps for the calculus of service curves for the aggregate as well as the residuated flows considered here.
The service curves can then be used for example to upper bound the delays of serving different flows.
The steps are the following.
\begin{enumerate}
  \item Compute individual service curves for the output ports (RR).
  \item Compute a service curve for the shared queues (packetized FIFO).
  \item Use the derived service curves to calculate maximum arrival curves for $B_{ij}$.
  \item Go to step 1) and iterate using the additional knowledge obtained.
\end{enumerate}

\newpage
\section{Semi-decentralized urban traffic control}

We summarize here the semi-decentralized approach we proposed in~\cite{FNHL15} for urban traffic control.
We based on the TUC (Traffic responsive Urban Control) strategy~\cite{Dia99,DPA02,DDAPBSL03}, where
antagonistic stages alternate the pass of way during a fixed or controlled cycle time, and we introduce a 
new contention time window inside the cycle time, where antagonistic stages can alternate a priority rule.
The priority rule gives green color to given stages (priority) and yellow color for other antagonistic stages
(no priority). By introducing this new time widow inside the cycle time, we reduce the red color, and by that
increase the capacity of the network junctions.
Such priority rule should be realizable by considering possible vehicle to vehicle (v2v) communications
where the vehicles entering to a junction from priority approaches can pass the junction without decelerating
while the ones entering from non priority approaches decelerate and yield the way.
We propose a linear quadratic model for the dynamics and the control of such a system, with an
illustration on a small grid network where the traffic is simulated with the traffic simulator SUMO 
(Simulation of Urban MObility). We compare our approach with the classical TUC strategy.
Let us consider the following notations.~~\\~~

\begin{tabular}{ll}
$c$			& cycle time duration, in seconds.\\
$k$			& discrete time index, corresponding to a duration of $kc$ sec.\\
$x_i(k)$		& number of cars on link $i$ at discrete time $k$.\\
$\bar{x}_i$		& constant nominal number of cars on link $i$.\\
$\Delta x_i(k)$		& $:= x_i(k) - \bar{x}_i$.\\
$s_i$			& saturation flow on link $i$.\\
$g_i(k)$		& green time duration for link $i$ during the $k$th cycle.\\
$\bar{g}_i$ 		& constant nominal green time duration for the stream coming from link $i$. \\
$\Delta g_i(k)$ 	& $:= g_i(k) - \bar{g}_i$. \\
$u_i(k)$		& $:= (g_i(k)/c) s_i$ average outflow from link $i$ during the $k$th cycle.\\
$d_i(k)$ 		& arrival demand flow to link $i$ at discrete time $k$. \\
$\alpha_{ij}$		& turning movement ratio from link $i$ to link $j$.
\end{tabular}~~\\~~

We recall below the principle of store and forward control strategies, like TUC.
The number of cars on link $i$ is updated as follows.
\begin{equation} \label{dyn0}
  x_i(k+1) = x_i(k) + d_i(k) + \sum_{j} \alpha_{ji} s_j g_j(k) - s_i g_i(k),
\end{equation}
where $j$ indexes here the upstream links of link $i$.
By introducing the nominal amounts, and by using vectorial notations, we get
\begin{equation}\label{dyn1}
  \Delta x(k+1) = \Delta x(k) + B \Delta g(k),
\end{equation}
where $B$ is a matrix built basing on the dynamics~(\ref{dyn0}) written on the whole network, and where 
the variations of the arrival demand flows on every link inside the cycle time are assumed to sum to zero.
Bounds for minimum green times and maximum storage capacity of links have also to be considered.

The criterion is the following, where $\lambda$ is a discount factor, and where an infinite time horizon is considered.
\begin{equation}\label{crit1}
  J = \min_{\Delta g} \frac{1}{2} \sum_{k=0}^{+\infty} \frac{1}{(1+\lambda)^k} \left( \|\Delta x(k)\|_Q^2 + \|\Delta g(k) \|_R^2\right),
\end{equation}
$Q$ and $R$ are non-negative definite, diagonal weighting matrices.
The first term on~(\ref{crit1}) aims to minimize the risk of over-saturation and the spill-back of link queues, while the second term
is used to influence the magnitude of the control.

\subsection*{Semi-decentralization}

We propose here to divide the red color time into two time periods: red and yellow.

\begin{figure}[htbp]
  \begin{center}
    \includegraphics[width=6cm]{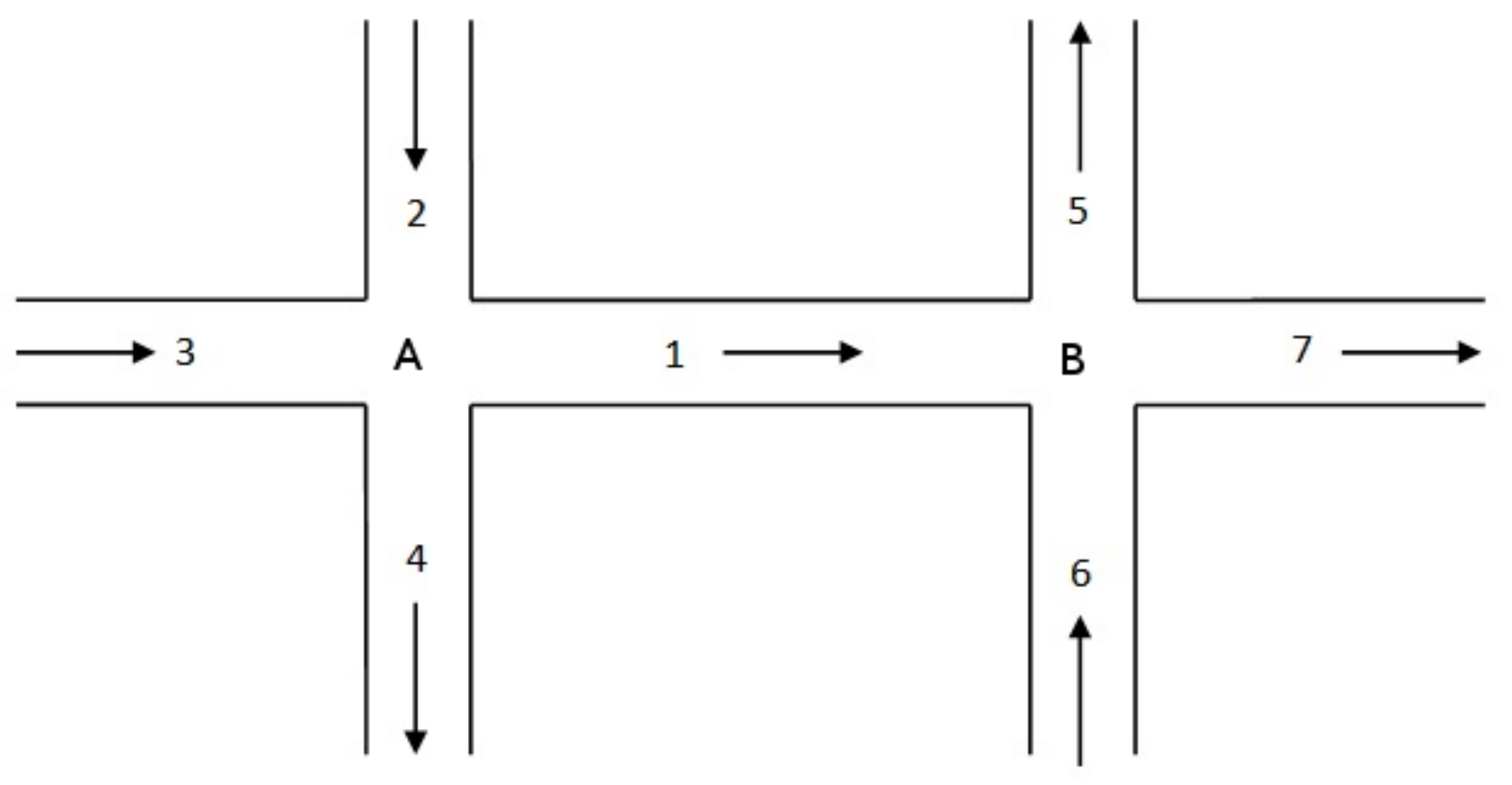} \hspace{1cm}
    \includegraphics[width=6cm]{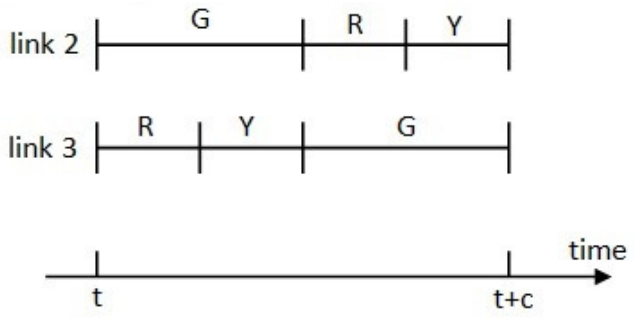}    
    \caption{At intersection A, we note G: green, R: red, Y: yellow.}
    \label{cycle-time}
  \end{center}
\end{figure}

By considering yellow time durations, we  need to choose three independent control variables, among six variables; see Figure~\ref{cycle-time}.
According to Figure~\ref{cycle-time}, we can for example consider the following three independent control variables.
\begin{itemize}
  \item $g_2(k)$ : green time duration for link 2, 
  \item $y_2(k)$ : yellow time duration for link 2, 
  \item $y_3(k)$ : yellow time duration for link 3.
\end{itemize}  
The other three dependent control variables are given as follows (see figure~\ref{cycle-time}). 
\begin{itemize}
  \item $r_2(k) = c - g_2(k) - y_2(k)$ : red time duration for link 2,
  \item $r_3(k) = g_2(k) - y_3(k)$ : red time duration for link 3,
  \item $g_3(k) = c - g_2(k)$ : green time duration for link 3.
\end{itemize}

\subsubsection*{The dynamics}

Let us consider the following additional notations.
\begin{itemize}
  \item $q^{\max}_{J}$ : capacity (maximum flow) of junction $J$.
  \item $Q_{ij}(k)$ : total flow going from link $i$ to link $j$ during the $k$th cycle.
  \item $Q^{out}_i$ : total flow exiting from link $i$ during the $k$th cycle.
  \item $\gamma_J$ : friction coefficient on junction $J$, with $0\leq \gamma \leq 1$.
\end{itemize}

We write the traffic dynamics on link 1 of Figure~\ref{cycle-time} with the new control model.

\begin{equation} \nonumber
  x_1(k+1) = x_1(k) + d_1(k) + Q_{21}(k) + Q_{31}(k) - Q^{out}_1(k),
\end{equation}
where $Q_{21}(k), Q_{31}(k)$ and $Q^{out}_1(k)$ are given in~(\ref{eqq1})-(\ref{eqq3}), and where we introduce a new parameter $\gamma_J$ 
(for junction $J$) which we call a friction coefficient, and which expresses the 
interaction between vehicles entering into a junction from antagonistic stages during the contention time window.
For example, in~(\ref{eqq1}), the flows of vehicles going from link~2 to link~1 during different time durations of 
the $k$th cycle are given as follows.
\begin{itemize}
 \item During $r_3(k) = g_2(k) - y_3(k)$, the flow is $\alpha_{21}s_2(g_2(k)-y_3(k))$, as usual.
 \item During $y_3(k)$, the flow is $\alpha_{21}s_2y_3(k)$ as usual, but multiplied by the friction coefficient
   $\gamma_A$ between the streams coming from link 2 (with green time) and link 3 (with yellow time),
   since the local control is activated with a priority rule setting. Link 2 has priority over link 3 during this time period.
 \item During $r_2(k)$, the flow is zero.
 \item During $y_2(k)$, the stream coming from link 3 has priority over the one coming from link 2.
   Therefore, the whole junction capacity $q_A^{\max}y_2(k)$ is used by the stream of link 3, and
   the remaining capacity $q_A^{\max}y_2(k) - s_3 y_2(k)$ is used by link 2. This flow is also multiplied
   by the coefficient friction $\gamma_A$ since the two streams pass through junction $A$ during the same time period.
\end{itemize}
\begin{align}
   Q_{21}(k) & = \alpha_{21} s_2 (g_2(k) - y_3(k)) + \gamma_A \alpha_{21} s_2 y_3(k) + \gamma_A (q_A^{\max} y_2(k) - s_3 y_2(k)). \label{eqq1} \\
   Q_{31}(k) & =  \alpha_{31} s_3 (c - g_2(k) - y_2(k)) + \gamma_A \alpha_{31} s_3 y_2(k) + \gamma_A (q_A^{\max} y_3(k) - s_2 y_3(k)). \label{eqq2} \\
   Q^{out}_1(k) & = s_1 (g_1(k) - y_6(k)) + \gamma_B s_1 y_6(k) + \gamma_B (q_B^{\max} y_1(k) - s_6 y_1(k)). \label{eqq3}
\end{align}

The dynamics~(\ref{eqq1})-(\ref{eqq3}) are still linear on the variables $x_i, g_i$ and $y_i$.
We notice here that the dynamics are written on the independent controls.
As in the classical TUC model, we consider nominal demands $\bar{d}_i$, nominal numbers of cars $\bar{x}_i$ and nominal
independent controls $\bar{g}_i$ and $\bar{y}_i$.
The choices of $\bar{x}$ and $\bar{g}$ can be done by the same way as in the classical TUC model. 
One way to choose $\bar{y}$ is to take  $\bar{y}_i = c - \bar{g}_i$. This is equivalent to fix the nominal red time to zero.
This choice can also be dependent on the junction design.
Then we can easily derive a linear dynamics similar to~(\ref{dyn1}).
The criterion~(\ref{crit1}) is written with the new (independent) control variables $\Delta g_i$.
We obtain a linear quadratic problem. The optimal control is derived by solving a Riccati equation as in the classical TUC model. 

\subsection*{Numerical results}

\begin{figure} \label{fig-sim}
  \begin{center}
    \begin{tabular}{c||c}
      \includegraphics[width=6cm]{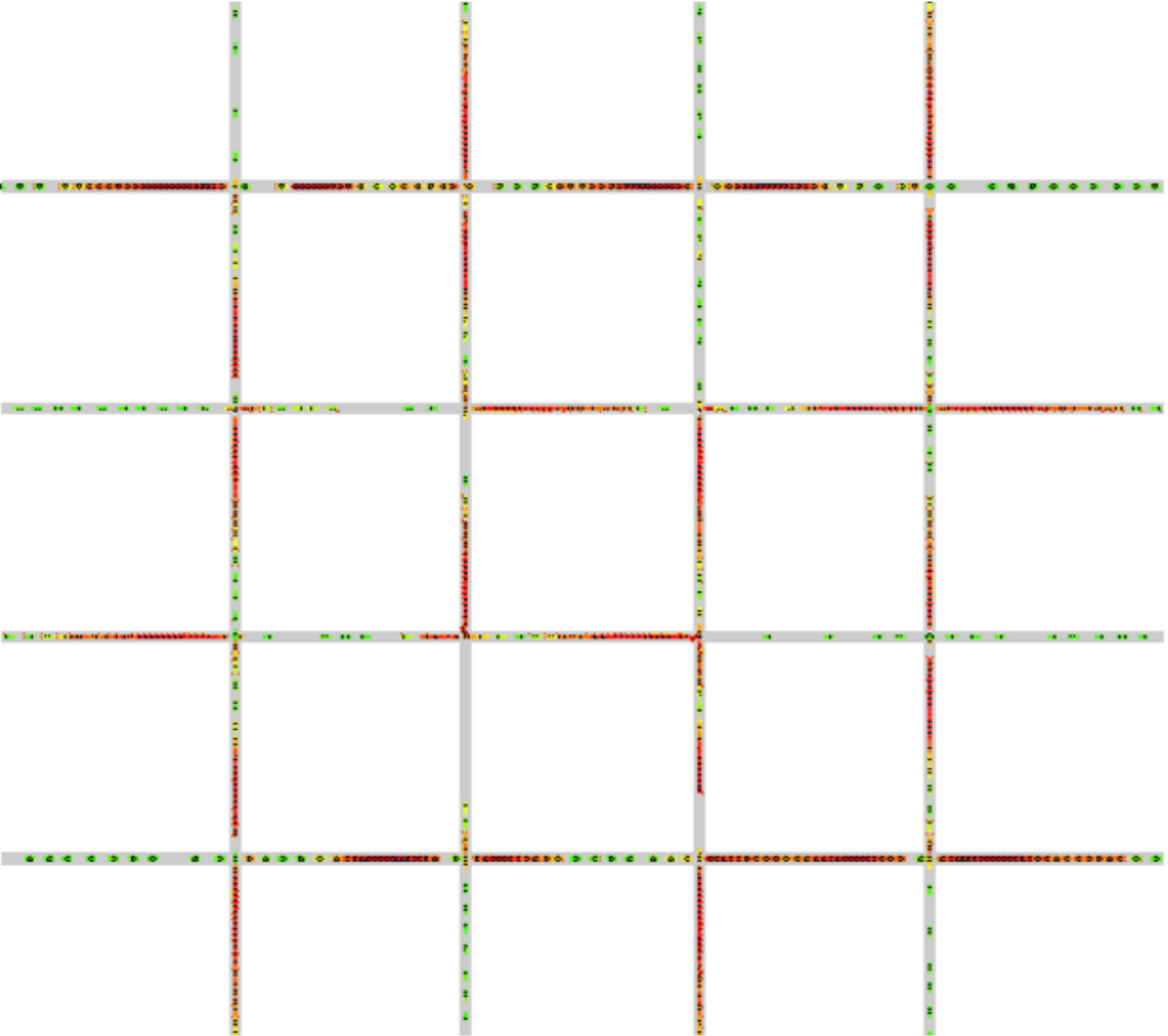} \hspace{1cm} & \hspace{1cm} \includegraphics[width=6cm]{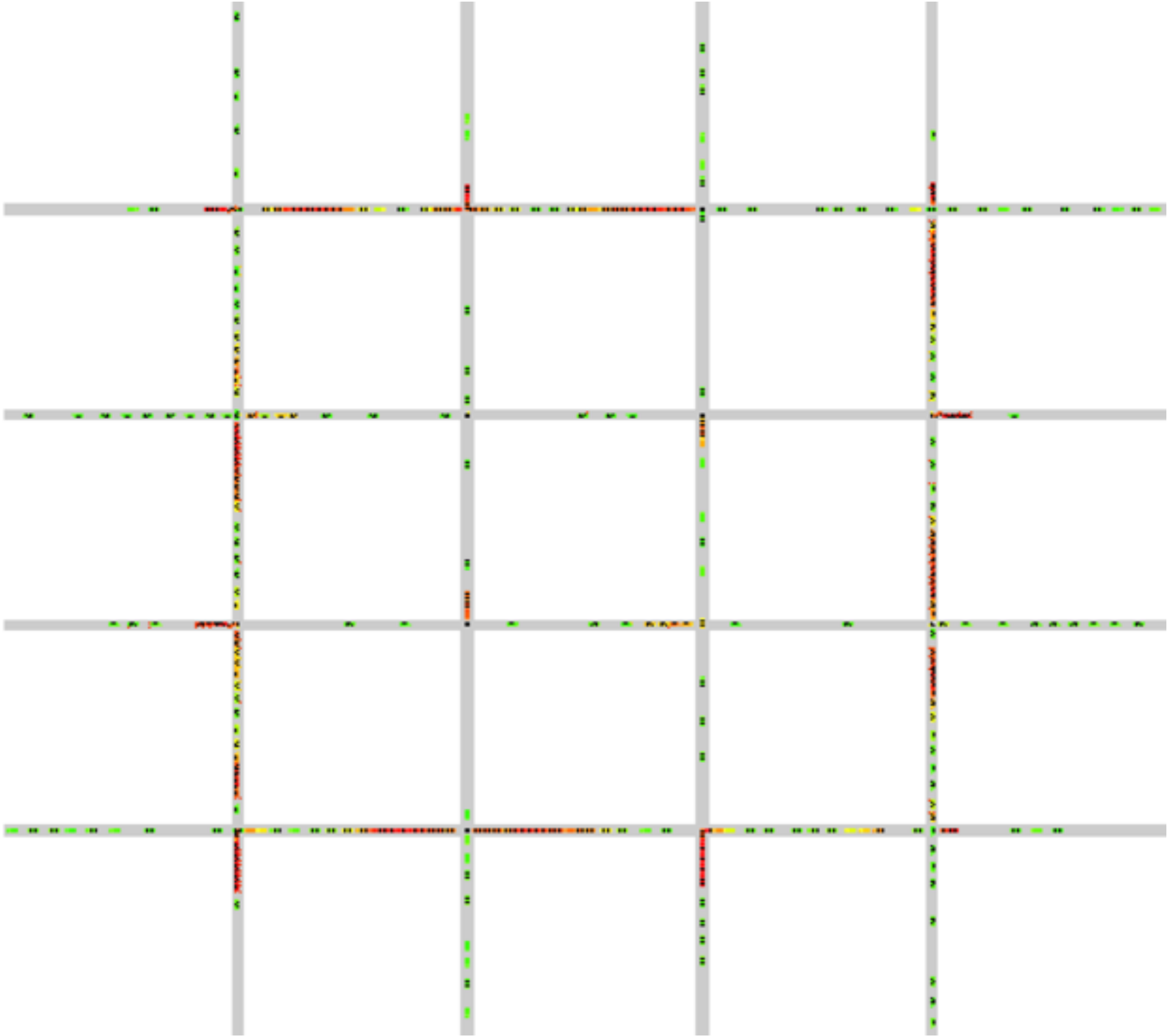}
    \end{tabular}
    \caption{The state of the traffic at the end of simulation. The colors of vehicles correspond to their speed
      (green: high speed, red: low speed).
      On the left side: Centralized TUC. On the right side: semi-centralized TUC.}
    \label{figNet}
  \end{center}
\end{figure}

\begin{figure}[htbp]
\begin{center}
\begin{tabular}{cc}
  \includegraphics[width=77mm]{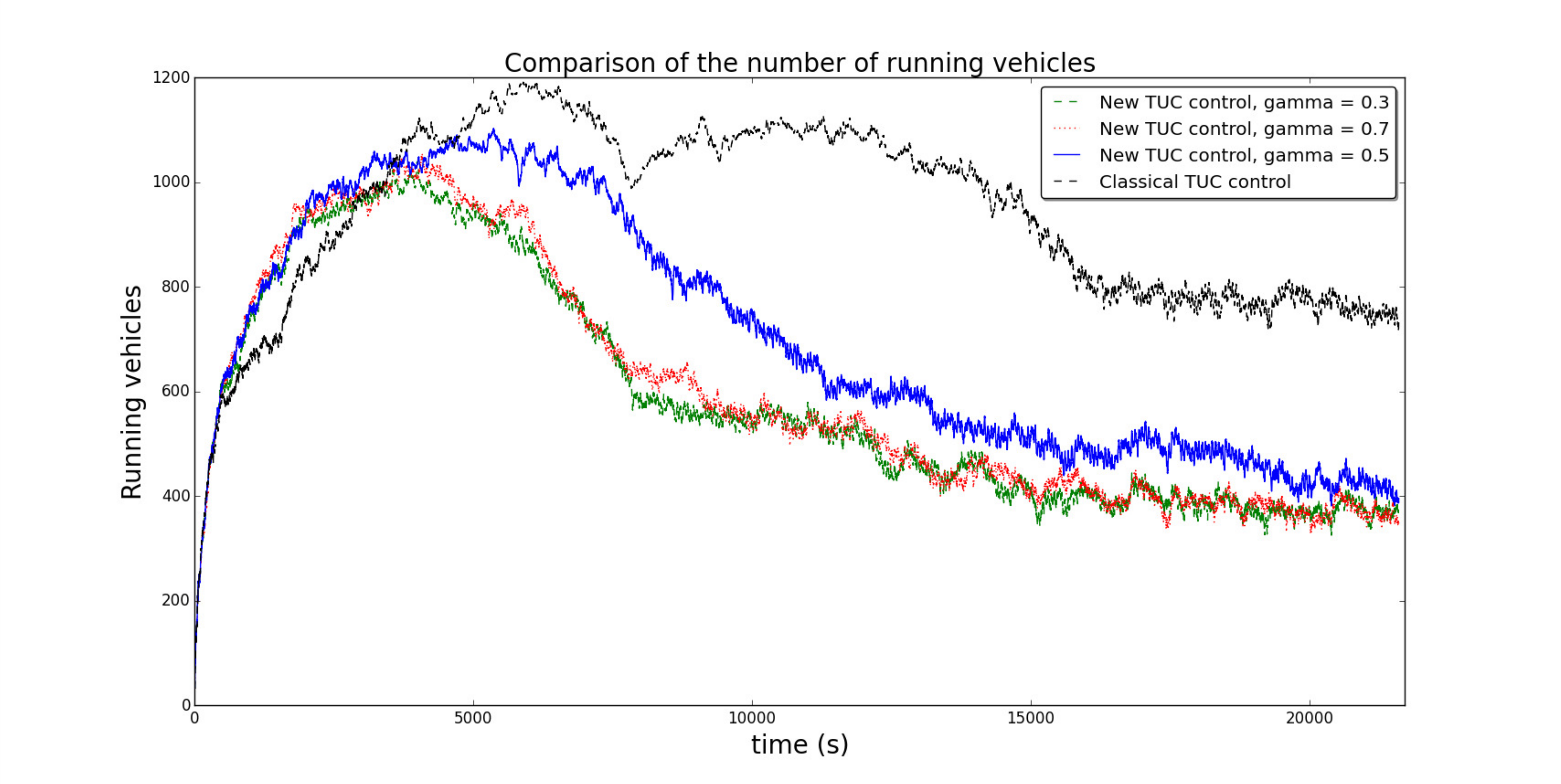} & \includegraphics[width=77mm]{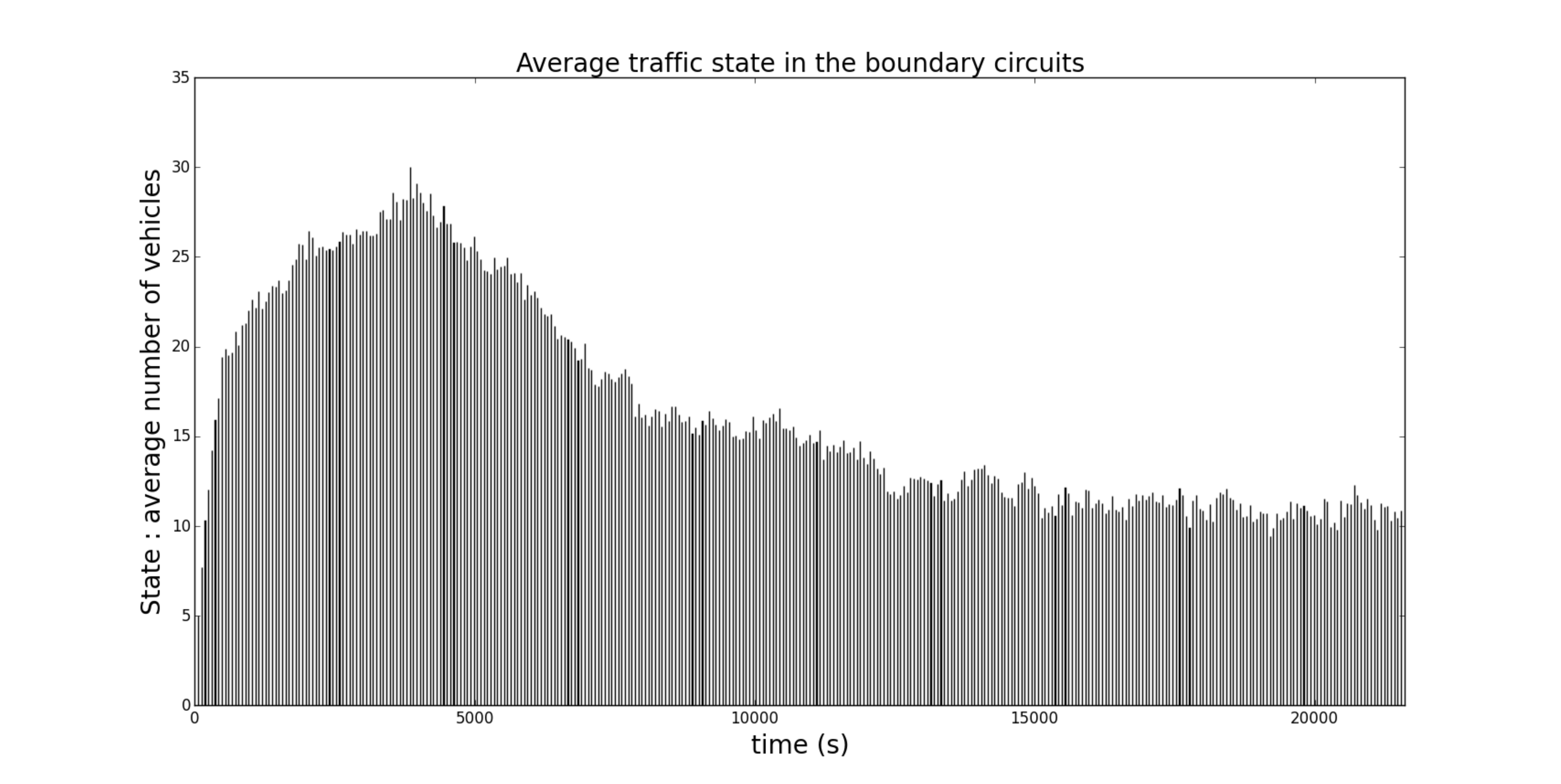} \\
  \includegraphics[width=77mm]{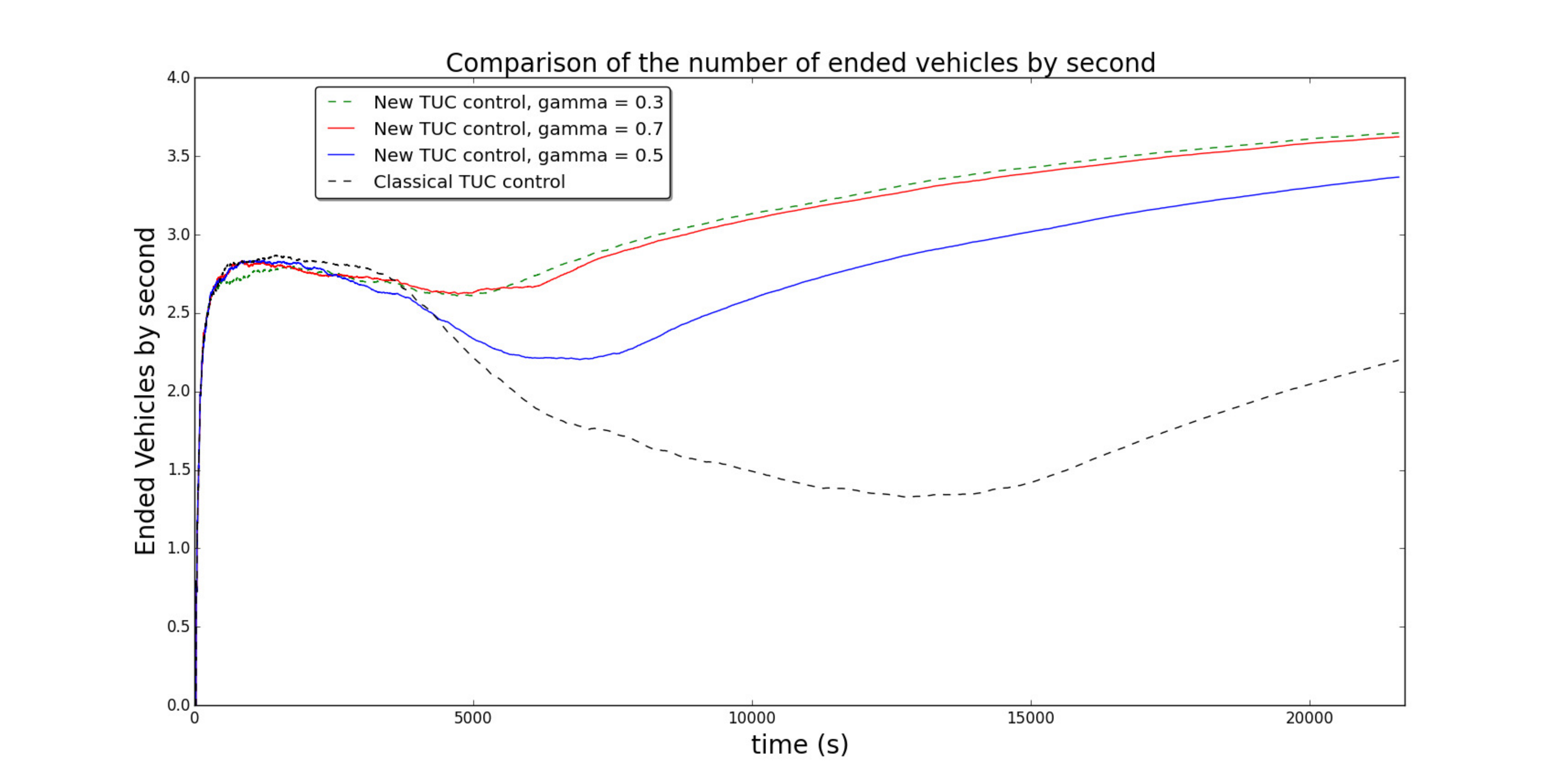} & \includegraphics[width=77mm]{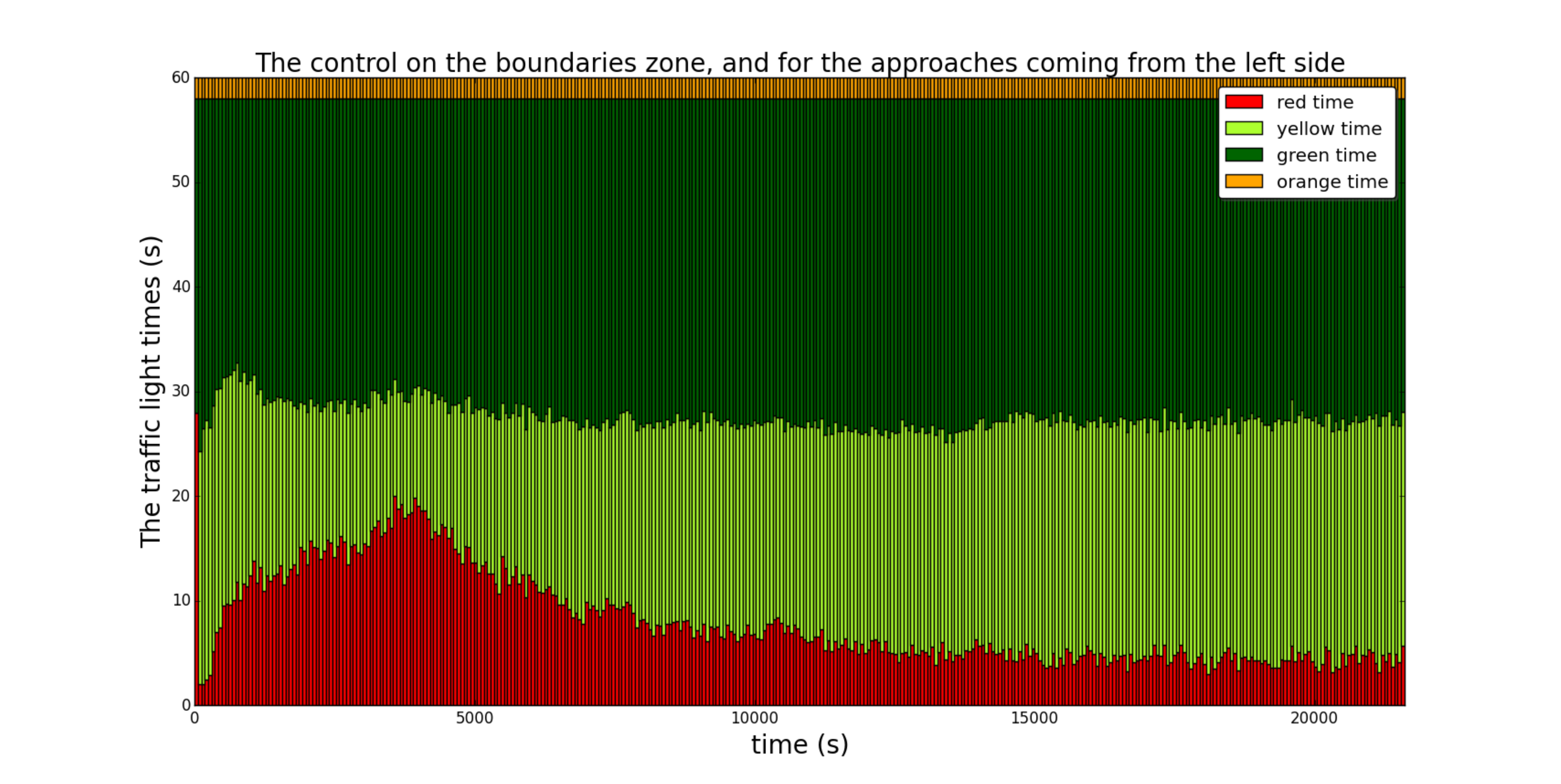} \\
  \includegraphics[width=77mm]{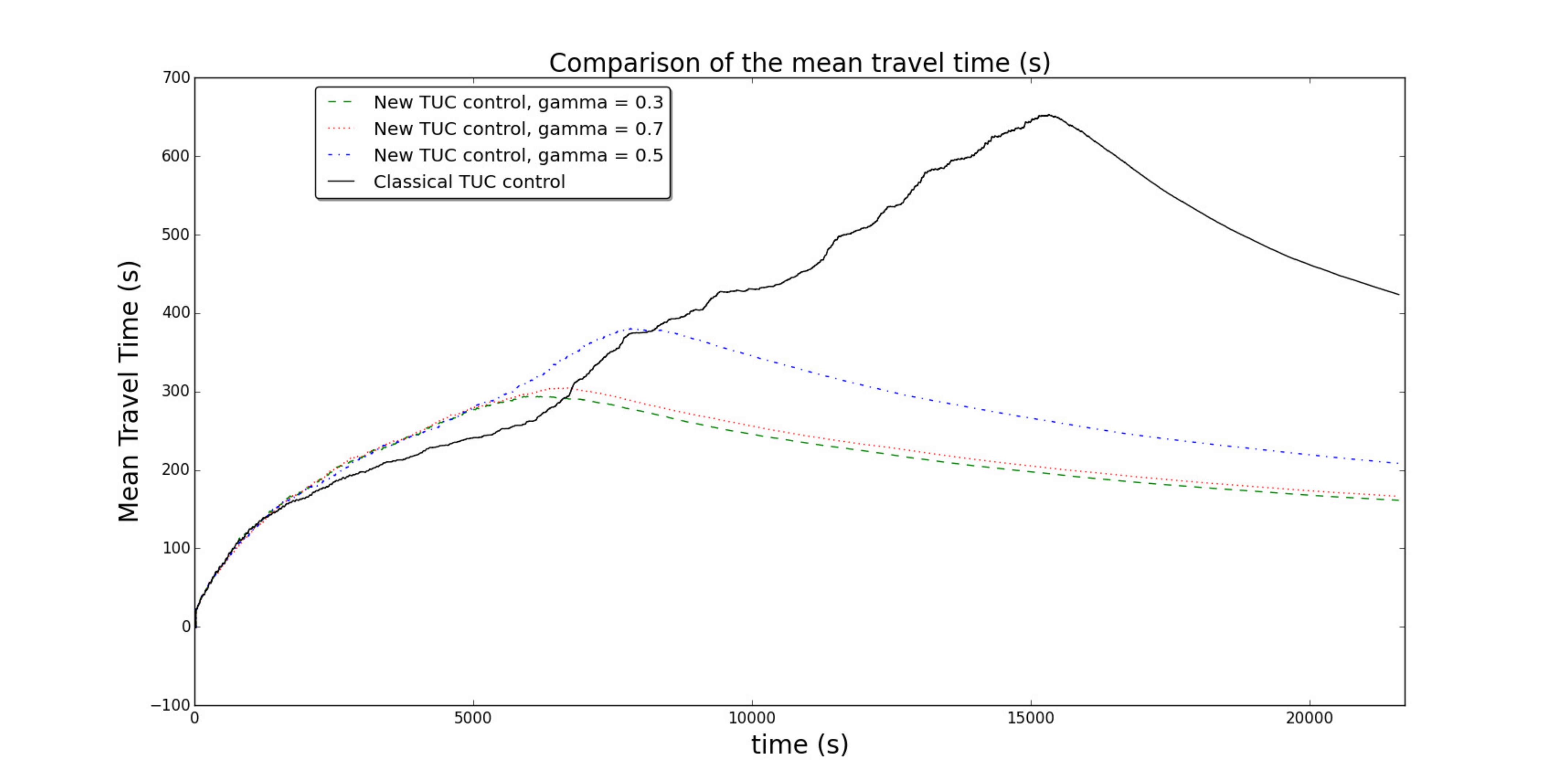} & \includegraphics[width=77mm]{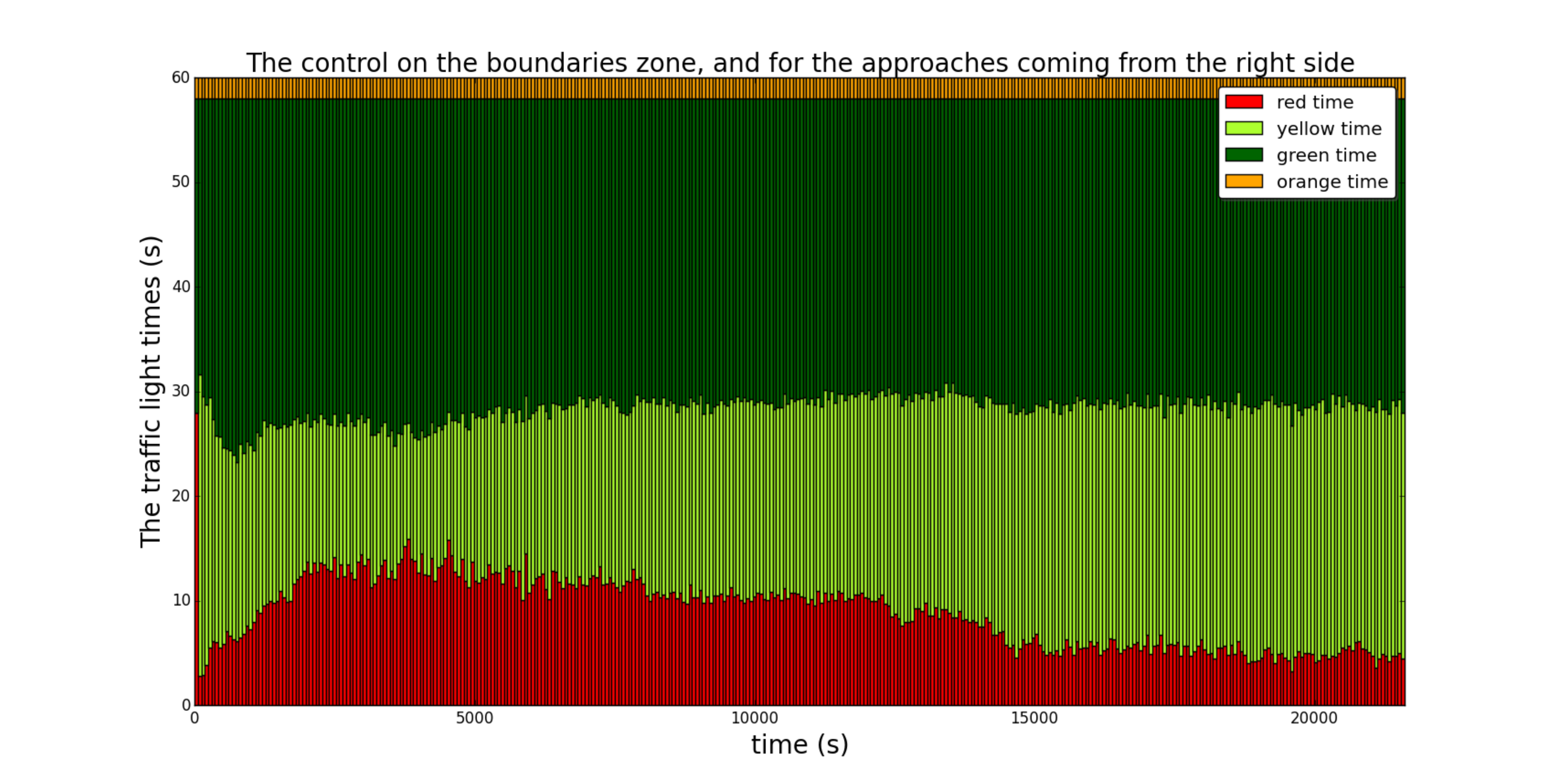}
\end{tabular}
\caption{Left side: Comparison of the classical TUC with the semi-decentralized TUC in terms of the number of running vehicles 
on the network, the flow of ended vehicles, and the average travel time through the network, respectively. 
Right side: the control in terms of the traffic light times into the cycle time.}
 \label{fig-all}  
\end{center}
\end{figure}

In this section, we apply the control model presented above, on a small grid network of Fig.~\ref{figNet}.
We consider here a cycle time of~60 seconds at all the junctions.
We use the traffic simulator SUMO for the traffic simulation with the interface Traci for the implementation of the control.
For more details on the implementation, please refer to~\cite{FNHL15}.
We present the results in Fig.~\ref{figNet} and Fig.~\ref{fig-all}.

On the left side of Fig.~\ref{fig-all}, we compare the classical TUC control with the our control model with different values
of the parameter $\gamma$. The comparison is made basing on three criterion: the number of running cars in the network, the
number of ended (served) cars, and the mean travel time.
We see clearly that our control improves the whole capacity of the network.

On the right side of Fig.~\ref{fig-all}, we give the average number of vehicles in the network (first row),
the control (in term of durations of the green, yellow and red times)
for the approaches coming from the left side of the circuit junctions (second row) and form right side of the circuit junctions (third row). 
An important result here is that the yellow time is almost fully used (i.e. the red time is almost zero)
in the case of free traffic flow, while the red time appears with important values in case of congestion.
This result is very important because it confirms the importance of the fact that the activation as well as the duration of the local control
(the contention time window with yellow times) are both controlled by the centralized control, which optimizes them in function of the
state of the traffic in the network.

\newpage
\section{Optimal-path selection for Vehicular Ad hoc NETworks}

The main reference of this section is~\cite{LSBF18}. We present here a work we have done in collaboration
with Mrs Souaad Lahlah\footnote{Mmrs Souaad Lahlah has passed two short time visitor internships at the 
Grettia Laboratory during the preparation of her PhD thesis.}, a PhD student from the University of Bejaia 
(Algeria) and with her colleagues from the same University.

We proposed a Proactive-optimal-path Selection with Coordinator Agents Assisted Routing (PSCAR) 
protocol for Vehicular Ad hoc NETworks (VANETs) in an urban environment.
We consider static nodes as coordinator agents at every intersection of the urban road network.
All the paths from any coordinator agent to another coordinator agent are known from all the agents,
since they are static. Therefore, the problem of seeking an optimal path to a given destination  node
is reduced to search an optimal path to the coordinator agent which is the nearest one to the destination node,
in such a way that we anticipate the movements of the destination.
Two criteria are used to select the optimal path: the distance and the car-density on the path.
The latter is estimated from the car-speeds on the links of the path, and from given fundamental traffic diagrams 
on the links. We used ns-2 (Network Simulator 2) for the evaluation of the communication performances of PSCAR,
and we based on SUMO (Simulation of Urban MObility) for the road traffic simulation. 
We compare PSCAR with existing solutions in term of packet delivery ratio, end-to-end delay, and network
overhead.

\subsection*{Description of PSCAR}

We assume that the cars are equipped with GPS. Vehicles periodically broadcast to their neighbors Hello-beacon packets reporting their
position coordinates and eventually other information such as speed, acceleration, etc. The nodes build and update neighbor cards basing on all 
the received Hello-bacon packets. On the other side, each Coordinator Agent (a static node deployed at an intersection) broadcasts periodically
a beacon message to announce its presence to its neighbors (vehicles or other Coordinator Agents).
Each coordinator agent maintains a path table with all the possible paths toward any other coordinator Agent.
Optimal paths are selected by the proposed Proactive-optimal-path Selection process, and they are updated periodically depending on the
variation of the car-density in the road network.

\subsection*{The Proactive-optimal-path Selection process}

The proactive-optimal-path selection process includes three selections to forward a data packet.
\begin{enumerate}
  \item from a forwarder node to a Coordinator Agent.
  \item from the first Coordinator Agent to the last Coordinator Agent. In this step, we select the shortest path (in which the vehicular traffic
     is estimated basing on a car-density estimation on the links; see below), and which is sufficiently connected to route the data packets to the destination.
  \item from the last Coordinator Agent to the destination node.
\end{enumerate}

\begin{figure}[htbp]
  \centering
  \includegraphics[scale=0.7]{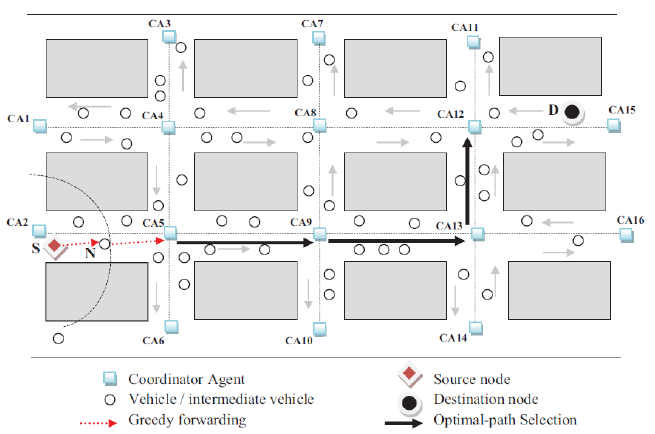} 
  \caption{Optimal-path Selection from a forwarder node to a Coordinator Agent.}
 \label{pscar1}  
\end{figure}

\subsection*{Local optimal and solutions}

In the optimal-path selection process based on the proposed traffic density estimation
approach, a local optimum can be reached instead of a global one. 
For such cases, we propose a recovery strategy which reroutes the data packet to the desired destination.
Two cases are distinguished here.
\begin{itemize}
  \item The source node may not be able to route data packets to the elected coordinator agent, because no sufficiently dense path 
    exists from the source node to the coordinator agent. In this case, the data packets are carrying by the source node until
    such a path exists. 
  \item The same situation can be presented at an intermediate node. In this case we also consider the possibility that
    the last visited coordinator agent selects a new path to the destination node.
\end{itemize}

\subsection*{Path update based on the car-density estimation approach}

\begin{figure}[htbp]
  \centering
  \includegraphics[scale=0.9]{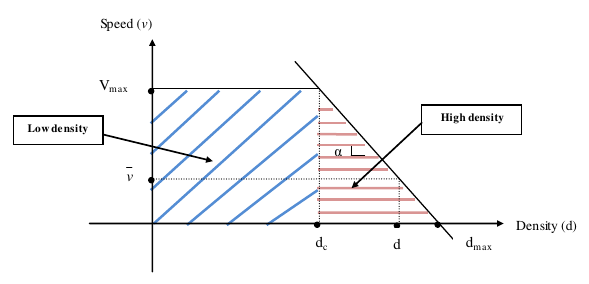} 
  \caption{Fundamental diagram of vehicular traffic.}
 \label{pscar2}  
\end{figure}

We propose here an approach for the estimation of the car-density of every link of the road network.
This estimation is used for the calculus of optimum routing paths of data packets in the road network.
The approach assumes known a fundamental traffic diagram giving the car-speed in function of the car-density
for every link. The coordinator agent at the initial extremity of a link estimates the average car-speed on link
by retrieving all the current car speeds on it, and by that deduces the average car-density on the link by means
of the fundamental diagram of traffic associated to the considered link.
Figure~\ref{pscar2} shows an example of a fundamental traffic diagram considered here. 

\begin{equation}
  v(d) = \min \left\{ v_{\max}, \frac{v_{\max}}{d_{\max} - d_{cr}} (d_{\max} - d)\right\}.
\end{equation}
$v, v_{\max}, d, d_{cr}$ and $ d_{\max}$ denote the average car-speed, the maximum car-speed, the average car-density,
the critical car-density, and the maximum car-density respectively, on the considered link.

\subsection*{Performance evaluations}

The protocol PSCAR proposed here is evaluated by simulation, and is compared to existing ones such as
GPCR (Lochert et al., 2005~\cite{LMFH05}) and A-STAR (Seet et al., 2004~\cite{See04}).
Simulations are performed here by using ns-2 for the simulation of communications and SUMO for the simulation of mobility.
We base our evaluations and comparisons on three main criteria: packet delivery ratio
(the fraction of the data packets successfully delivered to their destination vehicles), end-to-end delay (the average
time that a packet takes to traverse the network from the source to the destination), and network
overhead (the total bytes transmitted during the simulation, including beaconing messages, data packets, data packets
generated during the broadcast mechanism, etc.).

\begin{figure}[htbp]
  \centering
  \includegraphics[scale=0.4]{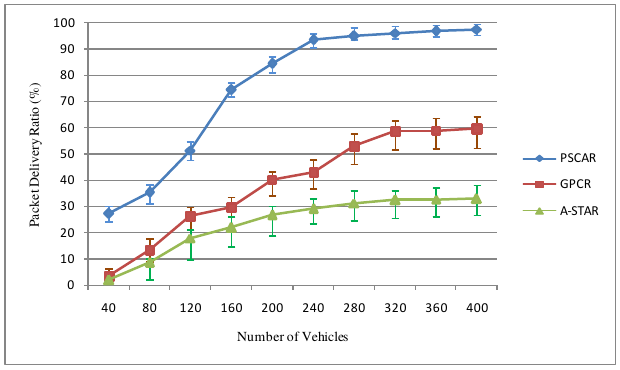} 
  \includegraphics[scale=0.4]{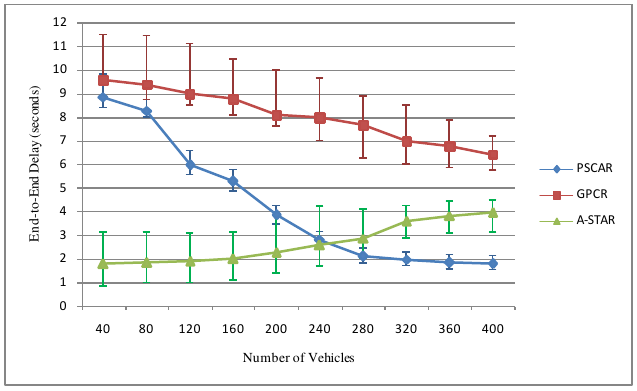}
  \caption{Left side: Packet Delivery Ratio over Number of Vehicles.
    Right side: End-to-End Delay over Number of Vehicles.}
 \label{pscar3-4}  
\end{figure}

\begin{figure}[htbp]
  \centering
  \includegraphics[scale=0.4]{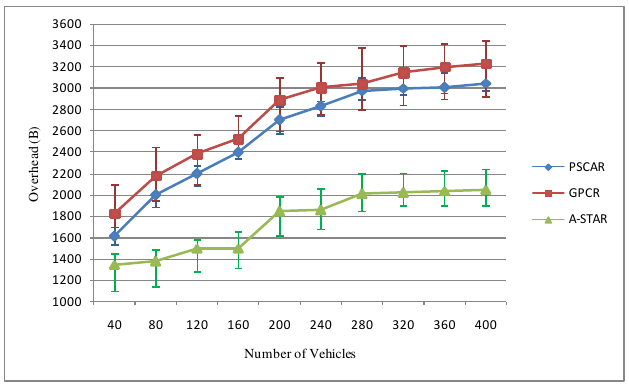}
  \caption{Overhead over Number of Vehicles.}
 \label{pscar5-6}  
\end{figure}

Figure~\ref{pscar3-4} compares PSCAR with GPCR and A-STAR basing on two criteria: the packet delivery ratio (left side)
and the end-to-end delay (right side) given in function of the number of vehicles in the network (car-density).
From the left side of Figure~\ref{pscar3-4}, we see that PSCAR outperforms GPCR and A-STAR in term of packet delivery ration.
The latter is 28\% (respectively 98\%) higher with PSCAR than with the maximum over GPCR and A-STAR, for low 
(respectively high) car-densities.
On the right side of Figure~\ref{pscar3-4}, we show that, in term of end-to-end delay, PSCAR outperforms both GCPR and A-STAR for high car-densities,
while it is outperformed by A-STAR in low densities.  
Indeed, when a local optimal occurs, PSCAR uses the carry and forward strategy to permit to a coordinator agent to choose another path.
This explains its high end-to-end delay at low densities.

On Figure~\ref{pscar5-6} PSCAR is compared to GPCR and A-STAR in term of network overhead 
(the number of beaconing messages generated during the construction of the neighborhood tables).
The figure shows that while PSCAR outperforms GPCR for all the car-densities, it is outperformed by A-STAR
for all the car-density. We consider that higher network overhead obtained here for PSCAR, compared to
A-STAR is a kind of price to pay for its good performances in term of packet delivery ration and end-to-end 
delays. Moreover, We consider that the network overhead performed by PSCAR is acceptable.

\subsection*{Conclusions and perspectives}

We proposed in this section a new position-based routing protocol for VANETs called
Proactive-optimal-path Selection with Coordinator Agent Assisted Routing (PSCAR).
PSCAR introduces a novel car-density estimation approach on the links of the road network,
which permits to assess link connectivity in the network.
Although the car-density estimation generates additional network overhead, PSCAR permits
better performances in term of packet delivery ration and end-to-end delay compared to
existing protocols such as GPCR and A-STAR.
Our perspectives in this direction of research is to improve the beaconing phase, perform
simulations with realistic street maps using for example TIGER database.

\newpage
\section{Simulation of cooperative ITS in urban networks}

The main reference of this section is~\cite{NFHL17}.
We summarize here an algorithm we proposed in~\cite{NFHL17} for urban traffic control including vehicle to infrastructure (v2i)
communications using  WAVE/IEEE 802.11p protocol. Both vehicles and traffic lights may be equipped with communication devices.
We assume mixed traffic with equipped and unequipped vehicles, as well as mixed
infrastructure with equipped and unequipped traffic lights.
The algorithm is first illustrated on an equipped junction, before extending it to a network with
equipped and unequipped junctions.
We used here the open source tool Veins~\cite{Sommer11} which couples the urban traffic simulator SUMO~\cite{sumo} with the 
communication simulator OMNET++~\cite{Omnet}.
We notice that we have modified VEINS in order to include TCP/IP support over IEEE 802.11p.

The algorithm uses the information transmitted by equipped cars to
equipped junctions, on their positions and speeds, in order to improve the split time of the
equipped traffic lights. We have chosen here to keep using traffic lights in order to guarantee 
the traffic safety. 
The main idea of the algorithm is to alternate the green light over the approaches of the equipped junction,
but in the case where no equipped vehicle is manifested from a fixed distance on the approach with green light,
and where at least one equipped vehicle is connected to the traffic light at a short distance on the approach with
red light, the light switches for a short time; and the test is repeated.
By that we ensure that with zero equipped vehicles, we alternate the green light as usual, while with a considerable rate 
of equipped vehicles we may improve the traffic. 
We show that the gain in road traffic performance is significant,
especially in the case with high penetration rates for equipped vehicles and junctions.

\subsection*{The algorithm}

We consider the communication protocol is IEEE 802.11p coupled with the Internet Protocol version~4 (IPv4)
and the Transmission Control Protocol (TCP). We assume that equipped junctions and cars are able to use those communication
capabilities. Moreover, we suppose that equipped cars are also able to localize themselves, for example with
GPS modules.
Let us present our new algorithm for urban traffic control at an equipped junction.
The algorithm includes three main steps: building a map, electing a vehicle and actuating the traffic light signal.

\textbf{Build a dynamic map.} 
An equipped junction builds a map of the connected vehicles coming and leaving the
junction. Similarly, an equipped vehicle builds a map of the equipped junctions approaching or leaving it in its communication range.
The built maps are dynamic and are updated periodically each time a message is received for the equipped junction, and triggered on
timer for the equipped vehicles.

\textbf{Vehicle election.}
An equipped junction elects periodically the lead vehicles on the approaching edges.
For simplicity and without loss of generality, we assume that the junctions have only two incoming edges, with one lane each.
The two incoming edges alternate priority of passing through the junction every $c/2$ time units where $c$ is the cycle time
of the traffic light. At most one vehicle among lead vehicles from incoming edges is elected with Algorithm~\ref{algor1} below. 
If the latter does not succeed to elect a vehicle, it tries again after a given time interval.

\begin{algorithm}
\label{algor1}
\caption{Vehicle Election}
\label{election}    
    \SetKwInOut{Input}{Input}
    \SetKwInOut{Output}{Output}
    \underline{function Elect} $(p,v,d_p,d_v,d_{min},\alpha)$\\
    \Input{ \begin{itemize}
	      \item $p$ is the identifier of the lead vehicle on the prioritized edge, and it is $None$ if no vehicle is detected on the prioritized edge
	      \item $v$ is the identifier of the lead vehicle on the non prioritized edge, and it is $None$ if no vehicle is detected on the non prioritized edge
	      \item $d_p$ represents the distance $p$ is to the junction, in case $p\neq None$,
	      \item $d_v$ represents the distance $v$ is to the junction, in case $v\neq None$,
	      \item $d_{min}>0$ is the minimum  distance to consider a vehicle close to the junction,
	      \item $\alpha > 1$ is a coefficient to weight the minimum distance.
	    \end{itemize}
      }
    \Output{ $p$ or $v$.}

    \uIf{($p \neq None$ and $v \neq None$ and $d_p > \alpha d_{min}$  and $d_v < d_{min}$) or ($p == None$ and $v \neq None$ and $d_v < d_{min}$)}
    {
	$elected=v$\;

    }
    \uElse{
        $elected=p$\;
    }
    return $elected$\;
    
\end{algorithm}

\textbf{Traffic light actuation.} 
The elected equipped vehicle sends a message to the equipped junction with its established TCP connection,
and sets the traffic light to green color for the edge on which it is moving.
A minimum duration for a given traffic light state is considered, during which 
the traffic light state has to remain unchanged to ensure stability.
A maximum duration is also considered, during which if no traffic light state switch has occurred, then
the state of the traffic light is automatically changed. This permits to avoid blockages.

We give below some properties of the proposed algorithm.
\begin{itemize}
  \item The algorithm assures safety since it works with traffic lights
    which never give green light to antagonistic stages.
  \item It is not necessary for a vehicle to be equipped to pass through a junction.
  \item Equipped vehicles are better served, and equipped junctions are better managed.
  \item In case of unequipped junctions, and in case of equipped junctions with no equipped vehicles,
    the algorithm is equivalent to an open loop control.
\end{itemize}

\subsection*{Simulation results}

We present statistical results of some runs with different seeds, first on an equipped junction, and then on a
small grid network of equipped and unequipped junctions.
For the scenario of one junction, we assume a junction of two incoming links with one lane each.
The traffic demand is varied as well as the ratio of equipped vehicles.
We simulate the vehicular traffic as well as the communications for 600 seconds.

The criteria we used here for performance evaluation of the communication are:
TCP data divided by simulation time, mean TCP end-to-end delay, and mean TCP throughput on road side unit.
The results are given on the first column of Tab.~\ref{tab-com}.

\begin{table}[htbp]
\centering
\begin{tabular}{|c|c|}
  \hline
   & \\
  \includegraphics[scale=0.72]{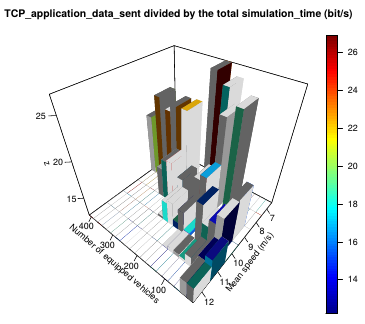} & \includegraphics[scale=0.72]{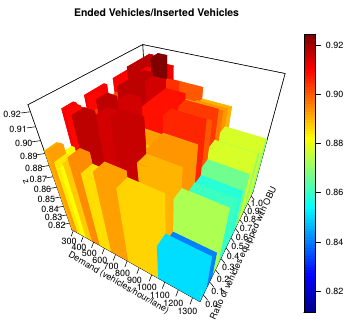} \\
  \hline
  & \\
  \includegraphics[scale=0.72]{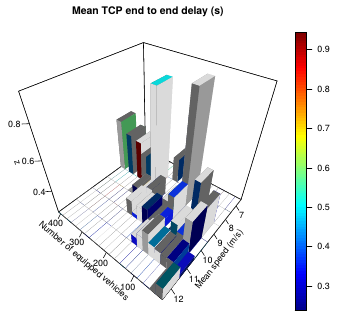} & \includegraphics[scale=0.72]{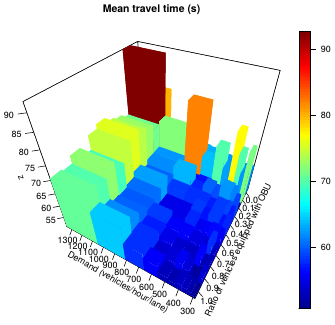} \\
  \hline
  & \\
  \includegraphics[scale=0.72]{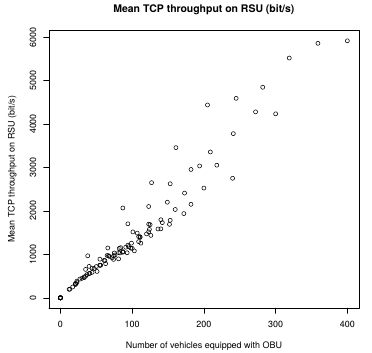} & \includegraphics[scale=0.72]{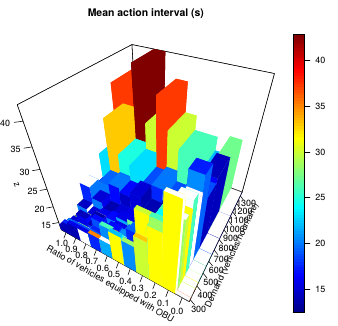} \\
  \hline
\end{tabular}
\caption{Results of simulation on one equipped junction.}
\label{tab-com}
\end{table}

From Table~\ref{tab-com} we see that the amount of TCP application data sent by the nodes increases as the mean vehicle speed
decreases. We ca interpret this observation by the fact that  as the mean vehicle speed is low, the communicating vehicles 
remain connected longer, and then they send more messages.
We can also see that the mean TCP end-to-end delay can be as high as 0.8 s. when vehicle speed is low (about 8 m/s).
Finally, the throughput on RSU is increasing linearly with the number of communicating vehicles.

The criteria used for the vehicular traffic are the ratio of the ended vehicles by the inserted vehicles, and the mean travel time.
The results are given on the second column of Tab.~\ref{tab-com}.
We see on Table~\ref{tab-com} that the ratio of the ended vehicles by the inserted vehicles increases when the demand decrease and when
the ratio of equipped vehicles increases. For a fixed demand, the mean travel time decreases as the number of equipped vehicles increases.
For the actuator, we define the mean action interval as the mean time interval between two consecutive changes of traffic lights state.
We observe on Table~\ref{tab-com} that the traffic light state is stable with a high traffic demand and a high equipped vehicle penetration rate.
while for a low demand and a high equipped vehicle penetration rate, the control is more reactive.

For the results on a small grid network, we consider 
a network of 4~horizontal and 4~vertical one-lane roads crossing at 16~junctions.
The direction of the traffic on the horizontal and on the vertical roads is alternated
from one road to another (first road traffic goes from east (north) to west (south), second road from
west (south) to east (north), and so on).

The results of simulation on the grid network are given in Table~\ref{tab:saturated1}, where
we compare our algorithm with an open loop fixed cycle traffic light program.
We notice that the open loop control case can also be realized by fixing the rate of
equipped vehicles to zero.
We see from Table~\ref{tab:saturated1} that with 100\% of equipped junction and 80\% of equipped vehicles,
the gain of our algorithm in terms of ended vehicles (resp. running vehicles and mean travel time)
can exceed 35\% (resp. 45\% and 35\%). 
A particular case is observed here for 20\% of equipped vehicles combined with 100\% of equipped
junctions, where the algorithm seems to be non efficient. 
We think that the small rate of equipped vehicles disturbs the control algorithm at every junction
and then deteriorates the traffic.

\begin{table}[htbp!]
      \begin{center}
      \caption{The results of simulation on the grid network.} 
      \label{tab:saturated1}
      \resizebox{\columnwidth}{!}{%
      
      \begin{tabular}{|c|c|c|c|c|}
      \hline 
      \diagbox{Equipped junctions}{Ended}{Penetration rate}&0\%   & 20\% & 50\% & 80\% \\
      \hline
      25\%& 1373$\pm$19 &1470$\pm$33& 1507$\pm$18 & 1484$\pm$19 \\ 
	  & (0$\pm$0)\% & \textbf{(+7.1$\pm$2.5)\%} & \textbf{(+9.8$\pm$2.4)\%} & \textbf{(+8.1$\pm$1.9)\%} \\
      \hline
      50\%& 1373 $\pm$19 & 1499$\pm$49 & 1583$\pm$19 & 1571$\pm$20\\		
	  & (0$\pm$0) \% & \textbf{(+9.2$\pm$3.4)\%} & \textbf{(+15.3$\pm$1.9)\%} & \textbf{(+14.5$\pm$2.3)} \%\\       
      \hline    
      100\%& 1373$\pm$19 & 1281$\pm$151 & 1805$\pm$49 & 1877$\pm$29\\		   
	  & (0$\pm$0)\% & \textcolor{red}{(-6.7$\pm$11.2)\%} & \textbf{(+31.5$\pm$4.1)\%} & \textcolor{green}{(+36.7$\pm$2.8)\%} \\
      \hline
      \hline 
      \diagbox{Equipped junctions}{Running}{Penetration rate}&0\%   & 20\% & 50\% & 80\% \\
      \hline
      25\% & 954$\pm$16 & 842$\pm$32 & 817$\pm$18 & 848$\pm$21 \\		     
      &(0$\pm$0\%) & \textbf{(-11.8$\pm$3.5\%)} & \textbf{(-14.4$\pm$2.7\%)} & \textbf{(-11.2$\pm$2.4\%)} \\
      \hline
      50\%  & 954$\pm$16 & 835$\pm$36 & 764$\pm$17 & 778$\pm$21 \\		     
      & (0$\pm$0\%) & \textbf{(-12.4$\pm$3.8\%)} & \textbf{(-19.9$\pm$2.1\%)} & \textbf{(-18.4$\pm$2.5\%)} \\
      \hline
      100\% & 954$\pm$16 & 962$\pm$80 & 583$\pm$43 & 517$\pm$27 \\		   
      & (0$\pm$0\%) & \textcolor{red}{(+0.9$\pm$8.7\%) } & \textbf{(-38.9$\pm$4.7\%)} & \textcolor{green}{(-45.8$\pm$2.9\%) } \\
      \hline
      \hline 
      \diagbox{Equipped junctions}{MTT(s)}{Penetration rate}&0\%   & 20\% & 50\% & 80\% \\
      \hline
      25\% & 413.9$\pm$1.8 & 381.9$\pm$4.3 & 376.0$\pm$3.1 & 380.0$\pm$3.4 \\		     
      & (0$\pm$0\%) & \textbf{(-7.7$\pm$1.1\%)} & \textbf{(-9.2$\pm$0.8\%)} & \textbf{(-8.2$\pm$0.9\%)} \\
      \hline
      50\% & 413.9$\pm$1.8 & 381.8$\pm$15.2 & 355.6$\pm$5.0 & 354.9$\pm$4.5 \\		     
      & (0$\pm$0\%) & \textbf{(-7.8$\pm$3.8\%)} & \textbf{(-14.1$\pm$1.5\%)} & \textbf{(-14.2$\pm$1.2\%)} \\
      \hline
      100\% & 413.9$\pm$1.8 & 399.6$\pm$30.6 & 302.0$\pm$6.9 & 281.2$\pm$4.6 \\		   
      & (0$\pm$0\%) & \textcolor{red}{(-3.4$\pm$7.4\%) } & \textbf{(-27.0$\pm$1.7\%)}  & \textcolor{green}{(-32.1$\pm$1.2\%) } \\
      \hline
      \end{tabular}
      
      }
      \end{center}
      \label{tab-net}
\end{table}


\newpage
This page is intentionally left blank

\part{Stochastic dynamic programming based modeling}
\label{part-stoch}

We present in this part our works related to the traffic modeling and control in
transportation networks, based on stochastic dynamic-programming systems.
The traffic models we present here are mainly based on dynamic programming systems associated
to stochastic optimal control of Markov chains.
This part is organized in four chapters.

Chapter~\ref{chap-stoch} provides an introduction to stochastic dynamic programming and its relationship with dynamic systems.
It reviews some existing results on the existence of asymptotic regimes of additive 1-homogeneous and monotone dynamic systems.
It focuses in particular on dynamic systems written in the form of stochastic dynamic programming systems associated to
optimal control problems of Markov chains. 

Chapter~\ref{chap-anticip} proposes a multi-anticipative vehicular traffic model.
The model is microscopic and is based on an existing model (\textit{the piecewise linear car-following model}~\cite{Far12}).
We extended the model of~\cite{Far12} in such a way that the drivers can take into account
more than one leader vehicle in their car-following behavior laws.
We show that the car dynamics are written as a stochastic dynamic programming system associated to optimal control of a Markov chain.
We give the interpretation of that in terms of vehicular traffic.
We derive analytically the macroscopic behavior emerging from the microscopic dynamics, by obtaining the fundamental traffic diagrams
giving the asymptotic traffic flow as a function of the average car-density.
Moreover, We show that the traffic is more smooth in the case where the drivers have such multi-anticipative behaviors.

Chapter~\ref{chap-stochtrain} proposes many extensions of the Max-plus algebra model of the train dynamics presented in Chapter~\ref{chap-maxtrain}.
We first present train dynamics where the effect of the passenger travel demand on the train dwell times is modeled naturally without any control.
In this case, a high level of the passenger travel demand at a given platform extends naturally the train dwell time at the platform.
We show that the dynamics are unstable in this case, and give interpretation of that fact in term of the train dynamics and of the effect
of the passenger demand on it.
We then propose a second model with a control law on the train dwell times, in such a way as to guarantee the stability of the train dynamics.
We show that the controlled train dynamics is written is this case as a dynamic programming system of optimal control problem of a Markov chain,
and give the interpretation of that in term of traffic. We show analytically that the train dynamics admit a stationary regime with a unique
average growth rate, interpreted here as the average asymptotic train time-headway.
We derive then by numerical simulation the asymptotic average train frequency as a function of the number of running trains on the metro line,
and of the passenger demand level; and by that, we derive the effect of the passenger travel demand on the train dynamics.

Finally, Chapter~\ref{chap-otherstoch} summarizes some of my other contributions on stochastic modeling and control of traffic.
Two works are briefly presented. 
First, we present an extension of an existing model for optimal routing in networks.
The objective of our extension is to introduce robustness of the routing against link failure, with an application for the guidance of the users of road networks.
A routing strategy is said to be robust here if the deterioration of its maximum value (probability of arrival) is minimized in case of link failure.
Second, we present a macroscopic model for multi-lane vehicular traffic, where we model the assignment and the dynamic reassignment of 
the traffic to different lanes of a highway stretch.


\chapter[Stochastic dyn. prog. - based systems]{Introduction to and reviews on stochastic dynamic programming - based systems}
\label{chap-stoch}

We give in this chapter some necessary reviews and theoretical tools for the models we present in the chapters of this part.
Our models are based on the dynamic programming systems associated to stochastic optimal control problems of Markov chains.
The reviews are organized in three sections.

First, we consider dynamic systems of the form 
$x(k) = f(x(k-1))$, where $f: \mathbb R^n \to \mathbb R^n$
is an additive homogeneous of degree~1 and monotone map.
We recall that under a kind of reducibility condition for $f$, the dynamic system admits a stationary regime with a unique
average growth rate. We give an extension of this result to dynamic systems of the form
$$x(k) = f(x(k), x(k-1), \ldots, x(k-m+1)),$$
where $f:\mathbb R^{n\times m} \to \mathbb R^n$.

Second, we consider the dynamics
\begin{equation}\label{dyn_stoch1}
  x(k+1)_i = \max_{u\in\mathcal U} \left( (M^u x(k))_i + c^u_i \right), \forall i \in\{1,2,\ldots, n\},
\end{equation}  
which are dynamic programming systems associated to control problems of Markov chains, with transition matrices $M^u$ and reward vectors $c^u$, for $u\in\mathcal U$.
We give some reviews on these dynamics, and then consider the extended dynamics 
$$x(k)_i = \max_{u\in\mathcal U} \left( (M^u x(k-1))_i + (N^u x(k))_i + c^u_i \right), \forall i \in\{1,2,\ldots, n\},$$
with possible implicit terms; and provide some results. 

Third, we consider another extension to the dynamic system~(\ref{dyn_stoch1}) to dynamic systems
associated to stochastic games:
$$x(k+1)_i = \min_{u\in\mathcal U}\max_{w\in\mathcal W} \left( (M^{uw} x(k))_i + c^{uw}_i \right), \forall i \in\{1,2,\ldots, n\}.$$
We give the necessary results of these dynamics, for our models.

\section{Non-expansive dynamic systems}
\label{sec-gdps}
We are concerned here with maps $\mathbf{f} : \mathbb R^n \to \mathbb R^n$.
Such a map is additive 1-homogeneous if
$$\forall x \in \mathbb R^n, \forall a \in \mathbb R, \mathbf{f}(a \textbf{1} + x) = a\textbf{1} + \mathbf{f}(x),$$
where $\textbf{1} \stackrel{def}{=} {}^t(1, 1, \ldots , 1)$. 
$\mathbf{f}$ is monotone if
$$\forall x, y \in \mathbb R^n, x \leq y \Rightarrow \mathbf{f}(x) \leq \mathbf{f}(y),$$
where $x \leq y$ means $x_i \leq y_i \forall i, 1 \leq i \leq n$.
If $\mathbf{f}$ is 1-homogeneous and monotone, then it is non-expansive (or 1-Lipschitz) for the sup. norm~\cite{CT80}, i. e.
$$\forall x, y \in \mathbb R^n, ||\mathbf{f}(x) - \mathbf{f}(y)||_{\infty} \leq ||x - y||_{\infty}.$$

A directed graph $\mathcal G(\mathbf{f})$ is associated to a non-expansive map $\mathbf{f}$.
$\mathcal G(\mathbf{f})$ is defined~\cite{GG99} by the
set of nodes $\{1, 2, \ldots,n\}$ and by a set of links such that there exists a link from
a node $i$ to a node $j$ if $\lim_{\eta\to\infty} \mathbf{f}_i(\eta e_j ) = \infty$, where $e_j$ is the $j$th
vector of the canonical basis of $\mathbb R^n$. 

We recall below an important result on the existence of asymptotic regimes and asymptotic growth rates 
for non-expansive dynamic systems (i.e. dynamic systems defined with non-expansive maps).
\begin{theorem}\cite{GG98b, GK95}\label{th-dps}
  If $\mathbf{f} : \mathbb R^n \to \mathbb R^n$ is 1-homogeneous and monotone and if $\mathcal G(\mathbf{f})$
  is strongly connected then $\mathbf{f}$ admits a unique additive eigenvalue, i.e.
  $\exists \! \mu \in \mathbb R, \exists x \in \mathbb R^n : \mathbf{f}(x) = \mu + x$.
  Moreover, the dynamic system $x(k+1) = \mathbf{f}(x(k))$ admits a stationary regime, with a unique 
  asymptotic average growth vector $\chi(\mathbf{f})$
  defined $\chi (\mathbf{f}) := \lim_{k \to \infty} \mathbf{f}^k(x)/k$. Furthermore, we have 
  $\chi (\mathbf{f}) = \mu \textbf{1}$.
\end{theorem}
 
For simplicity, we use in this article, for all the dynamic systems, the notation $x^k$ instead of $x(k)$.
We give in the following, a natural extension of Theorem~\ref{th-dps}, which will permit us to
consider dynamic systems of the form $x^k = \mathbf{f}(x^k,x^{k-1},\ldots,x^{k-m+1})$.
For that, let us consider $\mathbf{f}: \mathbb R^{m\times n} \to \mathbb R^n$
associating for  $(x^{(0)},x^{(1)},\cdots, x^{(m-1)})$, where $x^{(i)}$ are vectors in $\mathbb R^n$,
a column vector in $\mathbb R^n$.
$$\begin{array}{llll}
  \mathbf{f}: & \mathbb R^{m\times n} & \to & \mathbb R^n \\
     & (x^{(0)},x^{(1)},\cdots, x^{(m-1)}) & \mapsto & \mathbf{f}(x).
\end{array}$$
We denote by $\mathbf{f}_{(i)}, i=0,1,\ldots,m-1$, the following maps.
$$\begin{array}{llll}
  \mathbf{f}_{(i)}: & \mathbb R^{n} & \to & \mathbb R^n \\
     & x & \mapsto & \mathbf{f}_{(i)}(x) = \mathbf{f}(-\infty, \ldots,-\infty,x,-\infty, \ldots, -\infty).\\
     &   &         & \qquad \; \qquad \qquad \qquad \qquad \; \; \; \uparrow \\
     &   &         & \quad \qquad \qquad \qquad \qquad \; \; \; i^{\text{th}} \text{component} 
\end{array}$$
and by $\tilde{\mathbf{f}}$ the following map.
$$\begin{array}{llll}
  \tilde{\mathbf{f}}: & \mathbb R^{n} & \to & \mathbb R^n \\
             & x & \mapsto & \tilde{\mathbf{f}}(x) = \mathbf{f}(x,x, \ldots,x).
\end{array}$$

\begin{theorem}\label{th-dps2}
   If $\tilde{\mathbf{f}}$ is additive 1-homogeneous and monotone, and if $\mathcal G(\tilde{\mathbf{f}})$
   id strongly connected and $\mathcal G(\mathbf{f}_{(0)})$ is acyclic, then $\mathbf{f}$ admits a unique generalized
   additive eigenvalue $\mu > -\infty$ and an additive eigenvector $v > -\infty$, such that
   $\mathbf{f}(v, v-\mu, v-2\mu, \ldots, v-(m-1)\mu) = v$.
   Moreover, $\chi(\mathbf{f}) = \mu \textbf{1}$.
\end{theorem}

\proof
The proof consists in showing that the dynamic system $x^k = \mathbf{f}(x^k,x^{k-1},\ldots,x^{k-m+1})$
is equivalent to another dynamic system $z^k = \mathbf{h}(z^{k-1})$, where $\mathbf{h}$ is built from $\mathbf{f}$,
such that $\mathbf{h}$ satisfies additive 1-homogeneity, monotonicity and connectivity properties needed by Theorem~\ref{th-dps}.
We give here a sketch of the proof. 
Two steps are needed for the proof.
\begin{enumerate}
  \item Eliminate the dependence of $x^k$ on $x^{k-2}$, $x^{k-3}$, $\ldots$, $x^{k-m+1}$. 
    This is done by the well known state augmentation technique.
  \item Eliminate implicit terms (the dependence of $x^k$ on $x^k$).
    This is done by defining an order of updating the $n$ components of $x^k$, in
    such a way that no implicit term appears. This is possible because $\mathcal G(\mathbf{f}_{(0)})$ 
    is acyclic.
\end{enumerate}
\endproof

\section{Dynamic programming systems of stochastic optimal control problems}
\label{sec-dps}

We review here a particular non-expansive dynamic system encountered
in stochastic optimal control of Markov chains. The dynamic system is written
\begin{equation}\label{eq-dp2}
	  x_i^{k+1} = \max_{u\in\mathcal U} ([M^{u} x^k]_i + c^{u}_i), \quad \forall 1 \leq i \leq n.
\end{equation}
where $\mathcal U$ is a set of indexes corresponding to control actions; $M^{u}$,
for $u\in \mathcal U$, are stochastic matrices (i.e. satisfy $M^{u}_{ij} \geq 0$ and $M^{u} \textbf{1} = \textbf{1}$);
and $c^{u}$, for $u\in \mathcal U$, are reward vectors in $\mathbb R^n$.

(\ref{eq-dp2}) is the dynamic programming system associated to the stochastic optimal control of a Markov chain
with state $z\in \mathcal Z :=\{1,2,\ldots,n\}$, transition matrices $M^u, u\in\mathcal U$, associated rewards
$c^u_z, u\in\mathcal U, z\in\mathcal Z$, and final rewards $\phi_z, z\in\mathcal Z$.
\begin{equation}
  \max_{s\in\mathcal S} \mathbb E \left\{ \sum_{k=0}^{T-1} c_{z^k}^{u^k} + \phi_{z^T}\right\},
\end{equation}
where $\mathcal S$ is the set of feedback strategies $s: \mathcal Z \to \mathcal U$.

The variable $x_i^k, i=1,2,...,n, k\in\mathbb N$ is interpreted, in this case, as the function value
of the stochastic control problem.	

In the traffic models of the train dynamics we propose in Chapter~\ref{chap-stochtrain}, we are concerned with system~(\ref{eq-dp2}) above.
More precisely, we will consider dynamic systems with implicit terms~(\ref{eq-dp3}).
\begin{equation}\label{eq-dp3}
  x_i^{k} = \max_{u\in\mathcal U} ([M^{u} x^{k-1}]_i + [N^{u} x^k]_i + c^{u}_i), \quad \forall 1 \leq i \leq n,
\end{equation}
where $M^u, u\in\mathcal U$ and $N^u, u\in\mathcal U$ satisfy $M^u_{ij} \geq 0, N^u_{ij} \geq 0$, 
and $\sum_j (M^u_{ij}+N^u_{ij}) = 1, \forall i,u$. 
For the analysis of dynamic system~(\ref{eq-dp3}) we will use Theorem~\ref{th-dps2}.

We notice here that, in our traffic models, the interpretation of systems~(\ref{eq-dp2}) and~(\ref{eq-dp3}) will be different 
from the stochastic optimal control one, in the sense that, in our models, such systems model directly
the train dynamics, and are not derived from stochastic optimal control problems.

\section{Dynamic programming systems of stochastic games}
\label{sec-dpg}

We review here another particular non-expansive dynamic system encountered
in stochastic games on a controlled Markov chain, with two players (a minimizer and a maximizer). The dynamic system is written
\begin{equation}\label{eq-dpg1}
	  x_i^{k+1} = \min_{u\in\mathcal U}\max_{w\in\mathcal W} ([M^{uw} x^k]_i + c^{uw}_i), \quad \forall 1 \leq i \leq n.
\end{equation}
where $\mathcal U$ is a set of indexes corresponding to the control actions of the minimizer;
$\mathcal W$ is a set of indexes corresponding to the control actions of the maximizer; $M^{uw}$,
for $u\in \mathcal U, w\in\mathcal W$, are stochastic transition matrices (with $M^{uw}_{ij} \geq 0$ and $M^{uw} \textbf{1} = \textbf{1}$)
of the Markov chain;
and $c^{uw}$, for $u\in \mathcal U, w\in\mathcal W$, are pay-off vectors in $\mathbb R^n$, indexed by the controls of 
both the minimizer and the maximizer. 

(\ref{eq-dpg1}) is the dynamic programming system associated to the stochastic game on a controlled Markov chain
with state $z\in\mathcal Z := \{1,2,\ldots,n\}$, transition matrices $M^{uw}, u\in\mathcal U, w\in\mathcal W$, associated 
pay-offs $c^{uw}_z, u\in\mathcal U, w\in \mathcal W, z\in\mathcal Z$, and final pay-offs $\phi_z, z\in\mathcal Z$.
\begin{equation}
  \min \max |_{s\in \mathcal S} \mathbb E \left\{ \sum_{k=0}^{T-1} c_{z^k}^{u^k w^k} + \phi_{z^T}\right\},
\end{equation}
where $\mathcal S$ is the set of feedback strategies $s: \mathcal Z \to \mathcal U\times \mathcal W$.

The variable $x_i^k, i=1,2,...,n, k\in\mathbb N$ is interpreted, in this case, as the function value of the stochastic control game.

In the traffic model of the car dynamics we propose in Chapter~\ref{chap-anticip} the interpretation of system~(\ref{eq-dpg1}) will be different 
from the stochastic optimal control one, in the sense that, in our model, the dynamic system models directly
the car dynamics, and is not derived from a stochastic game problem.


\chapter[Multi-anticipation modeling]{Multi-anticipation piecewise-linear car-following modeling}
\label{chap-anticip} 

The main reference of this chapter is~\cite{FHL13}.
It consists in an extension of a microscopic traffic model we developed in~\cite{Far08,Far12}, called 
\textit{piecewise linear car-following}.
The extension permits to consider anticipation in driving. That is to say that the car-following law can take into
account more than one leader in the calculus of the reaction to stimulus.
By means of variational formulation, we characterize stability and stationary regimes to the model proposed here.
It is shown that, the multi-anticipative model realizes the same macroscopic behavior as the non-anticipative model,
in term of the stationary regime.
Nevertheless, in the transient traffic, the variance in car-velocities and accelerations is reduced in the case
of multi-anticipative driving, and the car-trajectories are smoothed.
  
In the sequel of this chapter, time is denoted by $t$, distance (car position) is denoted by $x$, and
$n$ numbers the cars, starting by the leader. We then consider the following notations.
\begin{itemize}
  \item $x(n,t)$: the cumulative traveled distance of car $n$, from time zero to time $t$.
  \item $y(n,t)$: the inter-vehicular distance $x(n-1,t)-x(n,t)$.
  \item $v(n,t)$: the velocity of car $n$ at time $t$.
\end{itemize}

In order to situate the model we present here with respect to the existing multi-anticipative models, 
and to understand the extension we do, let us first give a short review on multi-anticipative car-following models.

A straightforward multi-leader extension of the model of Chandler et al.~\cite{CHM58} is the Bexelius model~\cite{Bex68}
\begin{equation}\label{bexelius}
  \dot{v}(n,t+T)=\sum_{j=1}^m \alpha_j \Delta v^{(j)} (n,t),
\end{equation}
where $\dot{v}$ denotes the acceleration, $\alpha_j, j=1,2,\cdots,m$ are sensitivity parameters with respect to
the $j^{\text{th}}$ car ahead, and where $\Delta v^{(j)} (n,t)=v(n-j,t)-v(n,t)$.
This model is very simple but permits some mathematical analysis.

Hoogendoorn et al.~\cite{HOS06} have noted the non convenience of the additive form of Bexelius model~(\ref{bexelius}),
and proposed the modification
\begin{equation}\label{HoogBexel}
  \dot{v}(n,t+T)=\min_{1\leq j\leq m} \tilde{\alpha}_j \Delta v^{(j)}(n,t).
\end{equation}

Hoogendoorn et al.~\cite{HOS06} have also proposed a multi-anticipative generalization for the Helly model~\cite{Hel59,Hel61}
\begin{equation}\label{HoogHelly}
  \dot{v}(n,t+T)=\sum_{j=1}^{m_1} \alpha_j \Delta v^{(j)}(n,t) + \sum_{j=1}^{m_2} \beta_j [\Delta x^{(j)}(n,t)-S^j(n)],
\end{equation}
where $\Delta x^{(j)}(n,t)=x(n,t)-x(n-j,t)$.

Lenz et al.~\cite{LWS99} have generalized the Bando model~\cite{BHNS95} as follows.
\begin{equation}\label{lenz}
  \dot{v}(n,t)=\sum_{j=1}^m \kappa_j\left\{V\left(\frac{\Delta x^{(j)}(n,t)}{j}\right)-v(n,t)\right\},
\end{equation}
where $\kappa_j$ expresses the sensitivity of the $j$th leader.

\section{Review on the piecewise-linear car-following modeling}

We base here on the piecewise linear car-following model proposed in~\cite{Far12},
where the behavioral law $V_e$ giving the velocity of a vehicle as a function of the inter-distance $y$
between the vehicle and its leader, is approximated with a (min-max)-piecewise linear curve:
\begin{equation}\label{approx1}
  V_e(y)=\min_{u\in\mathcal U}\max_{w\in\mathcal W}\{\alpha_{uw} y + \beta_{uw}\},
\end{equation}
where $\alpha_{uw}$ and $\beta_{uw}$, for $(u,w)\in\mathcal U\times \mathcal W$, are parameters, and $\mathcal U$ and $\mathcal W$ are
two finite sets of indexes.
In~(\ref{approx1}), only one leader is considered.
The position of the considered car is then updated as follows.
\begin{equation}\label{Far1}
  x(n,t+1)=x(n,t)+ \min_{u\in\mathcal U}\max_{w\in\mathcal W} \{\alpha_{uw} (x(n-1,t)-x(n,t)) + \beta_{uw}\},
\end{equation}
System~(\ref{Far1}) is also written, for the traffic of $\nu$ cars $1,2, \cdots, \nu$ as follows.
\begin{equation}\label{plcf1}
  x_n(t+1)=\min_{u\in\mathcal U}\max_{w\in\mathcal W} \{[M^{uw} x(t)]_n + c^{uw}_n\}, \quad 1\leq n\leq \nu
\end{equation}
where $M^{uw}$ and $c^{uw}$, for $(u,w)\in\mathcal U\times\mathcal W$, are matrices and column-vectors respectively.

The dynamics~(\ref{plcf1}) have been interpreted in~\cite{Far12}, under the assumption
$\alpha_{uw}\in [0,1], \forall (u,w)\in\mathcal U\times\mathcal W$, as a dynamic programming system
associated to a stochastic game on a controlled Markov chain; see~\cite{Far12} for more details.

Two cases have been distinguished in~\cite{Far12}. 
\begin{itemize}
  \item The $\nu$ cars move on a ring road. In this case, $M^{uw}$ and $c^{uw}$ are given by
    $$M^{uw}=\begin{bmatrix}
        1-\alpha_{uw} & 0 & \cdots & \alpha_{uw}\\
        \alpha_{uw} & 1-\alpha_{uw} & & 0\\
        \vdots & \ddots & \ddots & \\
        0 & 0 & \alpha_{uw} & 1-\alpha_{uw}
      \end{bmatrix},$$
    and
    $$c^{uw}={}^t[\alpha_{uw}\nu/d+\beta_{uw},\;\; \beta_{uw},\; \cdots,\; \beta_{uw}],$$
    where $d$ denotes the car density on the road.
    
    The dynamics~(\ref{plcf1}) are stable under the condition $\alpha_{uw}\in[0,1]$, and the behavior law is realized
    at the stationary regime
    \begin{equation}
      \bar{v} = \min_{u\in\mathcal U}\max_{w\in\mathcal W}\{\alpha_{uw} \bar{y}+\beta_{uw}\},
    \end{equation}
    where $\bar{v}$ and $\bar{y}$ denote the asymptotic average car-velocity and inter-vehicular distance.
  \item The $\nu$ cars move on an ``open'' road, where the velocity $v_1(t)$ of the first car (the leader one)
    varies over time but it converges a stationary one. In this case, $M^{uw}$ and $c^{uw}$ are given by
    $$M^{uw}=\begin{bmatrix}
        1 & 0 & \cdots & 0\\
        \alpha_{uw} & 1-\alpha_{uw} & & 0\\
        \vdots & \ddots & \ddots & \\
        0 & 0 & \alpha_{uw} & 1-\alpha_{uw}
      \end{bmatrix},$$
    and
    $$c^{uw}(t)={}^t[v_1(t),\;\; \beta_{uw},\; \cdots,\; \beta_{uw},\;\; \beta_{uw}].$$
    Again, the dynamics~(\ref{plcf1}) is stable under the condition $\alpha_{uw}\in[0,1]$, and the reverse
    behavior law at the stationary regime is obtained as follows.
    \begin{equation}
       \bar{y}=\max_{u\in\mathcal U}\min_{w\in\mathcal W_u}\frac{v_1-\beta_{uw}}{\alpha_{uw}},
    \end{equation}
    where $v_1$ denotes the asymptotic velocity of the first car.
\end{itemize}

\section{Anticipation modeling}

We present in this section our extension of the model~(\ref{Far1}) to multi-anticipative driving, where
the velocity of a car is calculated depending on the inter-vehicular
distance with respect to a given number $m$ of leader cars, with $m \geq 2$.

Our multi-anticipative model is an extension for the piecewise linear car-following model~(\ref{Far1}).
We use a \textit{minimum} form as in~(\ref{HoogBexel}) (rather than an additive form as in~(\ref{bexelius})).
Moreover, we use a \textit{uniform} form for the sensitivity with respect to the inter-vehicular distance as
in~(\ref{lenz}) (the inter-vehicular distance with respect to the $j$th leader is divided by $j$).
We consider the following model.
\begin{equation}\label{extplcf1}
  x_n(t+1)=x_n(t)+\min_{1\leq j\leq m}(1+\lambda)^{j-1}\min_{u\in\mathcal U}\max_{w\in\mathcal W} \left\{\alpha_{uw}
                     \left(\frac{x_{n-j}(t)-x_n(t)}{j}\right)+\beta_{uw}\right\},
\end{equation}
where $m$ is the number of leaders taken into account in anticipation, and $\lambda\geq 0$ is a discount parameter
with respect to the leader index.
The dynamics~(\ref{extplcf1}) can be written simply
\begin{equation}\label{extplcf2}
  x_n(t+1)=x_n(t)+\min_{1\leq j\leq m}\min_{u\in\mathcal U}\max_{w\in\mathcal W} \left\{\alpha_{uwj}
                     \left(\frac{x_{n-j}(t)-x_n(t)}{j}\right)+\beta_{uwj}\right\},
\end{equation}
where $\forall (u,w)\in\mathcal U\times\mathcal W, (\alpha_{uwj})_j, 1\leq j\leq m$ are increasing positive sequences.

The interpretation of the \textit{minimum} operator with respect to the $j$th leader in~(\ref{extplcf2})
is that a car $n$ maximizes its velocity under the constraints
\begin{equation}\label{const}
   x_n(t+1)-x_n(t) \leq \min_{u\in\mathcal U}\max_{w\in\mathcal W} \left\{\alpha_{uwj}
                     \left(\frac{x_{n-j}(t)-x_n(t)}{j}\right)+\beta_{uwj}\right\}, \quad 1\leq j\leq m.
\end{equation}
Discounting with respect to the leader indexes, made by introducing the multiplicative term $(1+\lambda)^{j-1}$,
permits to favor closer leaders over distant ones.
In the case where $m=1$, we retrieve the classical piecewise linear car-following model (without anticipation).
If $\lambda=0$, then the cars respond equally to the stimulus of all the leaders $j$, with $j=1, 2, \cdots, m$.

The meaning of the anticipation here is the following.
For example, the information that a car $i$, $i=n-1, n-2, \cdots, \max(1,n-m)$, decelerates
at time $t$, is immediately transmitted to the car $n$ that reacts at time $t+1$ (with anticipation), instead of $t+i$ (without anticipation).

We study below in section~\ref{sec-stab-anticip} the stability of the car dynamics~(\ref{extplcf1}) and
characterize the existence of stationary regimes. Two cases are distinguished: traffic on a \textit{ring} road,
and traffic on an \textit{open} road. In both cases, we give the asymptotic car positions when stationary regimes exist.
Transient traffic for the car dynamics~(\ref{extplcf1}) is treated in section~\ref{sec-transient}.    

\section{Stability analysis and stationary regime}
\label{sec-stab-anticip}

As in~\cite{Far12,Far08}, we consider $\nu$ cars moving on a 1-lane road without passing.
We first study the case where the cars move on a \textit{ring} road, and then explore the \textit{open} road case.

\subsection{Multi-anticipation on a ring road}

Cars being moving on a ring road, the indexes $n-j$, in the dynamics~(\ref{extplcf1}),
are cyclic in the set $\{1,2,\cdots, \nu\}$.
The idea here is to summarize the two minimum operators in~(\ref{extplcf1}) in only
one minimum operator, and then retrieve the form of the piecewise linear car-following dynamics without anticipation.
Let us denote by $\mathcal Z$ the set of all pairs of indexes $(j,u)$, with $1\leq j\leq m$ and $u\in\mathcal U$.
\begin{equation} \nonumber
  \mathcal Z = \{z=(j,u), 1\leq j\leq m, u\in\mathcal U\}.
\end{equation}
The dynamics~(\ref{extplcf1}) is then written
\begin{equation}\label{extplcf3}
  x_n(t+1)=\min_{z\in\mathcal Z}\max_{w\in\mathcal W} \{[M^{zw} x(t)]_n + c^{zw}_n\}, \quad 1\leq n\leq \nu,
\end{equation}
where the matrices $M^{zw}=M^{juw}$ and the column vectors $c^{zw}=c^{juw}$ are given as follows.
$$\begin{array}{l}
M^{juw}=\begin{bmatrix}
        1-\alpha_{juw}/j & 0 & \cdots & \alpha_{juw}/j & 0 & 0\\
        0 & 1-\alpha_{juw}/j & \ddots & \ddots & \alpha_{juw}/j & 0\\
        \vdots & \ddots & \ddots & \ddots & \ddots & \alpha_{juw}/j\\
        \alpha_{juw}/j & \ddots & \ddots & \ddots & \ddots & \vdots\\
        0 & \alpha_{juw}/j & \ddots & \ddots & \ddots & 0\\
        0 & 0 & \alpha_{juw}/j & \cdots & 0 & 1-\alpha_{juw}/j
      \end{bmatrix}, \\~~\\
c^{juw}={}^t[(\alpha_{juw}/j)(\nu\bar{y})+\beta_{juw},\;\; \beta_{juw},\;\; \beta_{juw} \cdots,\; \beta_{juw}].
\end{array}$$

The dynamics~(\ref{extplcf3}) have the same form as~(\ref{plcf1}). It is then interpreted as
a dynamic programming system associated to a stochastic game on a controlled Markov chain.
The stability is guaranteed
under the condition $\alpha_{juw}\in[0,1], \forall (j,u,w)\in\{1,2,\cdots,m\}\times \mathcal U\times\mathcal W$;
see~\cite{Far12,Far08} for more details.
The stationary regime is characterized by the additive eigenvalue problem
\begin{equation}\label{eigpbm1}
  \bar{v}+x_n=\min_{1\leq j\leq m}\min_{u\in\mathcal U}\min_{w\in\mathcal W}\{[M^{juw} x]_n+c^{juw}_n\},\quad 1\leq n\leq \nu\;,
\end{equation}
where $\bar{v}$ is the asymptotic average car-speed, the same for all cars; and where the vector $x$ gives the asymptotic car-positions,
up to an additive constant.
The following result gives a solution for the system~(\ref{eigpbm1}).
\begin{theorem}\label{stationary1}
  If $\forall (j,u,w)\in\{1,2,\cdots,m\}\times\mathcal U\times\mathcal W, \alpha_{juw}\in[0,1]$, then the
  system~(\ref{eigpbm1}) admits a solution $(\bar{v},x)$ given by:
  \begin{align}
    \bar{v} = & \min_{u\in\mathcal U}\min_{w\in\mathcal W}\{\alpha_{1uw} \bar{y}+\beta_{1uw}\}, \nonumber \\
    x   = & {}^t[(\nu-1)\bar{y} \quad (\nu-2)\bar{y} \quad \ldots \quad \bar{y} \quad 0]. \nonumber
  \end{align}
\end{theorem}
\proof Available in~\cite{FHL13}. \endproof

Theorem~\ref{stationary1} tells that the asymptotic car positions at the ring road are uniform.
This is due to the assumption that all the cars have the same behavior.

\subsection{Multi-anticipation on an open road}

We assume here that cars move on an open road, and that the speed $v_1(t)$ of the first car is given over time.
In addition, in order to analyze the stability and the stationary regime of the car dynamics, we assume
that the velocity of the first car approaches a constant value $v_1$; i.e. 
$$\lim_{t\to+\infty} v_1(t)=v_1.$$
Moreover, the number of anticipation cars for the $n$th car is $\min(n-1,m)$, instead of being always $m$.

The dynamics~(\ref{extplcf1}) is then written
\begin{equation}\label{extplcf4}
  x_n(t+1)=\min_{1\leq j\leq m}\min_{u\in\mathcal U}\max_{w\in\mathcal W} \{[M^{juw} x(t)]_n + c^{juw}_n\}, \quad 1\leq n\leq \nu,
\end{equation}
where the matrices $M^{juw}$ and the vectors $c^{juw}$ are given by
$$\begin{array}{ll}
    M^{juw} = \left(\begin{matrix}
        1 & 0 & 0 & \ldots & & & \ldots & 0\\
        0 & 1 & 0 & 0 & \ldots & & & 0\\
        \vdots & & \ddots & & & & & \vdots\\
        0 & 0 & \ldots & 1 & 0 & 0 & \ldots & 0\\
        \alpha_{juw}/j & 0 & \ldots & 0 & 1-\alpha_{juw}/j & 0 \ldots & 0\\
        0 & \alpha_{juw}/j & 0 & \ldots & 0 & 1-\alpha_{juw}/j & & \vdots\\
        \vdots & 0 & \ddots & 0 & \ldots & 0 \ddots & 0\\
        0 & \ldots & 0 & \alpha_{juw}/j & 0 & \ldots & 0 & 1-\alpha_{juw}/j
      \end{matrix}\right)
    \begin{matrix}
      1\\ 2\\ \vdots\\ j\\ j+1\\ j+2\\ \vdots\\ \nu
    \end{matrix}
\end{array}$$
and
$$\begin{array}{ll}
    c^{juw} = \left(\begin{matrix}
                      v_1(t)\\ +\infty\\ \vdots\\ +\infty\\ \beta_{juw}\\ \beta_{juw}\\ \vdots\\ \beta_{juw}
                     \end{matrix}\right)
    \begin{matrix}
      1\\ 2\\ \vdots\\ j\\ j+1\\ j+2\\ \vdots\\ \nu
    \end{matrix}
\end{array}$$
The entries $(M^{juw})_{ik}$ and $(c^{juw})_{i}$ for $i,k \leq j$ do not play any role in the car dynamics
since $(c^{juw})_i=+\infty, \forall i\leq j$. Those entries correspond to anticipation of a car $i$ with
respect to its $j$th leader that does not exist since $i\leq j$.

The stationary regime is characterized as follows.
\begin{equation}\label{stationopen}
  \bar{v}+x_n=\min_{1\leq j\leq m}\min_{u\in\mathcal U}\max_{w\in\mathcal W}\{[M^{juw} x]_n+c^{juw}_n\},\quad 1\leq n\leq \nu\;.
\end{equation}
The following result gives a solution for the system~(\ref{stationopen}).
\begin{theorem}\label{opentheo}
  For all $y\in\mathbb R$ satisfying $\min_{u\in\mathcal U}\max_{w\in\mathcal W}(\alpha_{1uw}y+\beta_{1uw})=v_1$, the couple $(\bar{v},x)$
  is a solution for the system~(\ref{stationopen}), where $\bar{v}=v_1$ and $x$ is given up to an
  additive constant by
  \begin{equation}\label{eqx}
    x={}^t[(\nu-1)y, \quad (\nu-2)y, \quad \cdots, \quad y, \quad 0].
  \end{equation}
\end{theorem}
\proof Available in~\cite{FHL13}. \endproof

Theorem~\ref{opentheo} tells that the configuration of the car-positions at the stationary regime is uniform.
Moreover, it gives the asymptotic inter-vehicular distance as a function of the asymptotic speed $v_1$ of car~1.
An important remark here is that the car-velocity obtained here (in the multi-anticipation case) is the same as
the one obtained for the piecewise linear car-following model without anticipation~(\ref{Far1}).
Moreover, the \textit{optimal strategy} of driving at the stationary regime in the case of multi-anticipation
is to drive by taking into account only one leader: $(\bar{j},\bar{u},\bar{w})=(1,\bar{u},\bar{w})$.
Therefore, at the stationary regime, once the traffic is stabilized, it is not necessary for drivers to
take into consideration more than one leader.

Another important remark here is that even though the cars reduce
their approach in the multi-anticipative dynamics (due to the \textit{minimum} operator over the leader indices),
comparing to their dynamics without anticipation; as long as the cars approach the stationary
regime, where the traffic is stable, they retrieve what they have lost at the transient regime.
Therefore, by introducing the \textit{minimum} operator on the multi-anticipative dynamics, the traffic becomes smoother,
without decreasing the stationary car speed.

\section{Transient traffic}
\label{sec-transient}

Let us consider the same example taken in~\cite{Far12}, and adapt it to multi-anticipative driving.
We simulate the car-dynamics~(\ref{extplcf1}), with a time unit of 1/2~second (s), and a distance unit of 1~meter (m).
We keep the same parameters of the model as those of Example~1 of~\cite{Far12} 
(the parameters have been determined by approximating a given behavior law).

The car dynamics are simulated on a one lane road of about 10,000 meter. We vary in time the number of leaders and the
velocity of the first car, in order to show the effect of multi-anticipation on the transient traffic.
We give the results on Table~\ref{tab1}. 

The trajectories of the leader car are the same for all views of Table~\ref{tab1}.
We see that the car trajectories shown in Table~\ref{tab1} are smoothed in case of multi-anticipation driving.
Although in practice, the number of leaders taken into account by drivers is in general less than 5,
we simulated here the anticipation with up to 100 leaders; 
which could be possible if we consider communicating and/or automated cars.
   
\begin{table}
  \caption{Traffic on a 1-lane road. On the x-axis: time. On the y-axis: car-position. The number of cars 
               taken into account in anticipation are 1, 5, 10, 20, 50 and 100.
               The length of the road is 10,000 meter. The total simulation time is 500 seconds.}
  \begin{tabular}{cc}
     & \\  
    1 leader anticipation & 5 leaders anticipation \\
    \includegraphics[width=7.4cm]{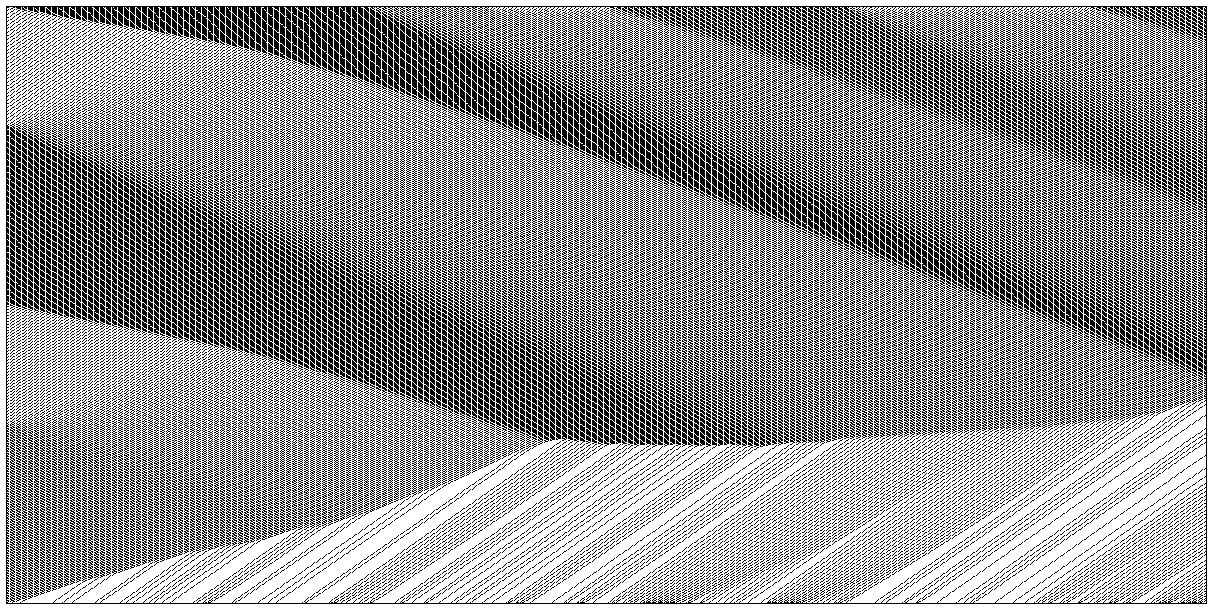} & \hspace{0cm} \includegraphics[width=7.4cm]{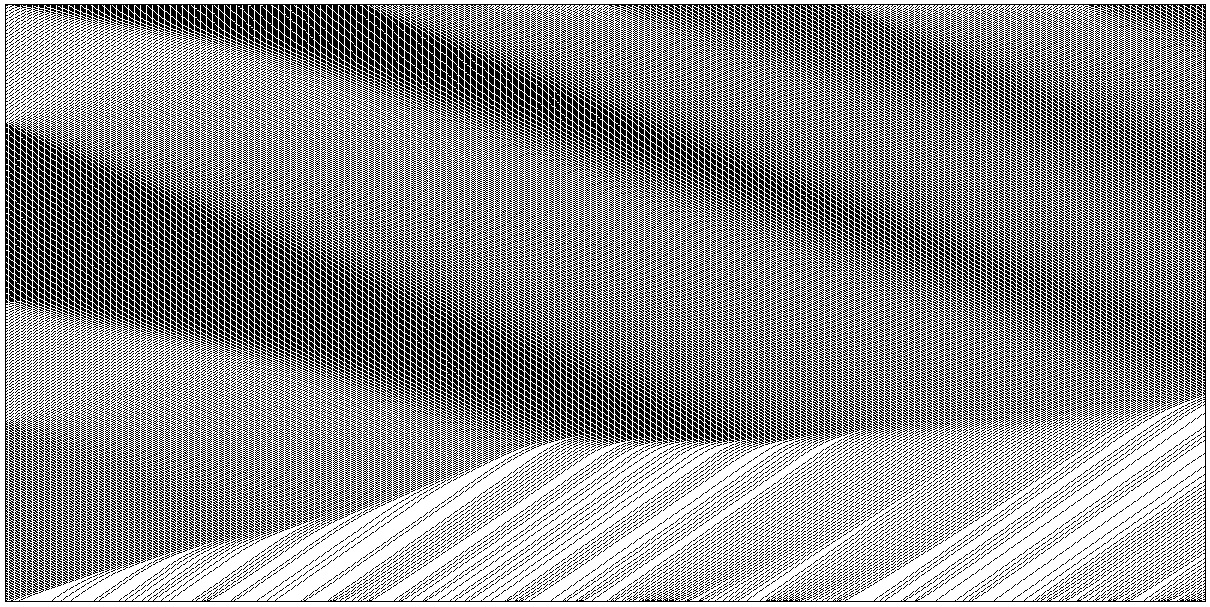} \\~~\\
    10 leaders anticipation & 20 leaders anticipation \\
    \includegraphics[width=7.4cm]{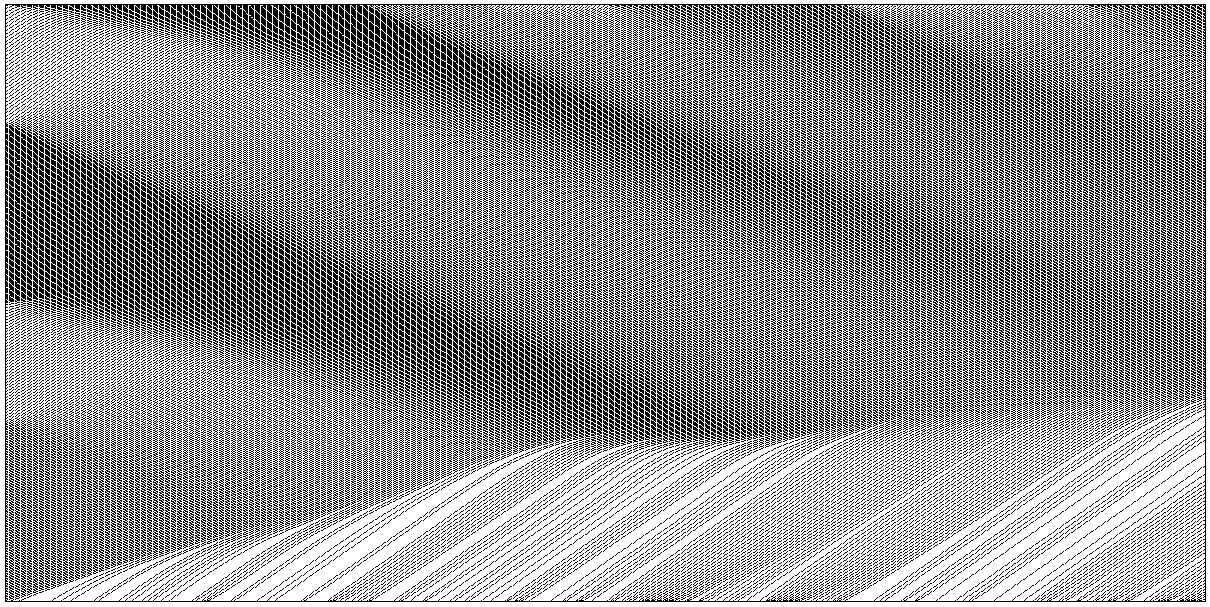} & \hspace{0cm} \includegraphics[width=7.4cm]{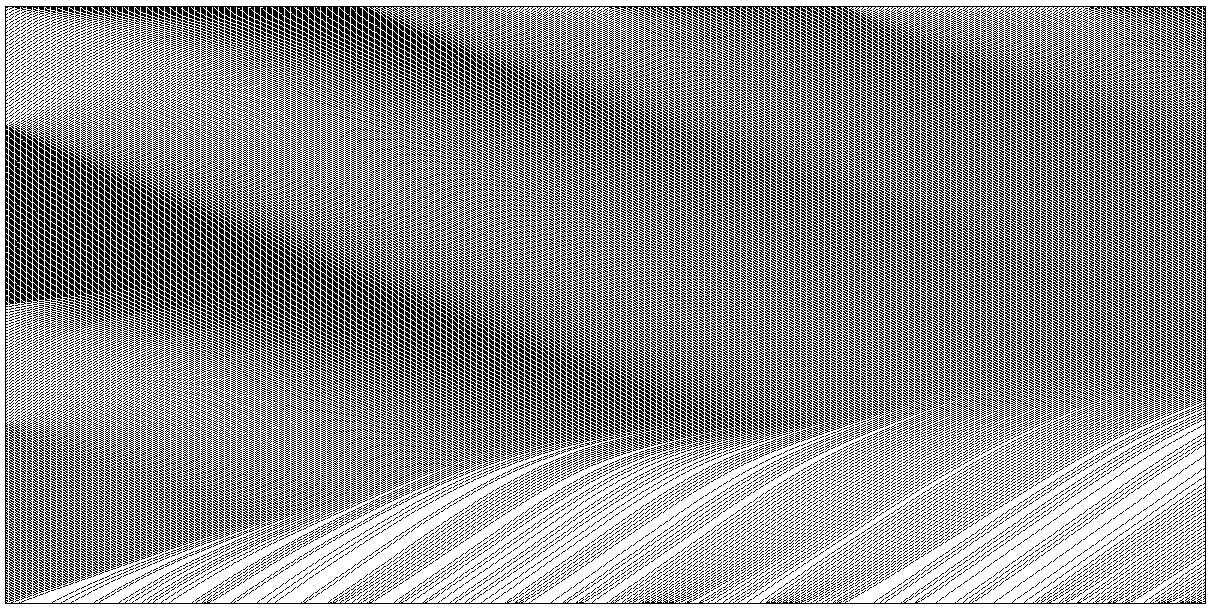} \\~~\\
    50 leaders anticipation & 100 leaders anticipation \\
    \includegraphics[width=7.4cm]{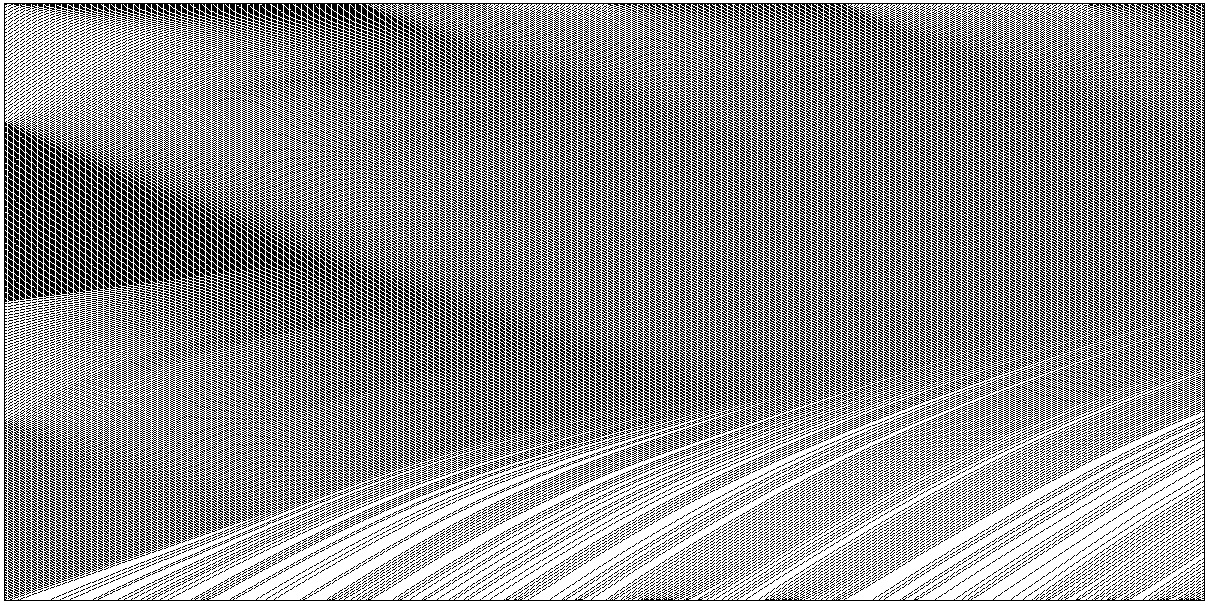} & \hspace{0cm} \includegraphics[width=7.4cm]{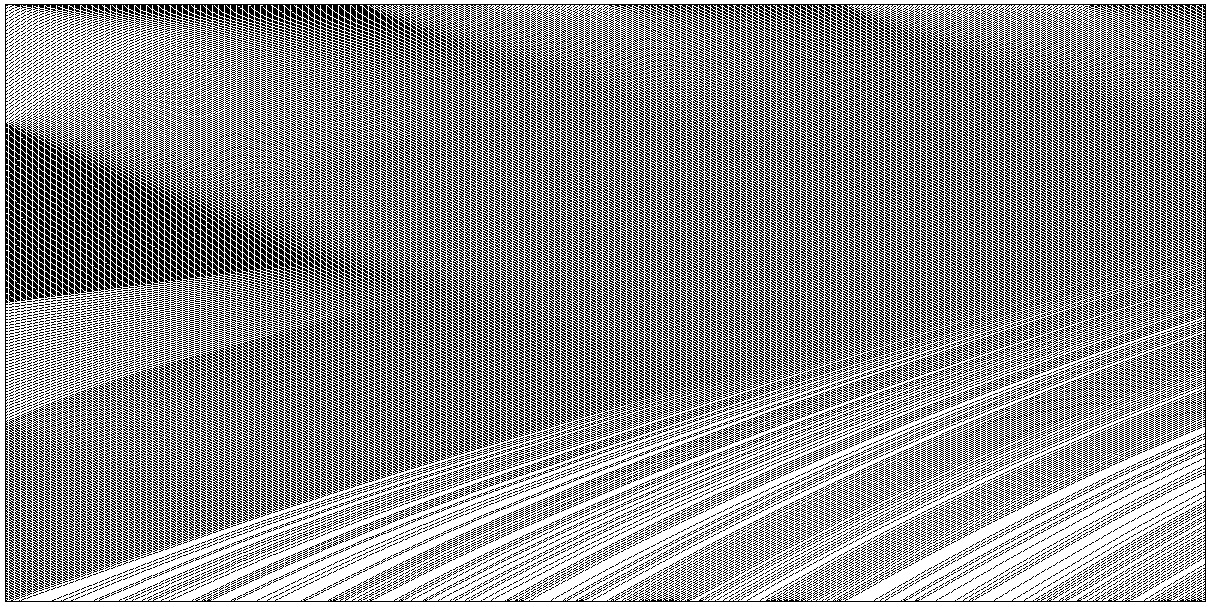}
  \end{tabular}  
  \label{tab1}
\end{table}

\newpage
This page is intentionally left blank

\chapter{Demand-dependent train dynamics control}
\label{chap-stochtrain}

The main references for this chapter are~\cite{FNHL18, FNHL17b, Far18}.
We present some extensions of the max-plus algebra model of the train dynamics, presented in chapter~\ref{chap-maxtrain}.
The extensions aim to make the train dwell times at platforms dependent on the passenger travel demand, and, in some models,
to control this dependence.
We first show in section~\ref{sec-tip} that if the dependence of the train dwell times on the passenger travel demand is
not controlled, then, the train dynamics are unstable.
In a second model (section~\ref{sec-stable}), we present a control model of the train dwell times based on the feedback of
the train dynamics as well as the passenger travel arrival rates.
We show analytically that the train dynamics are stable with the proposed control model, and derive by means of numerical
simulation, the asymptotic regimes of the train dynamics with the asymptotic average train time headways and frequencies,
as functions of the number of running trains and of the passenger travel demand.
Finally, we discuss some further extensions of the proposed modeling and control approach.

\section{Natural instability of the train dynamics}
\label{sec-tip}

We consider in this section, train dynamics that take into account the passenger travel demand.
We introduce this dependence by assuming that train dwell time $w_j$ at platform $j$ depends
on the passenger volume at the same platform $j$, which depends on the safe separation time
$g_j$ on the same platform (the time during which passengers are accumulated at the platform),
or equivalently on the train time headway $h_j$.

The model we present here extends naturally the model of chapter~\ref{chap-maxtrain},
without optimizing the effect of the passenger volumes on the train dynamics.
In other words, we assume that the dependence of the dwell times at platforms with the passenger
demand is uncontrolled. We will see that, in this case, the train dynamics are unstable.
We notice that we do not model the stocks of passengers at platforms or inside the trains.
We assume that the train dwell times at platforms depend directly on the average passenger arrival rates.
Also, we model here the effect of only embarking time of passengers into the trains, on the train dwell time.
We do not model the effect of the alighting time of passengers from the trains, on the train dwell time.

In order to model the effect of passengers on the train dwell times at platforms, we consider the following
additional constraint on the dwell time at platforms (in addition to constraints~(\ref{const1}) 
and~(\ref{const2}) of chapter~\ref{chap-maxtrain}).
\begin{equation}\label{eq-w1}
  w_j^k \geq \begin{cases}
                \frac{\lambda_j}{\alpha_j}\; g_j^k, & \text{ if } j \text{ indexes a platform}, \\
                0 & \text{ otherwise}.
             \end{cases}
\end{equation}
where
\begin{itemize}
    \item $\alpha_j$ is the total passenger upload rate from platform $j$ onto the trains, if $j$ indexes a platform;
           and $\alpha_j$ is zero otherwise.
    \item $\lambda_j := \sum_i \lambda_{ji}$ is the average rate of the total arrival flow of passengers to platform $j$, if $j$ indexes a platform;
	  and $\lambda_j$ is zero otherwise; where $\lambda_{ji}$ is the origin-destination arrival rates of passengers to platform $j$, with platform $i$ as destination.
\end{itemize}
$\lambda_j g_j^k$ gives the total arrival of passengers in the time interval $g_j^k$ during which passengers are accumulated onto platform $j$.
Therefore, $\lambda_j g_j^k / \alpha_j$ gives the minimum necessary train dwell time to upload the total arrival $\lambda_j g_j^k$.

By taking into account the additional constraint~(\ref{eq-w1}), the constraint~(\ref{const1}) is modified as follows.
\begin{equation}\label{new-constr}
    d_j^k \geq d_{j-1}^{k-b_{j}} + r_j +  \max \left\{ \underline{w}_j, \frac{\lambda_j}{\alpha_j} g_j^k \right\}.
\end{equation} 

We now modify the dynamic model~(\ref{eq-d1}) by taking into account the new constraint~(\ref{new-constr}).
We obtain, for nodes $j$ indexing platforms, the following.

\begin{equation}\label{eq-dd}
  d_j^k = \max \left\{
  \begin{array}{l}
     d_{j-1}^{k-b_{j}} + r_{j} + \underline{w}_j, \\ ~~\\
     \left(1+\frac{\lambda_j}{\alpha_j}\right) d_{j-1}^{k-b_{j}} - \left(\frac{\lambda_j}{\alpha_j}\right) d_j^{k-1} + \left(1+\frac{\lambda_j}{\alpha_j}\right) r_{j}, \\~~\\
     d_{j+1}^{k-\bar{b}_{j+1}} + \underline{s}_{j+1}.
  \end{array} \right.
\end{equation}
For non-platform nodes, the dynamics remain as in~(\ref{eq-d1}).

Let us notice that the dynamic system~(\ref{eq-dd}) has explicit and implicit terms.
Moreover, it can be written in the form~(\ref{eq-dp3}), i.e. as follows.
\begin{equation}\label{matrix_form}
  d_j^k = \max_{u\in\mathcal U} [(M^u d^{k-1})_j + (N^u d^{k})_j + c^u_j],
\end{equation}
where $M^u$ and $N^u$ are square matrices, and $c^u$ is a family of vectors, for $u\in \mathcal U$.
The matrices $N^u, u\in \mathcal U$ express implicit terms.
We notice that it suffices that $\exists j, \lambda_j > 0$ to have one of the matrices $M^u, u\in\mathcal U$ or $N^u, u\in\mathcal U$ not being 
sub-stochastic~\footnote{A matrix $A$ is sub-stochastic if $0 \leq A_{ij}\leq 1, \forall i,j$ and if $\sum_i A_{ij} \leq 1, \forall i$.}, 
since we have in this case $\alpha_j/(\alpha_j-\lambda_j) > 1$ and $-\lambda_j/(\alpha_j-\lambda_j) < 0$; see the dynamics~(\ref{eq-dd}).
Therefore, the dynamic system~(\ref{eq-dd}) cannot be seen as a dynamic programming system of a stochastic optimal control problem
of a Markov chain; see section~\ref{sec-dps} of chapter~\ref{chap-stoch}.

It is easy to see that if $m=0$ or $m=n$, then the dynamic system~(\ref{eq-dd}) is fully implicit (it is not triangular), and therefore, we can say that it admits an
asymptotic regime with a unique additive eigenvalue $h = +\infty$, which is also the average growth rate of the system, and which is the asymptotic average 
train time-headway.
This case corresponds to $0$ or $n$ running trains on the metro line. No train departure is possible for these two cases. We have the
average train frequency $f=0$ corresponding to the average time headway $h=+\infty$. 
We can also show that if $0 < m < n$, then the dynamic system~(\ref{eq-dd}) is triangular. That is to say that it is not fully implicit, and
there exists an order of updating the components of the state vector $d^k$, in such a way that no implicit term appears.

We know that the dynamic system~(\ref{eq-dd}) is unstable in general.
Indeed, as we noticed above, the matrices $M^u, u\in\mathcal U$ or $N^u, u\in\mathcal U$ are not sub-stochastic.
The consequence is that the dynamics are not non-expansive. In fact, they are additive 1-homogeneous but not monotone.
Many behaviors are possible for the dynamic system, depending on the parameters and on the initial
conditions (expansive behavior, chaotic behavior, etc.)
In practice, let us assume the $k^{\text{th}}$ arrival of a given train to platform $j$ is delayed.   
Then $g^k_j$ will increase by definition, and $w^k_j$ will also increase by application of~(\ref{eq-w1}).
By consequent of the increasing of $w^k_j$, the departure $d^k_j$ will be delayed, and the delay
of $d^k_j$ would be longer than the one of $a^k_j$, because, it accumulates the delay of $a^k_j$ 
and the increasing of $w^k_j$. Consequently, the arrival of the same train to platform $j+1$
(downstream of platform $j$) will be delayed longer comparing to its arrival to platform~$j$.
Therefore, the application of the control law~(\ref{eq-w1}) amplifies the train delays and propagates
them through the metro line. 
We notice that the instability of this kind of dynamics has already been pointed out; see for example~\cite{BCB91}.

\section{Stable dynamic programming model}
\label{sec-stable}

We propose here an adaptation of the train dynamics~(\ref{eq-dd}) in order to guarantee its stability.
As shown above, system~(\ref{eq-dd}) is unstable because of the relationship~(\ref{eq-w1}).
In order to deal with this instability, we propose to replace the train dwell time formula~(\ref{eq-w1}) by the following.
\begin{equation}\label{eq-w4}
  w_j^k \geq \begin{cases}
                \overline{w}_j - \frac{\theta_j^k}{\lambda_j^k/\alpha_j^k} g_j^k & \text{ if } j \text{ indexes a platform}, \\
                0 & \text{ otherwise}.
             \end{cases} 
\end{equation}
where the sign of the relationship between the dwell time $w_j^k$ and the safe separation time $g_j^k$ is inversed,
without inversing the relationship between the dwell time $w_j^k$ and the ratio $\lambda_j^k/\alpha_j^k$; and where
$\overline{w}_j$ (maximum dwell time at node $j$) and $\theta_j^k$ are control parameters.

The dynamics~(\ref{eq-dd}) are now rewritten, for nodes $j$ indexing platforms, as follows.
\begin{equation}\label{eq-dd2}
  d_j^k = \max \left\{
  \begin{array}{l}
     d_{j-1}^{k-b_{j}} + r_{j} + \underline{w}_j, \\~~ \\
     \left(1 - \delta_j^k\right) d_{j-1}^{k-b_{j}} + \delta_j^k d_j^{k-1} + \left(1 - \delta_j^k\right) r_{j} + \overline{w}_j, \\ ~~\\
     d_{j+1}^{k-\bar{b}_{j+1}} + \underline{s}_{j+1},
  \end{array} \right.
\end{equation}
where $\delta_j^k = \theta_j^k \alpha_j^k/\lambda_j^k, \forall j, k$.
For non-platform nodes, the dynamics remain as in~(\ref{eq-d1}).
In case where $\delta_j^k$ are independent of $k$ for every $j$,
the dynamic system (\ref{eq-dd2}) can be written under the form~(\ref{matrix_form}).

As for the dynamic system~(\ref{eq-dd}), in case where $m=0$ or $m=n$, the system~(\ref{eq-dd2}) is fully implicit, and admits an
asymptotic regime with a unique average asymptotic train time-headway $h = +\infty$ (no train movement is possible).
In case where $0 < m < n$, the dynamic system~(\ref{eq-dd2}) is triangular.
In the latter case, and if $\delta_j^k$ are independent of $k$ for every $j$, then the dynamic
system~(\ref{eq-dd2}) can be written under the form~(\ref{eq-dp3}).
If, in addition, $0 \leq \delta_j \leq 1, \forall j$, then $M^u$ and $N^u$ are sub-stochastic matrices,
i.e. satisfying $M^u_{ij} \geq 0, N^u_{ij} \geq 0$, $\sum_j (M^u_{ij}) \leq 1$ and $\sum_j (N^u_{ij}) \leq 1$.
Moreover, we have $\sum_j (M^u_{ij}+N^u_{ij}) = 1, \forall i,u$.
In this case,~(\ref{eq-dd2}) is a dynamic programming system of an optimal control problem of a Markov chain,
whose transition matrices and reward vectors
can be calculated from $M^u, N^u$ and $c^u, u\in\mathcal U$ (they are the transition matrices and reward vectors
corresponding to the equivalent explicit dynamic system obtained by solving the implicit terms of~(\ref{eq-dd2})).

\begin{theorem}\label{stable}
  If $0 < m < n$ and if $\delta_j^k$ are independent of $k$ for every $j$, and $0 \leq \delta_j \leq 1, \forall j$,
  then the dynamic system~(\ref{eq-dd2}) admits a stationary regime with a unique average growth rate $h$
  (interpreted here as the asymptotic average train time-headway), independent of the initial state $d^0$.
\end{theorem}

\begin{figure}[htpb]
      \centering
	  \includegraphics[scale = 0.42]{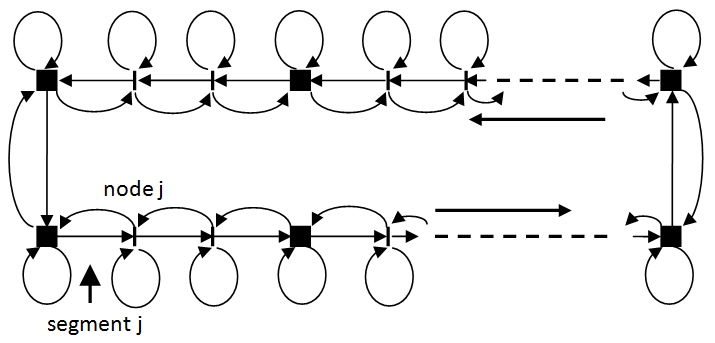}
      \caption{The graph associated to $\mathbf{f}$ in the proof of Theorem~\ref{stable}, which is also the same
         graph associated to the Max-plus linear system in the proof of Theorem~\ref{th-dp-mp}.}
      \label{graph2}
\end{figure}

\proof The proof uses Theorem~\ref{th-dps2}. 
  Let us denote by $\mathbf{f}$ the map associated to the dynamic system~(\ref{eq-dd2}), and use the notations
  $\mathbf{f}_{(i)}, i\in{0,1}$ and $\tilde{\mathbf{f}}$ as defined in section~\ref{sec-gdps} of chapter~\ref{chap-stoch}.
  Since $M^u$ and $N^u$ are sub-stochastic matrices, with $\sum_j (M^u_{ij}+N^u_{ij}) = 1, \forall i,u$, then
  $\mathbf{f}$ is additive 1-homogeneous and monotone.
  The graph associated to $\mathbf{f}$ is strongly connected; see Figure~\ref{graph2}.
  Indeed, this graph includes the one of the Max-plus linear dynamics (graph of Figure~\ref{graph1}), which is already 
  strongly connected. The graph associated to $\mathbf{f}_{(0)}$ is acyclic since $0<m<n$.
\endproof  

We do not yet have an analytic formula for the asymptotic train time-headway (we know that it coincides with the average growth rate $h$),
but Theorem~\ref{stable} guarantees its existence and its uniqueness.
Therefore, by iterating the dynamics~(\ref{eq-dd2}), one can approximate, for any fixed train density $\rho$, the associated asymptotic average
train time-headway $h(\rho)$ as follows.
\begin{equation}\label{hsim}
  h(\rho) \approx d^K_j/K, \forall j, \text{ for a large } K.
\end{equation}

In the following, we give a result (Theorem~\ref{th-dp-mp}) which tells us under which condition on the control parameters
$\bar{w}_j$ and $\delta_j$, the dynamic system~(\ref{eq-dd2}) is Max-plus linear.
We will use this result in section~\ref{sec-param}, in order to derive an approach for fixing the control parameters
$\bar{w}_j$ and $\delta_j$ in such a way that the effect of passenger arrivals on the train dynamics will be well modeled.

\begin{theorem}\label{th-dp-mp}
  Let $\tilde{h}$ be the asymptotic average growth rate of the Max-plus linear system~(\ref{eq-d1}).
  The dynamic programming system~(\ref{eq-dd2}) with parameters $\bar{w}_j = \tilde{h}, \forall j$, and $\delta_j = 1, \forall j$,
  is a Max-plus linear system, whose asymptotic average growth rate coincides with $\tilde{h}$.
\end{theorem}

\proof It is easy to see that if $\delta_j = 1, \forall j$, then system~(\ref{eq-dd2}) is a Max-plus linear
one whose associated graph has $n$ additional cycles (which are loop-cycles) comparing to the one associated
to system~(\ref{eq-d1}); see Figure~\ref{graph2}. Moreover, if $\bar{w}_j = \tilde{h}, \forall j$, then
the cycle mean of the loops are all equal to $\tilde{h}$.
All the other parameters corresponding to the characteristics of the metro line and of the
trains running on it, remain the same as the ones of system~(\ref{eq-d1}).
Therefore, as in Theorem~\ref{th-mpm}, the asymptotic average  time-headway $h$ is given by the maximum cycle mean
of the graph associated to the Max-plus linear system obtained from system~(\ref{eq-dd2}).
Four different elementary cycles are distinguished on that graph; see Figure~\ref{graph2}.
\begin{itemize}
  \item The Hamiltonian cycle in the direction of the train movements, with mean $\sum_j \underline{t}_j/m$.
  \item All the cycles of two links relying nodes $j-1$ and $j$, with mean $\underline{t}_j+\underline{s}_j$ each.
  \item The Hamiltonian cycle in reverse direction of the train dynamics, with mean
    $\sum_j \underline{s}_j/(n-m)$.
  \item The $n$ loop-cycles with mean $\tilde{h}$.
\end{itemize}
Hence
$$h = \max \left\{ \frac{\sum_j \underline{t}_j}{m}, \max_j (\underline{t}_j+\underline{s}_j),
   \frac{\sum_j \underline{s}_j}{n-m}, \tilde{h}\right\} = \tilde{h}.$$
\endproof

\section{How to fix the control parameters $\bar{w}_j$ and $\delta_j$}
\label{sec-param}

According to Theorem~\ref{th-dp-mp}, if we fix $(\bar{w}_j, \delta^k_j) = (\tilde{h}(\rho),1)$ in~(\ref{eq-dd2}), 
then we obtain a Max-plus linear train dynamics, which do not take into account the passenger demand.

We assume here that $\alpha_j^k, \lambda_j^k$ and $\delta_j^k$ are independent of $k$,
and then denote $\alpha_j, \lambda_j$ and $\delta_j$ respectively.
Let us consider the metro line as a server of passengers. We assume that 
the average arrival of passengers to a platform $j$ is $\lambda_j$, and that the average passenger service rate
at a platform $j$ is $\min(\alpha_j, \kappa /h)$, where $\kappa$ denotes the train capacity (maximum number of passengers in a train).
Under the Max-plus linear model~(\ref{eq-d1}) of the train dynamics, the asymptotic average service rate depends on the number $m$ of trains,
or equivalently on the train density $\rho$, since the asymptotic average train time headway is given as a function of $\rho$; see~(\ref{diag2}).
Therefore, the asymptotic average service rate is given by $\min(\alpha_j, \kappa /\tilde{h}(\rho))$.
Let us use the following notation.
\begin{equation}\label{lamtild}
  \tilde{\lambda}_j(\rho) := \min(\alpha_j, \kappa /\tilde{h}(\rho)).
\end{equation}
Therefore, the metro line as a server operating under the Max-plus linear model is stable under the condition 
$\lambda_j < \tilde{\lambda}_j(\rho), \forall j$.

The stability condition of the metro line seen as a server means that if $\lambda_j \leq \tilde{\lambda}_j(\rho), \forall j$, then the train dynamics is
Max-plus linear and can then serve passengers without adapting the train dwell times to the passenger arrival rates.
Basing on this remark, we propose an approach for fixing the parameters $\bar{w}_j$ and $\delta_j$ in function of the passenger arrivals, in such a way that
\begin{itemize}
  \item In case where $\lambda_j \leq \tilde{\lambda}_j(\rho), \forall j$, the dynamic system behaves as a Max-plus linear one.
    Therefore, the train dwell times are not constrained by the arrival rates of passengers.
  \item In case where $\exists j, \lambda_j > \tilde{\lambda}_j(\rho)$,  the system switches to a dynamic programming one.
    The train dwell times are constrained by the arrival rates of passengers.
\end{itemize} 

Let us fix $\bar{w}_j$ and $\delta_j$ as follows.
\begin{align}
  & \bar{w}_j(\rho) := \tilde{h}(\rho), \forall \rho, j. \label{param1} \\
  & \delta_j(\rho) := \frac{\tilde{\lambda}_j(\rho)}{\max\left(\lambda_j,\tilde{\lambda}_j(\rho)\right)}, \forall \rho, j. \label{param3}  
\end{align}
We notice here that fixing $\delta_j(\rho)$ as in~(\ref{param3}) is equivalent to fixing $\theta_j(\rho)$ as follows.
\begin{equation}
   \theta_j(\rho) := \frac{\tilde{\lambda}_j(\rho)}{\max\left(\lambda_j,\tilde{\lambda}_j(\rho)\right)} \; \frac{\lambda_j}{(\alpha_j-\lambda_j)}, \forall \rho, j. \label{param2}
\end{equation}

We then have $0\leq \delta_j(\rho) \leq 1$ by definition, and
\begin{itemize}
  \item In case where $\lambda_j \leq \tilde{\lambda}_j(\rho), \forall j$, we have $\delta_j(\rho) = 1, \forall j$, and the dynamic system is Max-plus linear,
     where the dwell times are not constrained by the arrival rates of passengers.
  \item In case where $\exists j, \lambda_j > \tilde{\lambda}_j(\rho)$, we have $\exists j, \delta_j(\rho) < 1$, and the system switches to a dynamic programming one. 
     The dwell times are constrained by the arrival rates of passengers.     
\end{itemize}

We then have the following result.
\begin{theorem}\label{th_stab}
  For any fixed value of the train density $\rho$ on the metro line,
  the dynamic system~(\ref{eq-dd2}) with parameters $\bar{w}_j$ and $\delta_j$ fixed dependent on $\rho$ as in~(\ref{param1}) and~(\ref{param3})
  respectively, 
  admits a stationary regime with a unique  asymptotic average growth rate $h(\rho)$, independent of the initial state $d^0$. 
  Moreover, we have $h(\rho) \geq \tilde{h}(\rho)$.
\end{theorem}
\proof The proof follows from Theorem~\ref{stable} and from all the arguments given above in this section,
in particular from $0\leq \delta_j(\rho)\leq 1, \forall \rho, j$. \endproof

\section{Derivation of the asymptotic average dwell and safe separation times}

From~(\ref{form2}) we have, $w = t - r$. From~(\ref{form3}) we have $t = (m/n) h$. Then 
we obtain $w = (m/n) h - r$. Therefore, if we have $h(\rho)$ by simulation (by means of~(\ref{hsim})), we deduce $w(\rho)$ as follows.
\begin{equation} \label{formw}
   w(\rho) =  (\rho / \bar{\rho})\; h(\rho) - r.
\end{equation}

Similarly, from~(\ref{form1}) we have $g = r + s$. From~(\ref{form3}) we have $s = ((n-m)/n) h$. Then 
we obtain $g = r + ((n-m)/n) h$. Therefore, from~(\ref{hsim})), we obtain
\begin{equation} \label{formg}
   g(\rho) = r + \left(1 - \rho / \bar{\rho} \right) h(\rho).
\end{equation}

In~\cite{Far18}, other results give the minimum number of trains to run on the metro line, in order that the train dynamics
behaves as Max-plus linear one, and then the passenger arrivals have no effect on the train dwell times and on the train dynamics.
Practical implementations of this approach have also been discussed in~\cite{Far18}.

\section{Numerical results}
\label{sec-sim}

We present here some numerical results.
We consider the symmetric case for the passenger arrivals, where the latter are the same at all platforms.
We vary the average rate in order to derive its effect on the train dynamics and on the traffic phase diagrams.
Table~\ref{tab_1}, Table~\ref{tab_2} and Table~\ref{tab_3} show the obtained results.
In the left side (respectively right side) of Table~\ref{tab_1} we see the increase of the train time-headways 
(respectively degradation of the train frequencies) due to increases 
of the passenger arrival rates.
In Table~\ref{tab_2}, we see the increase of the average dwell times $w$ at all the nodes, and the one at platforms, 
due to increases of the passenger arrival rates.
In Table~\ref{tab_3}, we see the degradation of the safe separation times $g$, due to increases of the passenger arrival rates.

\begin{table}[h]
\centering
\caption{Parameters of the line considered.}~~\\~~ 
\begin{tabular}{|l||l|}
  \hline
  Number of stations & 9 ($\Rightarrow $ 18 platforms)\\
  \hline
  Segment length & about 200 meters (m.) \\
  \hline
  Free train speed $v_{\text{run}}$ & 22 m/s (about 80 km/h) \\
  \hline 
  Train speed on terminus & 11 m/s (about 40 km/h) \\
  \hline
  Min. dwell time $\underline{w}$ & 20 seconds \\
  \hline
  Min. safety time $\underline{s}$ & 30 seconds \\
  \hline  
  Metro line length & 17.294 km \\
  \hline                          
  Passenger train capacity $\kappa$ & 500 passenger/train\\
  \hline  
  Passenger upload rate $\alpha$ & 30 passengers/s \\
  \hline  
\end{tabular}
\label{tab-param}
\end{table}

\begin{table}[thbp]
\centering
\caption{Asymptotic average train time-headway $h$ and frequency $f$ as functions of the number $m$ of running trains.
         The average passenger arrivals to platforms are given in the figure (denoted by $c$).} ~~\\~~
\begin{tabular}{|c|c|}
  \hline
   & \\
  \includegraphics[scale=0.32]{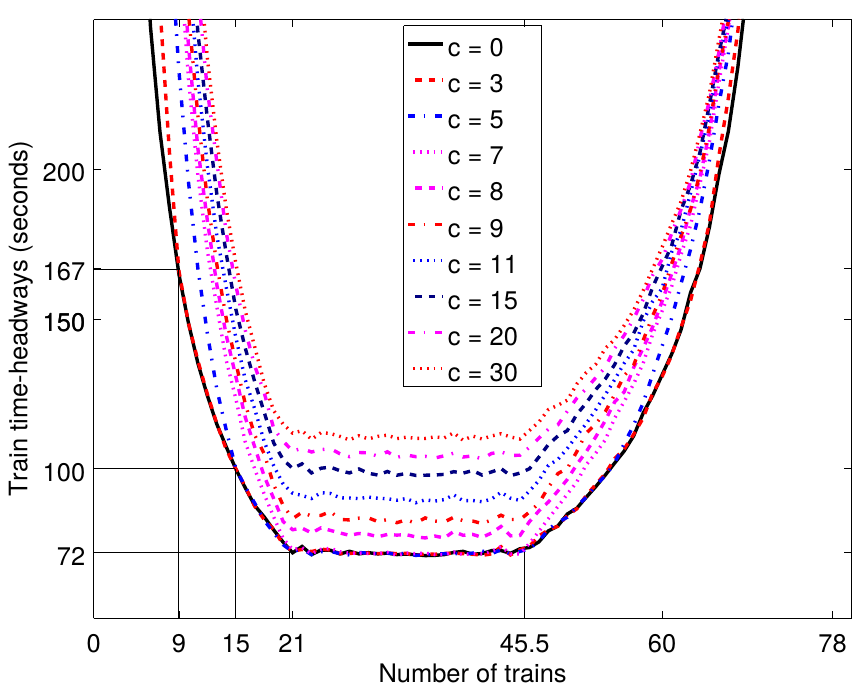} & \includegraphics[scale=0.32]{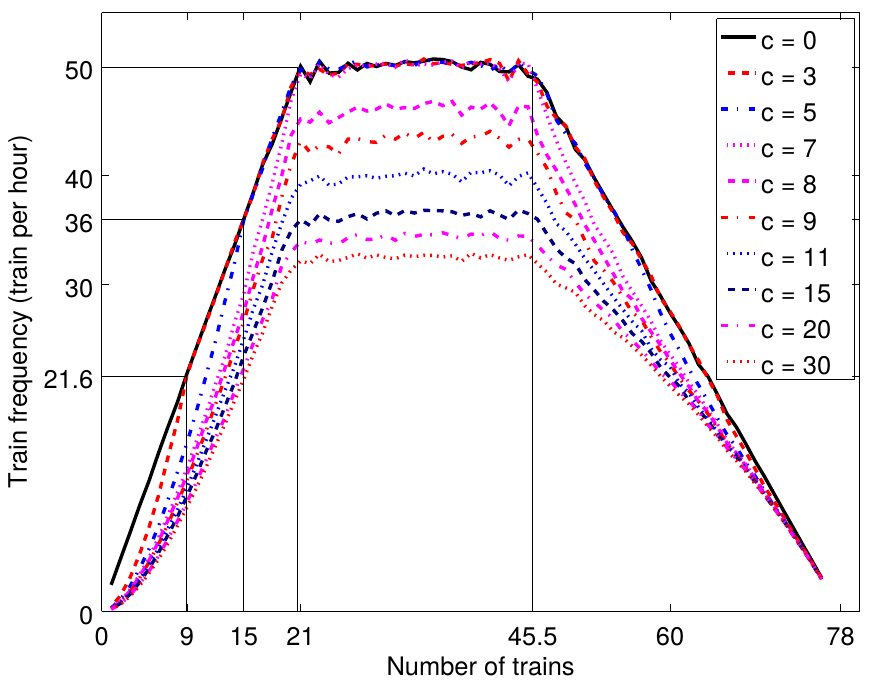} \\
  \hline
\end{tabular}
\label{tab_1}
\end{table}


\begin{table}[thbp]
\centering
\caption{Asymptotic average train dwell time $w$ at all nodes and dwell time at platforms; as functions of the number $m$ of running trains.
         The average passenger arrivals to platforms are given in the figure (denoted by $c$).}~~\\~~
\begin{tabular}{|c|c|}
  \hline
   & \\
  \includegraphics[scale=0.32]{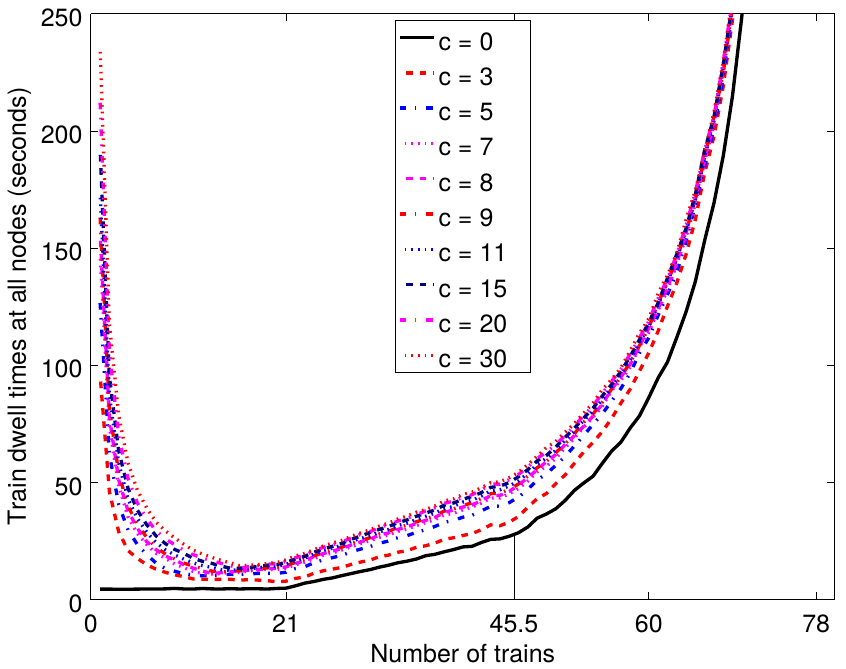} & \includegraphics[scale=0.32]{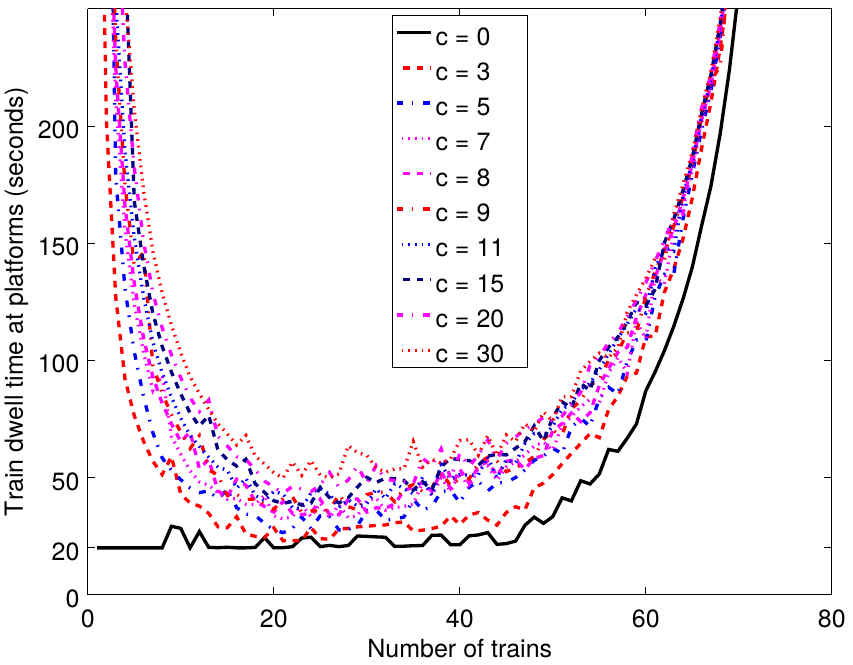} \\
  \hline
\end{tabular}
\label{tab_2}
\end{table}

\begin{table}[thbp]
  \centering  
  \caption{Asymptotic average train safe separation time $g$ as a function of the number $m$ of running trains.
         The average passenger arrivals to platforms are given in the figure (denoted by $c$).} ~~\\~~
  \begin{tabular}{|c|}
  \hline
   \\       
  \includegraphics[scale=0.32]{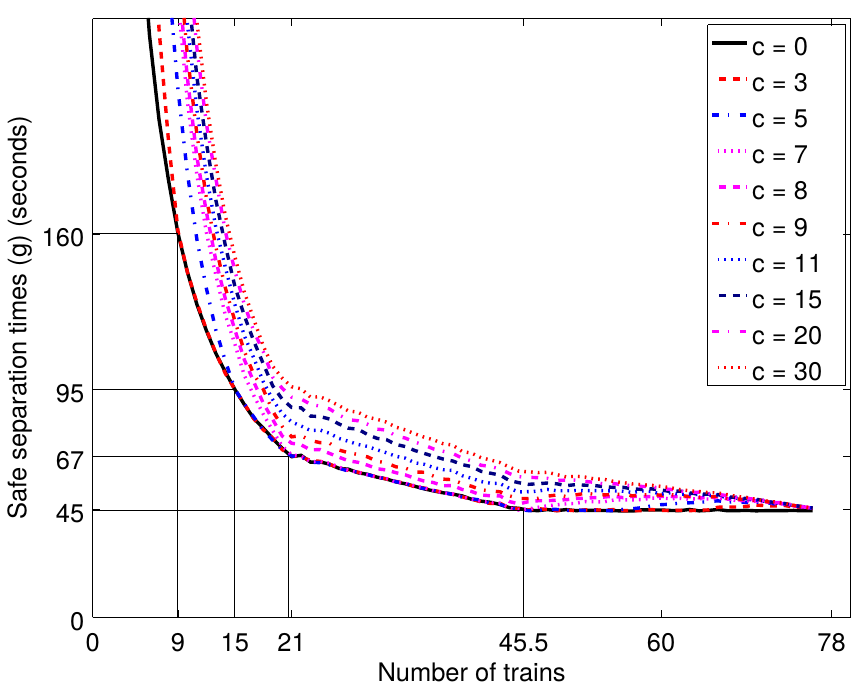} \\
  \hline
\end{tabular}
\label{tab_3}
\end{table}

\section{Conclusion and possible extensions}

Basing on the conclusions of the Max-plus linear model of chapter~\ref{chap-maxtrain}, we proposed here an extension to a stochastic
dynamic programming model, where the passenger arrivals are considered.
The model permits to take into account the effect of the passenger arrival demand on the train dwell times at platforms, and by that on the 
whole dynamics of the trains.
Numerous extensions of this approach are possible.
First, the derivation of analytic formulas for the asymptotic average train frequency, dwell time and safe separation time 
would allow a better comprehension of the traffic control model.
Second, the passenger demand being modeled here through average passenger arrival rates,
a dynamic model of the number of passengers at platforms an inside the trains would improve the traffic modeling, in particular by
taking into account the train and platform capacity limits.
Finally, the interaction between alighting and boarding passenger flows can also be modeled in order to 
take into account its effect on the train dwell times and on the train dynamics in general.


\chapter[Summary of my other contributions]{Summary of my other contributions on the stochastic traffic modeling and control} 
\label{chap-otherstoch}

This chapter summarizes some of my other contributions on the stochastic traffic modeling and control.
Two works are summarized.

First, we present a work done with Mrs Farida Manseur, a PhD student, on optimal and robust routing in road networks.
We based on an existing approach from~\cite{Sam12,Sam14,Sam14b} for optimal routing, and propose an extension to include a robustness aspect in
the routing and guidance of users of road networks. The robustness here is against link failure. Our approach permits to calculate
optimal itineraries with respect to minimum and reliable travel times, but also with respect to the robustness of the routing strategy
against link failure. Itineraries with several possible and interesting detours are more robust than others. We present the approach, with some analytical
results, and then give numerical simulation results where our approach is applied on the well known Sioux Falls network.
The main reference of this work is~\cite{Man17}.

Second, we present a macroscopic multi-lane traffic flow model which is based on a
lane assignment of the vehicular flows of road networks. 
The model assumes drivers associate respective specific utilities (traffic speed) for every lane.
The drivers then choose the lane with the highest utility.
The model is expressed by a system of conservation laws with a smooth but implicitly
defined flux function. First we explore on two data-sets how traffic data
supports the fact that traffic speed constitutes an explanatory variable of lane assignment.
Second, we address the problem of discretization of the model. Finally, we give some directions
for future research. The main reference of this work is~\cite{Far13}.

\newpage
\section{Robust routing in road networks}
\label{sec-robust}

The main reference of this work is~\cite{Man17}.
We summarize here a new adaptive algorithm for optimal and robust guidance for the users of road networks.
The algorithm extends an existing routing algorithm known as the Stochastic On Time Arrival (SOTA)~\cite{Sam12,Sam14,Sam14b,Nie09,Nie12,Fan05}.
The latter is one of the most appropriate algorithm taking into account the variability of travel times through the road networks.
It permits the derivation of the maximum cumulative probability distribution of the time arrival to a given destination in the network,
from which one can get the most reliable origin-destination paths under given travel time budgets.

The extension we proposed here introduces robustness against link and path failures in the routing.
The idea is to favor itineraries with possible and reliable alternative diversions, in case of link failures, 
with respect to itineraries without or with less reliable alternatives. 
We consider both static and dynamic versions of the algorithm. 
In the static version, the traffic dynamics is not taken into account. 
We first show some interesting properties of the algorithm in this case.
Then, we consider the dynamic version of the algorithm, where we interface our routing algorithm in the traffic simulator
SUMO (Simulation of Urban Mobility) in order to take into account the car-dynamics.
By means of some routing scenarios, we show the effectiveness of the proposed algorithm.

\subsection{Review of the SOTA approach}

The main references for the existing SOTA approach are~\cite{Fan05,FN06a,FN06b,Nie09,Nie12,Sam12,Sam14,Sam14b}.

A road network is represented by a graph $G\left(N,A\right)$, where $N$ is the set of nodes, with
$\left|N\right|=n$, and $A$ is the set of arcs, with
$\left|A\right|=a$. The set of successor and predecessor nodes of a given node $i$ are denoted by 
$\Gamma^{+1}\left(i\right)=\left\{j, \left(i,j\right){\in}A \right\}$ and 
$\Gamma^{-1}\left(i\right)=\left\{k, \left(k,i\right){\in}A \right\}$ respectively.
Travel times through the links $(i,j)$ of the network are stochastic with associated probability distribution functions $p_{ij}$.
SOTA formulation of the routing problem in the network consists in finding the best routing
strategy from any starting node $i,\left(i=1,2,{\dots}{\dots},n\right)$,
that maximizes the probability of arriving to a given destination node, denoted
$d$, within a time budget $t$.

Given a node $i{\in}N$ and a time budget $t$, 
$u_i\left(t\right)$ denotes the maximum probability for a user to arrive to destination node $d$ within a time $t$, parting from node $i$.
The maximum probabilities $u_i\left(t\right)$ are written as follows.

\begin{align}
  & u_i \left( t \right) = \max_{j \in \Gamma^{+1}(i)} \int_0^t p_{\mathit{ij}}\left(w\right)u_j\left(t-w\right)\mathit{dw}, \forall i \in N\setminus\{d\} ,j \in \Gamma^{+1}\left(i\right),0{\leq}t{\leq}T  \label{sota1}\\
  & u_d\left(t\right)=1, 0{\leq}t{\leq}T \label{sota2}\\
  & s_i\left(t\right)=\mathit{arg} \max_{j \in \Gamma^{+1}(i)} \int_0^t p_{\mathit{ij}}\left(w\right)u_j\left(t-w\right)\mathit{dw}, \forall i \in N\setminus \{d\}, j \in \Gamma^{+1}\left(i\right),0{\leq}t{\leq}T \label{sota3}
\end{align}
where $T$ is the maximum time budget, and where $s_i\left(t\right)$ is the optimal successor node for the traveler being at node $i$.

\subsection{Robust routing modeling}

A routing strategy is said here to be robust if it minimizes the deterioration of its maximum value
calculated before the depart at the origin, against eventual reconfigurations of the network that may
be due to accidents, works, etc.
Taking into account the fact that one or many links of the selected optimal path may fail during the travel,
the users may then be sensitive to path changing. Indeed, the users may prefer paths with efficient alternative detours,
with respect to paths without, or with less efficient detours, even with a loss in the average travel time, and/or in
its reliability. In order to take into account such behaviors, we propose a model that
includes the existence as well as the performance of detours for selected paths, in the calculus of the travel time
reliability (i.e. the probability of reaching a destination node). This new way of calculating travel time reliability
guarantees a kind of robustness of the guidance strategies. That is to say that the travel time reliability associated
to the obtained optimal guidance strategy is not likely to change, however associated adaptive paths change during the
travel.

\textit{Optimality} refers to the robustness of the strategy. The idea is to replace the maximum
operator in formula~(\ref{sota1}) by a weighted mean over a chosen number of successor nodes. Instead of calculating
$u_i\left(t\right)$ basing on the successor node giving the maximum value of
$u_i\left(t\right)$, we propose here to consider also other successor nodes of $i$, and we rather calculate
$u_i\left(t\right)$ basing on a weighted mean over a number of successor nodes of $i$. Let us consider the following notation.
\begin{equation}\label{eq-a1}
    A_{\mathit{ij}}\left(t\right)=\int_0^t p_{\mathit{ij}}\left( w \right) u_j\left(t-w\right)\mathit{dw},{\forall}i{\in}N\setminus \{d\},
                j{\in}\Gamma^{+1}\left(i\right),{\forall}0{\leq}t{\leq}T.
\end{equation}

We denote by $A_i(t)$ the vector $A_i(t) = (A_{i1}(t), A_{i2}(t), \ldots, A_{in_i}(t))$, where $n_i$ is the number of successor nodes 
of node $i$ in the graph.
We then define $n$ maps $S_i, i=1,\ldots, n$ as follows.
$$\begin{array}{llll}
     S_i: & \mathbb R^{n_i} & \to & \mathbb R^{n_i} \\
          & A_i(t) & \mapsto & S_i(A_i(t)),
  \end{array}, \forall i \in \{1,2, \ldots, n\}.$$
where $S_i(A_i(t))$ is the vector whose components are the same as those of $A_i(t)$, but sorted in a non increasing order.
$S_{ij}(A_i(t))$ denotes here the $j^{\text{th}}$ component $\left( S_i(A_i(t))\right)_j$ of vector $S_i(A_i(t))$.

We then rewrite the probability for a user to reach the destination node $d$ from node $i$ in a time budget $t$, as follows.
\begin{align}
  & u_i\left(t\right)=\sum _{p=1}^m\psi_p S_{ip}(A_i(t)) ,{\forall}i{\neq}d,0{\leq}t{\leq}T, \label{eq-sum} \\
  & u_d\left(t\right)=1,0{\leq}t{\leq}T, \label{eq-sum2}
\end{align}
where $m$ is a parameter giving the number of successor nodes taken into account in the sum of
formula~(\ref{eq-sum}), and $\psi _p$ are non increasing weighting coefficients satisfying
\begin{equation*}
\psi _p{\geq}0,{\forall}p{\in}\left\{1,2,{\dots},m\right\},\sum _{p=1}^m\psi _p=1,\text{ and } \psi _1{\geq}\psi
_2{\geq}{\dots}{\geq}\psi _m.
\end{equation*}

In the examples we give below for illustration, we simply take $m=2$.
In order that formula~(\ref{eq-sum}) will have a meaning, $\psi _p$ have to
be chosen such that $\psi_1~\geq~\psi_2~\geq \dots~\geq~\psi_m$. That is to say that
$\psi _p$ decrease as $S_{ip}(A_i(t))$ decrease with respect to $p$.
This dependence of $\psi _p$ on $A_{\mathit{ip}}\left(t\right)$ makes the model non-trivial.
Indeed, instead of taking the
maximum over $A_{\mathit{ip}}\left(t\right)$, with respect to successors
$p$ of $i$, as in formula~(\ref{sota1}), we take a weighted mean in
formula~(\ref{eq-sum}), where the weights are in the same order as the one of the quantities
$S_{ip}(A_i(t))$. Therefore, we need to first sort the quantities
$A_{\mathit{ip}}\left(t\right)$, before applying the mean operator. So the model~(\ref{eq-sum}) needs
more operations than the model~(\ref{sota1}). Finally, let us notice that if $m=1$, or if
$m>1$ and $\psi _p=0,{\forall}p{\geq}2$, then the model~(\ref{sota1})-(\ref{sota2})
coincides with the model~(\ref{eq-sum})-(\ref{eq-sum2}). 
Therefore, the model~(\ref{eq-sum})-(\ref{eq-sum2}) extends the model~(\ref{sota1})-(\ref{sota2}).

The optimal guidance strategy is determined by the sequence of
successor nodes $s_i\left(t\right)$ as follows.
\begin{equation}
  s_i\left(t\right)=\mathit{arg}\underset{j{\in}\Gamma^{+1}\left(i\right)}{\max}\left(A_{\mathit{ij}}\left(t\right)\right),i{\in}N. \label{eq-sum3} 
\end{equation}
We notice here that although formula~(\ref{eq-sum3}) resembles to formula~(\ref{sota3}), the resulted successor nodes from the two
formulas are not necessarily the same, since the maximized quantities in both formulas 
are calculated differently. 

\subsection{Static routing in Sioux Falls network}

We present here an implementation of the robust routing algorithm on the well known Sioux Falls network.
We assume that the link travel times on the network are drawn from bi-variate Gamma distribution.

\begin{figure}[htbp]
   \centering
   \includegraphics[scale=0.8]{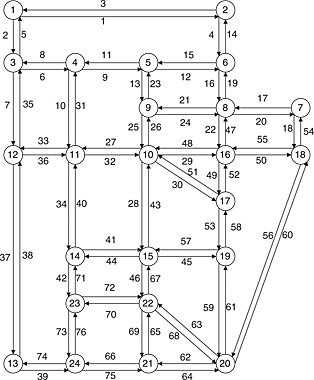}
   \caption{Sioux falls test network. The nodes as well as the links are numbered in this figure.}
   \label{sioux1}
\end{figure}   

To reach the destination node~10 parting from node~1, we have~2979 elementary paths. We apply
the robust routing model~(\ref{eq-sum})-(\ref{eq-sum2}) and derive the probabilities
$u_i(t)$ for all origin nodes $i$ of the network. 
Let us take the network of Figure~\ref{sioux1}, where we remove link 37. By that, we penalize (in term of robustness) the passage by node 12, since
we decrease the number of successor nodes of it. Therefore, node 12 will have only two, instead of three successor nodes.
At node~3, coming from node 1, we have three routing actions: go to successor node 12, which will give us only
two routing options at the next step (go to successor node 11 or back to 3), or go to successor node 4, which will give
us three routing options at the next step (go to successor nodes 5 or 11, or back to 3). Indeed, if we chose node 12 as
successor of node 3, then at node 12, if link 36 fails, we have to back to node 3. However, if we chose node 4 as
successor of node 3, then at node 4, if one of the links 9 or 10 fails, we have the possibility to change
and take a detour using the other link; see Figure~\ref{sioux1}. Therefore, passing by node 4 allows us
more options than passing by node 12.
We will see below that the propose robust routing model~(\ref{eq-sum})-(\ref{eq-sum2}) is able to take into account such robustness
criterion in the selection of the optimal routing strategy.

\begin{figure}[htbp]
   \centering
   \includegraphics[scale=0.45]{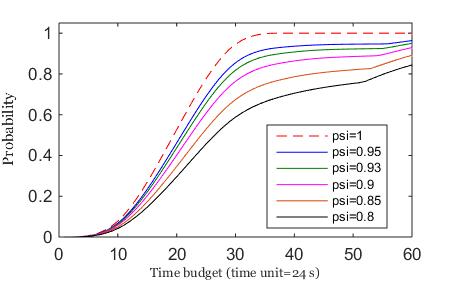}
   \includegraphics[scale=0.45]{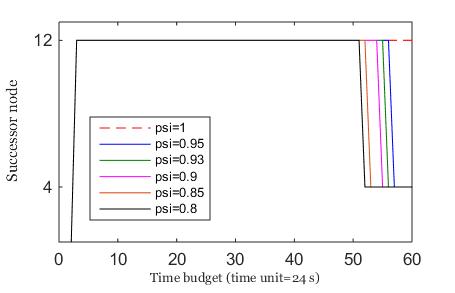}
   \caption{Left side: the robust-optimal probability of arriving on time, in function of the time budget, and for different values of $\psi $.
            Right side: the robust-optimal routing policy on node~3, in function of the time budget, and for different values of $\psi $.}
   \label{proba2}
\end{figure}   

As we take here $m=2$, we have $\psi_1$ and $\psi_2$ with $\psi_1\geq \psi_2$ and $\psi_1+\psi_2=1$.
For simplicity, we note $\psi := \psi_1 \geq 1/2$, and then $\psi_2$ is simply given by $1-\psi$.
From Figure~\ref{proba2},
the optimal policies obtained by the proposed robust routing model~(\ref{eq-sum})-(\ref{eq-sum2}) for the two cases
of $\psi = 1$ and $\psi {\in}(1/2,1)$ are clearly different. For the case
$\psi =1$, we get the same optimal policy for all considered values of time budgets.
However, with low values of $\psi $ ($\psi =0.95$, respectively 0.93, 0.9, 0.85, 0.8), where path-robustness is considered, we see that,
with a time budget greater than or equal to 57 time units (a time unit is 24 seconds here) (respectively 56, 55, 53 time units), the 
optimal successor node of node 3 is node 4 rather than node 12, even though paths passing through node 12 have lower average travel time comparing to
those passing through node~4. For example, for $\psi =0.95$, we see that with a time budget
lower than $56$ time units, the robust routing model prefers paths passing through node~12 (see blue
line on the right side of Figure~\ref{proba2}). However, with a time budget higher than~56 time units, the optimal policy on node~3 changes, and node~4
becomes the robust-optimal successor node. That means that, node~12 which has only two successor nodes is penalized, i.e. it gets low values
$u_i(t)$. Therefore, paths that pass through that node i.e. paths
with small number of alternatives or detours have low probability to be selected as robust-optimal paths.

More details on this approach are available in~\cite{Man17}, including 
\begin{itemize}
  \item Evaluation of the price of robustness, in term of travel time reliability and in term of travel time budget.
  \item An approach for fixing the parameter $\psi$ in such a way to utilize the time budget, first, to satisfy the travel time reliability,
    and then utilize the remaining time budget for robustness of the routing strategy.
  \item Robust dynamic routing with an interface with the traffic simulator SUMO.
\end{itemize}

\newpage
\section{Logit lane assigment model}

We summarize in this section a macroscopic traffic model for traffic flow assignment onto the different lanes
of a road. The model generalizes the ideas of~\cite{LK09}. 
It assumes that at a macroscopic scale, it is possible to neglect the microscopic
details of the lane-change maneuvers, and to model only the result of those maneuvers.
Moreover, the model is described as a user equilibrium, where road users choose the lanes on which to run
in function of the car-speed. The model considers a stochastic utility for the lane choice process.

Cars move in one direction on a multi-lane road, where lane change is allowed, and passing is permitted by
lane change. We assume that drivers are classified in a number of classes according to their destination.
This assumption is important because the assignment of drivers on different lanes may depend only on their
classes. For example, approaching a divergent, some drivers have to change lane in order to take the desired
destination, independent of the car density on each lane. Let us use the following notations.
\begin{itemize}
  \item $I$ the set of lanes indexed by $i$.
  \item $D$ the set of driver classes indexed $d$.
  \item $I^d$ the set of lanes accessible to user class $d$.
  \item $D_i$ the set of user classes that can move on lane $i$.
  \item $\rho^d(x,t)$ the car-density of class users $d$ at location $x$ and time $t$.
  \item $\rho_i(x,t)$ the car-density at location $x$ on lane $i$ at time $t$.
  \item $\rho^d_i(x,t)$ the car density of class users $d$ at location $x$ on lane $i$ at time $t$.
\end{itemize}

Lane utilities are user-class specific, and depend on the car-speed on every lane, which 
depends on the car-density, according to the fundamental traffic diagram of each lane.
We notice that a massive choice of a lane with a high utility has the
effect of lowering the car speed on that lane, and by that lowering its utility.
This should induce an equilibrium. We notice that the car-densities are conservative by user class,
but not by lane.

The utility of a lane $i$ for a user of class $d$ is assumed to be
\begin{equation}
  U_i^d := v_i + \theta_i^d + \xi_i^d,
\end{equation}
where $v_i = V_i(\rho_i)$ is the car-speed on lane $i$ given by the fundamental traffic diagram on the lane, 
$\theta_i^d$ is a constant expressing the preference of  class $d$ users for lane $i$,
and $\xi_i^d$ is the Gumbel random variable expressing the stochastic elements of the lane choice.
We assume that the users choose the lane with the highest utility.

The car-density is split on the different lanes according to the Logit lane assignment model
\begin{equation}\label{eq-assign}
   \rho_i^d = \rho^d \frac{\exp \left( \left\{ V_i\left( \sum_{d\in D} \rho_i^d\right)\right\} / \nu\right)}{\sum_{j\in I^d} \exp \left( \left\{ V_j\left( \sum_{d\in D} \rho_j^d\right)\right\} / \nu\right)},
      \forall d\in D, \forall i\in I^d,
\end{equation}
with $\rho^d = \sum_{j\in I^D} \rho_j^d$ and $\nu$ is a sensitivity parameter.

We then obtain a fixed point system on the variables $\rho_i^d$; see~\cite{KL12}.
System~(\ref{eq-assign}) is an assignment problem with increasing costs. 
Therefore, an equivalent formulation by an optimization problem is obtained by Beckmann transformation.
\begin{align}
   & \max_{\left(\rho_i^d\right)_{d,i}} \quad  \sum_{i\in I^d} \int_0^{\rho_i} V_i(r) dr + \sum_{d\in D, i\in I^d} \rho_i^d \mu_i^d - v \sum_{d\in D, i\in I^d} \rho^d H\left(\frac{\rho_i^d}{\rho^d}\right),  \label{eq-beck} \\
   & \quad \mid \rho_d = \sum_{i\in I^d} \rho_i^d, \forall d\in D, \label{eq-beck2}
\end{align}
where $H(x) \stackrel{\text{\tiny def}}{=} x (\ln x -1)$ is the negentropy function.

On can also show that the necessary optimality conditions of~(\ref{eq-beck})-(\ref{eq-beck2}) are equivalent to~(\ref{eq-assign}).
The partial densities $\rho_i^d$ and $\rho_i$ can then be expressed by implicit but smooth functions of the densities per destination.
Let us denote $\bar{\rho} \stackrel{\text{\tiny def}}{=} \left(\rho^d\right)_{d\in D}$. Then we write
\begin{equation}
  \rho_i^d = \mathcal R_i^d(\bar{\rho}), \quad \rho_i = \sum_{\delta\in D_i} \rho_i^{\delta} = \mathcal R_i(\bar{\rho}).
\end{equation}
The flows $q^d$ for densities by destinations $d$ are
\begin{equation}
   q^d \stackrel{\text{\tiny def}}{=} \mathcal Q^d(\bar{\rho}) = \sum_{i\in I_d} \rho_i^d V_i(\rho_i) = \sum_{i\in I_d} \mathcal R_i^d(\bar{\rho}) V_i\left( \mathcal R_i(\bar{\rho})\right).
\end{equation}
We then define the flow function $\bar{\mathcal Q} \stackrel{\text{\tiny def}}{=} \left( \mathcal Q^d\right)_{d\in D}$.

The Logit lane assignment model is then written as follows.
\begin{equation}\label{eq-logit}
   \partial_t \bar{\rho} + \partial_x \bar{\mathcal Q}(\bar{\rho}) = 0.
\end{equation}
It is not easy to solve the system~(\ref{eq-logit}) analytically, for different reasons.
\begin{itemize}
   \item The flow functions are defined implicitly.
   \item Explicit estimation of the gradient of the flow function is impossible.
   \item Analytic estimation of the eigenvalues of the flow function is impossible.
\end{itemize}

\subsection{Testing the lane choice model on data}

First tests we have done are to check whether the lane speed is indeed a relevant explanatory
variable of lane assignment and whether the Logit assignment is compatible with
measurements.

\subsubsection*{MARIUS data}

The objective here is to determine from observations in what proportions a lane choice is
explained by the speed or by the random component of the utility (i.e. by unknown
explanatory variables).
The study is conducted for various levels of density, in order to assess the importance of this
level.

\begin{figure}
  \centering
  \includegraphics[scale=0.5]{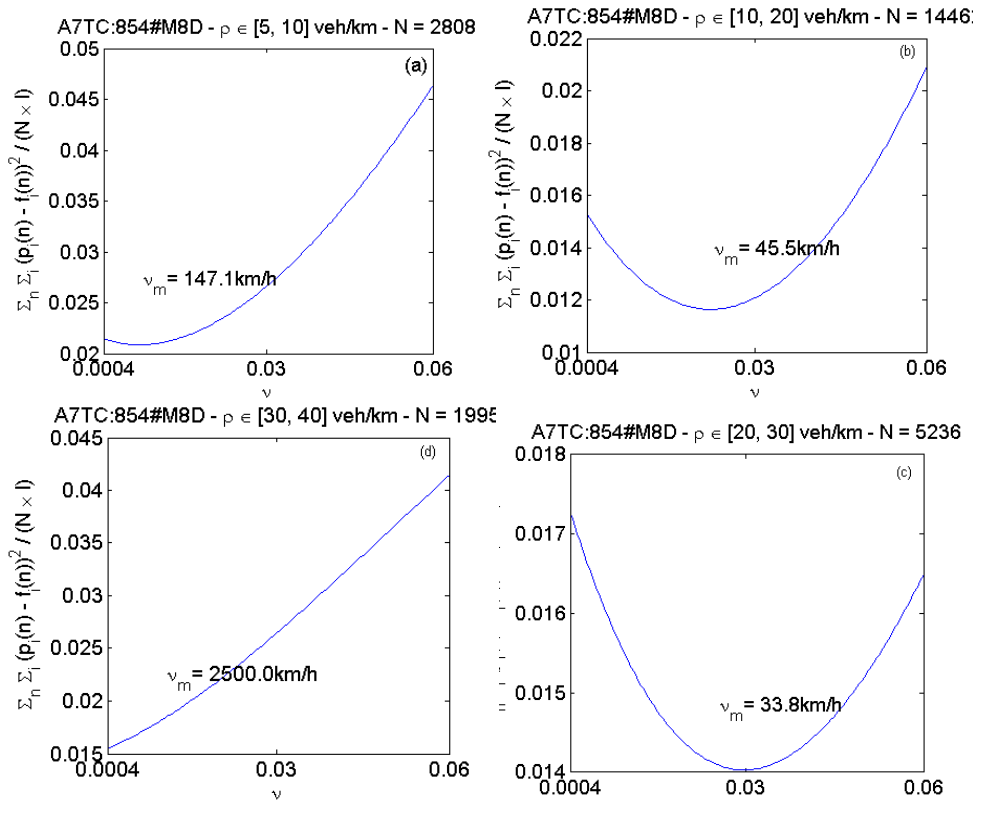}
  \caption{Estimates of $\nu$. Horizontal axis: $\nu$; vertical axis: normalized quadratic error.} 
  \label{fig-logit1}
\end{figure}

\begin{figure}
  \centering
  \includegraphics[scale=0.7]{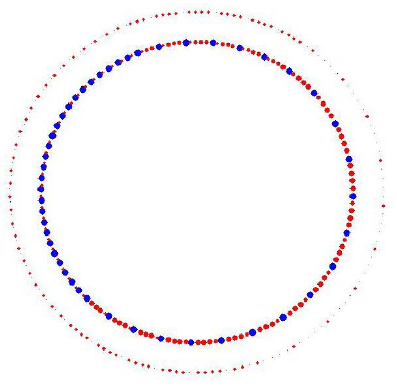}
  \caption{Position of vehicle groups on a ring road with two lanes. Interior ring: lane~1. The size of a dot corresponds to the proportion of the vehicle group on each lane.}
  \label{fig-logit2}
\end{figure} 

In the homogeneous case, where $I_d = I, \forall d$, the Logit  model  is based  on the  assumption  of the utility
$U_i = v_i + \xi_i$, where $v_i$ is the known part of the utility (the speed), and $\xi_i$ is the unknown part of
the explanatory variables. 
The parameter $\nu$ is proportional to the standard deviation of the distributions of the variables $\xi_i$.
The experimental distribution is compared to the theoretical one for values of $v$ ranging between $17$ km/h and $2500$ km/h.
The best value is the one minimizing the sum of the squares of the differences between the observed and the theoretical
probabilities; see Figure~\ref{fig-logit1}.

From these preliminary results, the speed is not an explanatory variable when the traffic is fully free-flow or congested.
However it is an explanatory variable for intermediate densities.
Moreover, for intermediate densities, the stochastic assignment model fits the data  better  than  the  deterministic  one (which  corresponds to $\nu=0$). 

Further analysis are presented in~\cite{Far13} on the Boulevard Périférique dataset.
Moreover, a number of numerical schemes are given in~\cite{Far13} (Lax-Friedrichs, Euler-Lagrange and Lagrange).
Figure~\ref{fig-logit2} shows the distribution of two classes of cars onto a ring road with two lanes, under the Lagrangian scheme 
presented in~\cite{Far13}.
We assumed here that the classe~1 users (blue color) move only on lane~1 ($\theta^1_1 = +\infty$ and $\theta^1_2 = 0$),
while users of classe~2 have the same preference to move on lanes~1 and~2 ($\theta^2_1 = \theta^2_2$).

Finally, let us notice that the model presented here admits extensions of the schemes to merges and diverges, following~\cite{LK08} for example.
Further investigations on experimental aspects are likely to improve the physical aspects of the model.
At the end, the resulting schemes should be complemented with an efficient parametric estimation method 
and should be used to model an infrastructure such as Boulevard Périphérique, or on
the NGSIM data, in order to validate the Logit lane assignment model.

\newpage
This page is intentionally left blank
 
\chapter{Conclusions and perspectives}

In conclusion, the work presented in this dissertation is partly in continuation with my previous Ph.D. and post-doctoral work,
and partly oriented towards new directions and approaches including all the novelties of the field. We believe that the study 
of transport systems and mobility in general will become even more important in the future. The mathematical modeling, 
control and numerical simulation of these systems are and will remain necessary for the comprehension,
the optimization, and the anticipation of the different phenomena and emerging behaviors of these systems. 
The \textit{dynamic systems} approach adopted here is very effective and useful for understanding the physics and dynamics of
the studied systems. However, as started in the works of this dissertation, the models and control strategies should be adapted 
with the information and communication technologies, the big data techniques of analysis, the digitalization, the automation and
even the robotics.
In parallel, new models and control strategies should rethink this area of modeling and optimal management of complex transport
systems in its new dimensions, regardless of the existing models and control strategies.

In terms of modeling scale, the microscopic scale would take more part for the modeling of the mobility thanks to the growing
equipment of the mobile units and the transport infrastructures. However, the macroscopic scale will remain indispensable for
the modeling of complex phenomena. One of the concrete examples is the macroscopic modeling of the vehicle assignment on the
lanes of a road that we presented here in chapter~\ref{chap-otherstoch}, see also~\cite{Far13}. Indeed, lane change and 
vehicle passing on a road are very difficult to model in a microscopic scale, because, inter alia, of the non-linear dynamics.
The behavior that we understand and that we are able to reproduce with the least errors is rather macroscopic (flows assigned 
to the different lanes according to the average speed on each lane, seeking an equilibrium).

In terms of the decision-making levels for traffic and mobility management, we will have more and more possibilities to consider 
the operational level thanks to the increased capacity for observing the state of the system, made possible by the availability
and variety of sensors. Thus, interesting data streams would arrive in real time, and could therefore be used for operational 
management. On the other hand, the aggregation of all this information, combined with predictive models will always be 
interesting for management at higher levels (tactical and strategic). Therefore, we will need to imagine intelligent 
management combining multiple levels of decisions, and using information available at all the levels. A concrete example is the one
of the \textit{semi-decentralized} regulation for urban road traffic, we proposed in chapter~\ref{chap-other-det},
see also~\cite{FNHL15}.  


\backmatter

\newpage
This page is intentionally left blank

\listoffigures 
\listoftables  
 
\end{document}